\documentclass[3p]{elsarticle}
\makeatletter
\def\ps@pprintTitle{%
   \let\@oddhead\@empty
   \let\@evenhead\@empty
   \let\@oddfoot\@empty
   \let\@evenfoot\@oddfoot
}
\makeatother
\usepackage{amsmath,amsfonts,amsthm,amssymb,eucal, mathtools} 
\usepackage{nicefrac} 
\DeclareMathAlphabet{\bbol}{U}{bbold}{m}{n} 
\usepackage{graphicx,subfig,tikz} 
\usepackage{multirow} 
\usepackage{hhline} 
\usepackage{color,xcolor} 
\usetikzlibrary{fadings} 
\usetikzlibrary{matrix, shapes.geometric, arrows.meta} 
\usepackage{appendix} 

\usepackage[breakable]{tcolorbox}
\newcommand{\dblock}[1]{\begin{minipage}{\widthof{a #1}}\begin{tcolorbox}[colback=white,colframe=black,right=0cm,left=0cm] #1\end{tcolorbox}\end{minipage}}

\newtheorem{thm}{Theorem}[section]
\newtheorem{prop}[thm]{Proposition}

\theoremstyle{definition}
\newtheorem{dfn}[thm]{Definition}
\newtheorem{oss}[thm]{Remark}

\usepackage[norelsize,ruled,vlined,titlenumbered]{algorithm2e}
\SetKwFor{For}{for}{}{end for}
\SetKwRepeat{Do}{do}{while}

\definecolor{celeste}{HTML}{0072BD}
\definecolor{rosso}{HTML}{A2142F}
\newcommand{\norm}[1]{\lVert#1\rVert}

\newcommand{\cJ}{\mathcal{J}}

\newcommand{\cM}{\mathcal{M}}
\newcommand{\cB}{\mathcal{B}}
\newcommand{\cL}{\mathcal{L}}

\newcommand{\cR}{\mathcal{R}}
\newcommand{\RR}{\mathbb{R}}
\newcommand{\NN}{\mathbb{N}}
\newcommand{\cF}{\mathcal{F}}

\DeclareMathOperator{\supp}{supp}
\newcommand{\NS}{N$_2$S}
\newcommand{\NNSS}{N$_2$S$_2$}
\DeclareMathOperator{\argmax}{arg max}

\newcommand{\defeq}{\coloneqq}
\newcommand{\floor}[1]{\left\lfloor#1\right\rfloor}
\newcommand{\m}[1]{\mathbf{#1}}
\newcommand{\range}[1]{1, \ldots, #1}
\DeclareMathOperator{\argmin}{argmin}
\newcommand{\scalar}[2]{\langle #1, #2 \rangle}
\def\RR{{\mathbb{R}}}

\def\p{{\bf p}}

\def\vv{{\bf v}}
\def \w {{\bf w}}

\def \y {{\bf y}}

\definecolor{dred}{rgb}{0.92,0,0}
\definecolor{dblue}{rgb}{0,0,0.7}
\definecolor{dgreen}{rgb}{0,0.5,0}
\definecolor{dgray}{rgb}{0.4,0.4,0.4} 
\definecolor{ddgray}{rgb}{0.3,0.3,0.3}
\definecolor{lgray}{rgb}{0.92,0.92,0.92}
\definecolor{gray}{rgb}{0.8,0.8,0.8}

\definecolor{LightCyan}{rgb}{0.88,1,1}

\usepackage{subfiles}
\begin{document}

\begin{frontmatter}
\title{Local spline refinement driven by fault jump estimates for scattered data approximation}

\author[unifi]{Cesare Bracco}
\ead{cesare.bracco@unifi.it}
\author[unifi]{Carlotta Giannelli}
\ead{carlotta.giannelli@unifi.it}
\author[unifi]{Francesco Patrizi\corref{cor1}}
\ead{francesco.patrizi@unifi.it}
\author[unifi]{Alessandra Sestini}
\ead{alessandra.sestini@unifi.it}
\address[unifi]{Department of Mathematics and Computer Science ``Ulisse Dini'',\\ University of Florence, Viale Giovanni Battista Morgagni 67/a, 50134 Florence, Italy}
\cortext[cor1]{Corresponding Author}

\begin{abstract}
We present new fault jump estimates to guide local refinement in surface approximation schemes with adaptive spline constructions. The proposed approach is based on the idea that, since discontinuities in the data should naturally correspond to sharp variations in the reconstructed surface, the location and size of  jumps detected in the input point cloud should drive the mesh refinement algorithm. To exploit the possibility of inserting local meshlines in one or the other coordinate direction, as suggested by the jump estimates, we propose a quasi-interpolation (QI) scheme based on locally refined B-splines (LR B-splines). Particular attention is devoted to the construction of the local operator of the LR B-spline QI scheme,  which properly adapts the spline approximation according to the nature and density of the scattered data configuration. A selection of numerical examples outlines the performance of the method on synthetic and real datasets characterized by different geographical features.
\end{abstract}

\end{frontmatter}

\section{Introduction}
Computer aided design methods offer accurate and flexible shape representations using trimmed tensor product B-spline surfaces. However, due to the inherent rigidity of tensor product B-spline constructions, achieving localized refinement requires the identification of a suitable non-standard B-spline-like basis across various unstructured mesh configurations.
Several constructions are available in
the literature, including T-splines \cite{li2012}, (truncated) hierarchical B-splines (THB-splines) \cite{thb}, and locally refined B-splines (LR B-spline) \cite{tor}. In practice, starting with an initial tensor product grid, adaptive spline refinement is usually performed by simply considering a dyadic splitting of the selected mesh elements. 
Few attempts were proposed so far to construct locally refinable extensions of tensor product constructions
which enable anisotropic local knot insertion in one or the other direction. In principle, both T-splines
and LR-splines allow the insertion of local knot segments in arbitrary regions and directions of the parametric
domain. However, since these adaptive spline models do not necessarily guarantee the linear independence
of the considered blending functions \cite{li2012,bressan1,lindep}, the focus of the research communities so far was mainly devoted
to the identification of admissible mesh configurations for a proper basis construction. In particular, different refinement
algorithms that generate LR meshes with the needed properties for LR B-spline spaces were recently proposed \cite{bressan2,N2S2,EG}. In the hierarchical spline framework instead, the linear independence of the basis construction is straightforward
\cite{vuong2011} but the required nestedness of the tensor product spline spaces, considered at different refinement
levels, precludes anisotropic refinement schemes along arbitrary directions. In order to overcome these limitations, partially nested
hierarchical B-spline refinement was proposed as Patchwork B-splines (PB-splines) \cite{PBsplines1,PBsplines2} by necessarily penalizing
the inherent simplicity of hierarchical spline constructions.

Discontinuity detection algorithms have been widely considered in different application settings to properly identify salient features in the input data and drive surface reconstruction processes. For example, fault indicators based on radial basis function interpolation were introduced in \cite{romani2019,rbf1}. Detection and recovery of discontinuity curves from scattered data were also previously addressed in \cite{bozzini2013,bozzini2014} relying on an intermediate gridded data approximation. Ordinary and gradient fault indicators based on numerical differentiation formulas were proposed in \cite{bracco2018,bracco2019} to work directly on unstructured input data without the need of additional data processing. More recently, the use of the so-called {null rules} was proposed for the design of fault detection algorithms \cite{BGC23}, by also considering its application to (adaptive) surface reconstruction with THB-splines. 
Such reconstruction employs a Quasi-Interpolation (QI) methodology, for a comprehensive overview, refer to the recent book \cite{Buh22}. A QI is composed of two ingredients essentially: a basis for a finite dimensional functional space such as, e.g., B-splines or radial functions, and linear functionals, of low computational costs, providing coefficients for the basis elements. A popular choice for the latter are the discrete Lagrange type functionals, which are based on a finite number of function evaluations. Discrete QIs have also been extensively used, for instance, in formulating  numerical quadrature rules, see for example \cite{sab, mazzia, degliesposti}. In this work as well we consider discrete Lagrange type quasi-interpolation for developing an adaptive scattered data approximation method based on Locally Refined (LR) B-splines. LR spaces have been already exploited for similar purposes in \cite{skytt2, skytt3, skytt, skytt4, gael}, in particular for surface reconstruction from terrain and bathymetry datasets. However, in these references, the approximation is carried out through the so-called multilevel B-spline methods rather than quasi-interpolation. Furthermore, the adaptivity of their schemes rely on error evaluation, which requires the computation of the approximate surface at each refinement iteration. The repeated assembling and evaluation of such intermediate approximating surfaces constitute the most expensive steps and may result highly time-consuming for large point clouds. Conversely, and similar to the approach in \cite{BGC23} for THB splines, our method separates adaptive mesh generation from  reconstruction, eliminating the need for intermediate constructions and evaluations during the refining process.
Additionally, thanks to the new jump estimates, we achieve a more comprehensive understanding of both ordinary and gradient faults. This additional information is effectively employed to enhance the compression capability of the scheme without compromising the fidelity of the produced reconstruction. Moreover, LR B-splines allow for anisotropic refinements, as local meshlines can be inserted in only one coordinate direction as opposed to THB splines. This property is exploited to further reduce degrees of freedom while sustaining the same accuracy at straight fault portions.
Particular attention is also devoted to the construction of the local operator of the LR B-spline QI scheme in such a way that it properly adapts to the nature and density of the scattered data configuration. A selection of numerical examples outlines the performance of the method on synthetic and real datasets characterized by different geographical features.

The structure of the paper is as follows. Section~\ref{sec:preliminaries} briefly recalls the definition and key properties of LR B-splines in the bivariate setting and the fault detection approach for scattered data introduced in \cite{bracco2018, faults19} and here adopted. The new jump estimates are presented in Section~\ref{sec:jstim}, while Section~\ref{sec:qi} introduces the algorithm here developed to produce the final adaptive LR quasi-interpolating spline for scattered data approximation. The numerical examples are then presented in Section~\ref{sec:exm}. Finally, Section~\ref{sec:closure} concludes the paper.


\section{Preliminaries}\label{sec:preliminaries}
Let $f:\Omega \subset \mathbb{R}^2\rightarrow \RR$ be a function on a closed finite domain $\Omega$, and $X\subset\Omega$ be a set of scattered data-sites, with associated function values $f(X)\defeq\{f(\pmb{x}) \,:\, \pmb{x} \in X\}$. The goal of the paper is to define an adaptive spline quasi-interpolant which approximates the function $f$ by using its values at the points in $X$. The main idea is that, since discontinuities of the function or of its gradient correspond to the main sharp variations of the function, the refinement of the spline space used in the approximation can be driven by the location and size of such jumps. As spline space, we choose the one spanned by LR B-splines, which, allowing local anisotropic refinement, provides the flexibility we need. In this section we recall the basics of LR B-splines and the technique we employ to find the location of the discontinuities.                                                  
\subsection{Locally Refined B-splines}\label{sec:LR}
Let $\cM$ be an open tensor mesh on $\Omega$ with respect to a bidegree $\m{p} = (p_1, p_2)$, that is, with boundary meshlines of full multiplicities, that is, with corresponding knot value repeated $p_k + 1$ times with $k = 1$ for vertical meshlines and $k = 2$ for horizontal meshlines, respectively, and with simple internal meshlines, i.e., of multiplicity one (with no repetitions of their knot values). Consider a sub-mesh in $\cM$, denoted by $\cM_B$, composed of $p_1 + 2$ vertical meshlines and $p_2 + 2$ horizontal meshlines, respectively, counting the meshline multiplicities. Such meshlines can be parametrized as $\{x_i\} \times [y_1, y_{p_2 + 2}]$ and $[x_1, x_{p_1 + 2}] \times \{y_j\}$ with $\m{x} \defeq (x_i)_{i=1}^{p_1 + 2}$ and $\m{y} \defeq (y_j)_{j=1}^{p_2 + 2}$ such that $x_i \leq x_{i + 1}$ and $y_j \leq y_{j + 1}$ for all $i = \range{p_1 + 1}$ and $j = \range{p_2 + 1}$. $\m{x}$ and $\m{y}$ are local knot vectors in the two directions from which we can defined a tensor product B-spline $B \defeq B[\m{x}, \m{y}]$. Thus, there is a one-to-one correspondence between the sub-meshes, such as $\cM_B$, that one can identify on $\cM$ and the tensor product B-splines that can possibly be defined on $\cM$. Among all of them, there is the special class of minimal support B-splines.

\begin{dfn}\label{def:MSBsplines}
We say that $B$ has minimal support on $\cM$ if no line in $\cM\setminus \cM_B$ traverses the interior of the support of $B$, $\text{int}(\supp B)$, entirely.
\end{dfn}
Therefore, what is commonly called (tensor product) B-spline set on a tensor mesh $\cM$ is the class of minimal support B-splines on $\cM$.
Given now an open tensor mesh $\cM$ with associated B-spline set $\cB$, assume to insert a new meshline $\gamma$ in $\cM$, with endpoints on $\cM$ and long enough to traverse entirely the support of one B-spline $B \in \cB$ but not necessarily traversing the entire domain $\Omega$ (in this case we would say that $\gamma$ is a local meshline). Let $\cM' \defeq \cM \cup \{\gamma\}$ be the new mesh. By construction, there is at least one B-spline $B = B[\m{x}, \m{y}] \in \cB$ that has not minimal support on $\cM'$. Suppose $\gamma$ is a vertical meshline, $\gamma = \{\widehat{x}\} \times [a, b]$ with $a \leq y_1$ and $b \geq y_{p_2 + 2}$. Then $\widehat{x}$ is not in $\m{x}$. By knot insertion of $\widehat{x}$ in $\m{x}$ we can express $B[\m{x}, \m{y}]$ in terms of two B-splines $B[\m{x}_1, \m{y}]$ and $B[\m{x}_2, \m{y}]$ with minimal support on $\cM'$. This procedure is the key idea of the Locally Refined (LR) B-splines and meshes, which are defined as follows.

\begin{figure}
\centering
\subfloat[]{
\begin{tikzpicture}[scale=3.6]
\filldraw[fill=red!50!white,draw=red, line width=4pt] (.2,.4) -- (.8,.4) -- (.8,1) -- (.2,1) -- cycle; 
\draw[red,line width=4pt] (.4,.4) -- (.4,1);
\draw[red,line width=4pt] (.6,.4) -- (.6,1);
\draw[red,line width=4pt] (.2,.6) -- (.8,.6);
\draw[red,line width=4pt] (.2,.8) -- (.8,.8);
\draw[step=.2] (0,.2) grid (1,1.2);
\draw (.5,.2) -- (.5,.8);
\draw[dashed] (.2,.7) -- (.8,.7);
\end{tikzpicture}
}
\subfloat[]{
\begin{tikzpicture}[scale=3.6]
\filldraw[fill=green!50!black,opacity=.5,draw=green!50!black, line width=4pt] (.2,.6) -- (.8,.6) -- (.8,1) -- (.2,1) -- cycle; 
\draw[green!50!black,line width=4pt,opacity=.5] (.4,.6) -- (.4,1);
\draw[green!50!black,line width=4pt,opacity=.5] (.6,.6) -- (.6,1);
\draw[green!50!black,line width=4pt,opacity=.5] (.2,.7) -- (.8,.7);
\draw[green!50!black,line width=4pt,opacity=.5] (.2,.8) -- (.8,.8);

\filldraw[fill=orange!90!black,opacity=.5,draw=orange!90!black, line width=4pt] (.2,.4) -- (.8,.4) -- (.8,.8) -- (.2,.8) -- cycle; 
\draw[orange!90!black,line width=4pt,opacity=.5] (.4,.4) -- (.4,.8);
\draw[orange!90!black,line width=4pt,opacity=.5] (.6,.4) -- (.6,.8);
\draw[orange!90!black,line width=4pt,opacity=.5] (.2,.7) -- (.8,.7);
\draw[orange!90!black,line width=4pt,opacity=.5] (.2,.6) -- (.8,.6);
\draw[step=.2] (0,.2) grid (1,1.2);
\draw (.5,.2) -- (.5,.8);
\draw (.2,.7) -- (.8,.7);
\end{tikzpicture}
}
\subfloat[]{
\begin{tikzpicture}[scale=3.6]
\filldraw[fill=green!50!black,opacity=.5,draw=green!50!black, line width=4pt] (.2,.6) -- (.8,.6) -- (.8,1) -- (.2,1) -- cycle; 
\draw[green!50!black,line width=4pt,opacity=.5] (.4,.6) -- (.4,1);
\draw[green!50!black,line width=4pt,opacity=.5] (.6,.6) -- (.6,1);
\draw[green!50!black,line width=4pt,opacity=.5] (.2,.7) -- (.8,.7);
\draw[green!50!black,line width=4pt,opacity=.5] (.2,.8) -- (.8,.8);

\filldraw[fill=cyan!90!black,opacity=.5,draw=cyan!90!black, line width=4pt] (.4,.4) -- (.8,.4) -- (.8,.8) -- (.4,.8) -- cycle; 
\draw[cyan!90!black,line width=4pt,opacity=.5] (.6,.4) -- (.6,.8);
\draw[cyan!90!black,line width=4pt,opacity=.5] (.5,.4) -- (.5,.8);
\draw[cyan!90!black,line width=4pt,opacity=.5] (.4,.7) -- (.8,.7);
\draw[cyan!90!black,line width=4pt,opacity=.5] (.4,.6) -- (.8,.6);

\filldraw[fill=magenta!90!black,opacity=.5,draw=magenta!90!black, line width=4pt] (.2,.4) -- (.6,.4) -- (.6,.8) -- (.2,.8) -- cycle; 
\draw[magenta!90!black,line width=4pt,opacity=.5] (.4,.4) -- (.4,.8);
\draw[magenta!90!black,line width=4pt,opacity=.5] (.5,.4) -- (.5,.8);
\draw[magenta!90!black,line width=4pt,opacity=.5] (.2,.7) -- (.6,.7);
\draw[magenta!90!black,line width=4pt,opacity=.5] (.2,.6) -- (.6,.6);

\draw[step=.2] (0,.2) grid (1,1.2);
\draw (.5,.2) -- (.5,.8);
\draw (.2,.7) -- (.8,.7);
\end{tikzpicture}
}
\subfloat[]{
\begin{tikzpicture}[scale=3.6]
\filldraw[fill = yellow!80!black, opacity=.5, draw = yellow!80!black, line width = 4pt] (.4,.4) -- (1,.4) -- (1,1) -- (.4,1) -- cycle; 
\draw[yellow!80!black, line width=4pt, opacity=.5] (.6,.4) -- (.6,1);
\draw[yellow!80!black, line width=4pt, opacity=.5] (.8,.4) -- (.8,1);
\draw[yellow!80!black, line width=4pt, opacity=.5] (.4,.6) -- (1,.6);
\draw[yellow!80!black, line width=4pt, opacity=.5] (.4,.8) -- (1,.8);

\filldraw[fill=cyan!90!black,opacity=.5,draw=cyan!90!black, line width=4pt] (.4,.4) -- (.8,.4) -- (.8,.8) -- (.4,.8) -- cycle; 
\draw[cyan!90!black,line width=4pt,opacity=.5] (.6,.4) -- (.6,.8);
\draw[cyan!90!black,line width=4pt,opacity=.5] (.5,.4) -- (.5,.8);
\draw[cyan!90!black,line width=4pt,opacity=.5] (.4,.7) -- (.8,.7);
\draw[cyan!90!black,line width=4pt,opacity=.5] (.4,.6) -- (.8,.6);

\draw[step=.2] (0,.2) grid (1,1.2);
\draw (.5,.2) -- (.5,.8);
\draw (.2,.7) -- (.8,.7);
\end{tikzpicture}
}
\caption{Replacement of an LR B-spline due to the refinement of the LR mesh. Consider the LR mesh $\cM$ reported in figure (a) and the bidegree $\m{p} = (2, 2)$. Let $B[\m{x},\m{y}]$ be the LR B-spline on $\cM$ whose support and tensor mesh $\cM_B$ are highlighted in figure (a). Let us insert a horizontal meshline $\gamma$ (dashed in figure (a)). Such $\gamma$ is traversing $\supp B$ and $B[\m{x},\m{y}]$ is replaced by the B-splines $B[\m{x},\m{y}_1]$ and $B[\m{x},\m{y}_2]$ involved in the knot insertion. In figure (b) we see the supports and tensor meshes of the latter on the new LR mesh. In particular we see that $B[\m{x},\m{y}_1]$ (the bottom B-spline in figure (b)) has not minimal support on the LR mesh as there is a vertical meshline traversing its support which is not in its associated tensor mesh. Thus $B[\m{x},\m{y}_1]$ is replaced as well, via knot insertion, by two other B-splines $B[\m{x}_1,\m{y}_1], B[\m{x}_2, \m{y}_1]$. Therefore, in the end, we move from $B[\m{x}, \m{y}]$, which was defined on $\cM$, to $B[\m{x}_1, \m{y}_1], B[\m{x}_2,\m{y}_1], B[\m{x}, \m{y}_2]$ on the new LR mesh $\cM' \defeq \cM \cup \{\gamma\}$. The supports and tensor meshes of the latter are represented in figure (c). In figure (d) we see that $\cM'$ has not the \NS~property as $B[\m{x}_2, \m{y}_1]$ is nested into another LR B-spline on the mesh and both are defined on simple knots, i.e., their associated tensor meshes have all distinct meshlines, of multiplicity one.}\label{fig:LRgeneration}
\end{figure}
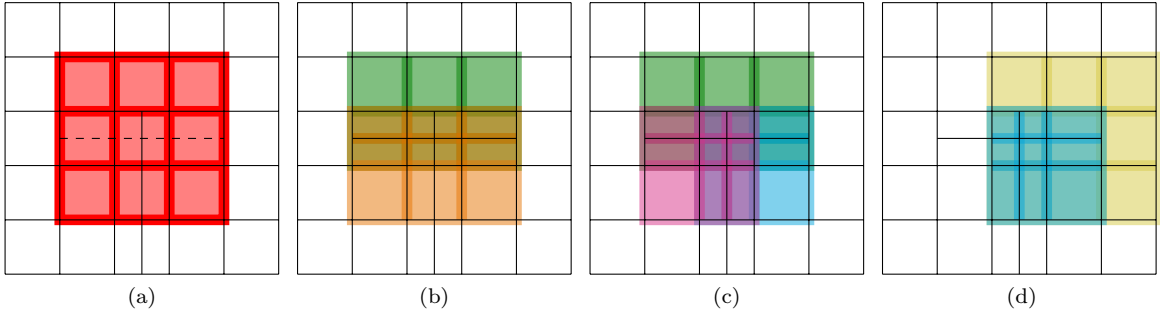

\begin{dfn}\label{def:LRBsplines}
An LR mesh $\cM'$ on $\Omega$ is either a tensor mesh or it is recursively obtained by inserting a meshline $\gamma$ in a previous LR mesh $\cM$, traversing the support of at least one of the LR B-splines on $\cM$. The LR B-spline set on $\cM'$ is either directly the B-spline set on $\cM'$ if $\cM'$ is a tensor mesh, or it is formed by replacing LR B-splines on $\cM$ that do not have minimal support on $\cM'$ with the minimal support B-splines on $\cM'$ obtained from the knot insertion procedure.
\end{dfn} 
It is important to underline that, despite a given LR mesh could be constructed with different meshline insertion orderings, the LR B-spline set is well-defined, that is, independent of which order is considered, as shown in \cite[Theorem 3.4]{tor}. Figures \ref{fig:LRgeneration} (a)--(c) show how an LR B-spline is replaced when the underlying LR mesh is refined.

Furthermore, although being defined on locally refined meshes, the LR B-splines preserve most of the properties enjoyed by B-splines on tensor meshes, such as, non-negativity, locality of the supports, being piecewise polynomials and the so-called two-scale relations with only non-negative coefficients \cite{effortless}, that is, the possibility to be expressed by LR B-splines on finer LR meshes using non-negative coefficients (which are provided by the knot insertion procedure). Moreover, they form a partition of unity when suitable scaled with positive weights \cite[Section 7]{tor}.

However, while B-splines on tensor meshes are always locally linearly independent, the almost total absence of constraints in the refinement process of LR meshes may create linear dependence relations in the LR B-spline sets \cite{lindep}. Nevertheless, LR B-splines can be as well locally linearly independent. This property has been indeed characterized as follows in \cite[Theorem 4]{bressan2}.

\begin{thm}
The following statements are equivalent:
\begin{enumerate}
\item The LR B-splines on the LR mesh $\cM$ are locally linearly independent.
\item The LR B-splines on $\cM$ form a partition of unity, without the use of scaling weights.
\item Each box $\beta$ of $\cM$ is contained in exactly $(p_1 + 1)(p_2 + 1)$ LR B-spline supports.
\item There exists no pair of LR B-splines $B_1, B_2$ on $\cM$ defined on local knot vectors made of simple knots with $B_1 \neq B_2$ and $\supp B_1 \subseteq \supp B_2$.
\end{enumerate}
\end{thm}
A box $\beta$ of $\cM$ verifying condition 3. is said non-overloaded and a mesh $\cM$ for which condition 4. is satisfied is said to have the Non-Nested-Support (\NS) property. In Figure \ref{fig:LRgeneration} (d) we see an example of LR mesh without the \NS~property.

\begin{oss}
The \NS~ property has slightly different (more complicated) definitions in \cite{EG} and \cite{bressan2}, where it was introduced. This is because we are considering a simplified framework where boundary meshlines have full multiplicity and the internal meshlines are simple in the LR meshes.
\end{oss}

Since the introduction of LR B-splines, in order to reduce the approximation error, several refinement strategies have been introduced to automatically refine the LR mesh where it is needed. The best known and used are called minimum span, full span, structured mesh, all introduced in \cite{johannessen}, hierarchical LR mesh \cite{bressan2}, effective grading \cite{EG} and Non-Nested-Support-Structured (\NNSS) mesh \cite{N2S2}. In this paper we focus on the last as it is the only strategy among them which has the \NS~ property and allows for anisotropic refinements. The local linear independence guarantees the polynomial reproduction \cite[Proposition 4.4]{N2S2}, a feature that we shall exploit in the surface reconstruction scheme presented in Section \ref{sec:qi}.

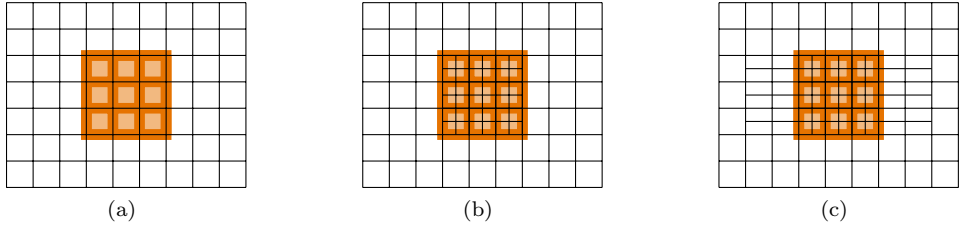
\begin{figure}
\centering
\subfloat[]{
\begin{tikzpicture}[scale = 0.35]
\fill[orange!90!black, opacity=.5] (4, 4) rectangle (7, 7);
\draw[orange!90!black, line width=4pt] (4, 4) -- (4, 7) -- (7, 7) -- (7, 4) -- cycle;
\draw[orange!90!black, step = 1, line width=4pt] (4, 4) grid (7, 7);
\draw[step = 1] (1, 2) grid (10, 9);
\end{tikzpicture}}\qquad\qquad
\subfloat[]{
\begin{tikzpicture}[scale = 0.35]
\fill[orange!90!black, opacity=.5] (4, 4) rectangle (7, 7);
\draw[orange!90!black, line width=4pt] (4, 4) -- (4, 7) -- (7, 7) -- (7, 4) -- cycle;
\draw[orange!90!black, step = 1, line width=4pt] (4, 4) grid (7, 7);
\draw[step = 1] (1, 2) grid (10, 9);
\draw[step = 0.5] (4, 4) grid (7, 7);
\end{tikzpicture}}\qquad\qquad
\subfloat[]{
\begin{tikzpicture}[scale = 0.35] 
\fill[orange!90!black, opacity=.5] (4, 4) rectangle (7, 7);
\draw[orange!90!black, line width=4pt] (4, 4) -- (4, 7) -- (7, 7) -- (7, 4) -- cycle;
\draw[orange!90!black, step = 1, line width=4pt] (4, 4) grid (7, 7);
\draw[step = 1] (1, 2) grid (10, 9);
\draw[step = 0.5] (4, 4) grid (7, 7);
\foreach \y in {0,...,2}{
\draw (2, {4.5 + \y}) -- (4, {4.5 + \y}); 
\draw (7, {4.5 + \y}) -- (9, {4.5 + \y});
};
\end{tikzpicture}}
\caption{\NNSS~refinement of a marked B-spline. Assume bidegree $\m{p} = (2, 2)$. In figure (a) we have picked a B-spline on a tensor mesh to be refined. Its support and associated tensor mesh are highlighted with colours. The first step of each iteration of the \NNSS~refinement is to halve all the non-empty knot intervals, in both directions for isotropic refinement, only in one of the two for anisotropic refinement. This step is pictured in figure (b), in the isotropic case. The next step re-instate the \NS~property on the new mesh by suitably extending the vertical or the horizontal meshlines on the new mesh, in order to prevent nesting of the LR B-spline supports. In figure (c) we show the result of this second step, when extending the horizontal meshlines. Figure (c) is also the final mesh produced by the \NNSS~refinement procedure.}\label{fig:N2S2iteration}
\end{figure}

The \NNSS~mesh strategy is a function-based strategy, i.e., the input is a set of LR B-splines that need to be refined rather than a set of mesh boxes, and it is composed of two steps. In the first step, we halve all the non-empty knot intervals of the selected LR B-splines, which simply results in a dyadic refinement of all the boxes in the tensor meshes associated to the marked LR B-splines. In case we want to refine (some of) the LR B-splines anisotropically, we would halve only the knot intervals in one of the two directions. This first step may create nesting in the LR mesh, that is, we may have a pair of LR B-splines defined on simple knots whose supports are one inside the other. The second step re-instate the \NS~ property. If $B_1$ is nested in $B_2$ we extend all the meshlines of $\cM_{B_1}$ in one of the two directions, to cross entirely $\supp B_2$. This eliminates the nesting of $B_1$ as $B_2$ is refined into other LR B-splines. However, solving the nesting of $B_1$ may trigger nesting in other LR B-splines. Anyways, it has been show \cite[Theorem 3.5]{N2S2} that this nesting propagation ends in a finite number of iterations. Figure \ref{fig:N2S2iteration} outlines the steps of the \NNSS~refinement of a B-spline.

In Section \ref{sec:qi}, we shall see how to generate an LR mesh according to the jump estimates described in Section \ref{sec:jstim}. We then present a discrete quasi-interpolation method formulated in the related LR spline space. More precisely, given the scattered data-sites $X$ and the corresponding set of values $f(X)$, we shall define a spline, denoted as $\mathfrak{Q}f$, in the space spanned by the LR B-spline set $\cL$ defined on the created LR mesh, which is a good approximation of the data in some sense, e.g., with respect to the Root-Mean-Square Error (RMSE). $\mathfrak{Q}f$ has therefore expression
\begin{equation} \label{QIdef}
\mathfrak{Q}f(x, y) \defeq \sum_{B \in \cL} q_B B(x, y)\,,
\end{equation}
for some set of coefficients $\{q_B\}_{B \in \cL}$. As the LR mesh will be finer in proximity of large variations of the scattered data (jumps), there will be more basis functions where it is needed to capture the behaviour of the surface, resembling the outcome that one would get refining the mesh based on an error indicator. \\
Each coefficient $q_B$ is computed by solving a local polynomial least squares problem with the addition of a smoothing correction term, within the support of $B$.  As we shall prove in Proposition \ref{prop:linearreproduction} at the end of Section \ref{sec:solve}, for degrees $\p \ge (1,1)$, the resulting QI shall have linear precision under the N$_2$S property of the mesh. We underline that such reproduction of linear polynomials is essential for an accurate representation of flat areas in the reconstructed surface, as it prevents the occurrence of artefacts, such as bumps. On the other hand, the smoothing term considered in our local scheme hinders the achievement of an higher order polynomial reproduction and consequently it reduces the convergence order of the method. However, it ensures the existence of the QI and moreover provides an overall more gentle and reasonable behaviour of the surface reducing also the risk of overfitting. For real world applications, these aspects should be favoured compared to the convergence rate, as often data amount and density cannot be changed so that no convergence can be achieved at all. 
\subsection{Fault Detection}\label{sec:fault}
Let $\cF^O$ and $\cF^G$ be unknown ordinary and gradient, respectively, fault curves in $\Omega$.  Setting $\cF \defeq \cF^O \cup \cF^G$, let $f: \Omega \to \RR$ be a function on $\Omega$ that is smooth (at least $C^2$) in $\Omega\setminus \mathcal{F}$ and with  finite jump discontinuities across the curves of $\cF^O$. Similarly, let $\nabla f: \Omega \setminus \cF \to \RR$ be the gradient of $f$ and assume that it has finite jump discontinuities across the curves of $\cF^G$. In this section we summarize the available strategies for the direct detection of the ordinary and gradient fault curves in $\Omega$, that is, $\cF^O$ and $\cF^G$, from the scattered evaluation of $f$. According to \cite{bracco2019}, this is done by defining two indicators based on Minimal Numerical Differentiation Formulas (MNDFs).

MNDFs, introduced in \cite{davydov2018}, generalize finite difference formulas: they approximate the value of any linear differential operator applied to multivariate functions by using their values at some scattered data points, as follows. Let $D$ be a bivariate linear differential operator of order $k$, i.e.,
\begin{equation} \label{diffop} Df(\pmb{x}) \defeq 
\sum_{\alpha \in \NN^2, |\alpha| \le k} c_\alpha(\pmb{x})\partial^\alpha f(\pmb{x}),\quad\text{with}\quad
\partial^\alpha f\defeq\frac{\partial^\alpha f}{\partial x_1^{\alpha_1}\partial x_2^{\alpha_2}}\,.  
\end{equation}
Given $\pmb{x}\in\Omega$ and a proximity set of points $X_{\pmb{x}} \defeq \{ \pmb{x}_i, i=1,\ldots,N_{\pmb{x}}\} \subset \Omega$ and the corresponding values $f(X_{\pmb{x}})$, a numerical differentiation formula is an approximation of $Df(\pmb{x})$ of the form
\begin{equation}\label{diff1}
\widehat{D}f(\pmb{x})\defeq\sum_{i=1}^{N_{\pmb{x}}} w_i f(\pmb{x}_i),
\end{equation}
where the weights $w_i$ depend on $\pmb{x}$ and are chosen such that the formula is exact for any function in a certain finite dimensional space. 
A natural choice is requiring that \eqref{diff1} is exact whenever $f$ belongs to the space $\Pi_q$ of bivariate polynomials of order less than or equal to $q$ (i.e., of total degree  $q-1$). When $q>k$ and $f$ is sufficiently smooth, it can be proved that the error of the formula goes to zero as the diameter of $X_{\pmb{x}}$ goes to zero \cite[Proposition $4$]{davydov2018}. Of course we need $N_{\pmb{x}} \ge q(q+1)/2$ in order to be able to find a formula of exactness order $q$. However, if $N_{\pmb{x}}>q(q+1)/2$ the formula is not unique, and
then it is also required that the weights minimize the weighted $\ell_2$-(semi)norm, that is,  the weight vector $\bar{\bar\w}$ should verify
\begin{equation*}\label{wnorm2}
\|\bar{\bar\w}\|_{2,\mu}=\argmin\{\|\w\|_{2,\mu}:Dp(\pmb{x})=\widehat Dp(\pmb{x}) \text{ for any }p\in\Pi_{q} \}\text{ with }
\|\w\|_{2,\mu} \defeq\left(\sum_{i=1}^{N_{\pmb{x}}}  w_i^2 \| { \pmb{x}_i}-{\pmb{x}}\|_2^{2\mu}\right)^{\frac{1}{2}}.
\end{equation*}
We denote by $\ell_2$-MNDFs the MNDFs with this constraint on the weights. In the following, for simplicity we will omit the double overline. 

For each point $\pmb{x}_i\in X$, let us consider the $\ell_2$-MNDF $\widehat\nabla f (\pmb{x}_i)$ approximating the gradient $\nabla f(\pmb{x}_i)$ with polynomial exactness of order $q=2$, and obtained setting $\mu=3$
\begin{equation} \label{indicnum} 
\widehat\nabla f (\pmb{x}_i) := \sum_{j \in J_i } \w_j f(\pmb{x}_j)\,, 
\end{equation}
using the points of the proximity set $X_i\defeq\{\pmb{x}_j:\,j\in J_i\}\subseteq X$ consisting of the $N_i\defeq 2\dim\Pi_q$ points of $X$ closest to $\pmb{x}_i$.
Note that each $ \w_j:= (w_{1,j},w_{2,j})$ is a vector of $2$ weights, one for each component of the gradient. We define the ordinary fault indicator at $\pmb{x}_i$ as 
\begin{equation} \label{indic}
I_{\nabla}(\pmb{x}_i,X_i)\defeq
\frac{\| \widehat\nabla f (\pmb{x}_i)\|_2}{ \left\|\sum_{ j \in J_i } | \w_j|\,  \| \pmb{x}_j - \pmb{x}_i \|_2 \right\|_2}\,, 
\end{equation}
where $| \w_j|\defeq (|w_{1,j}|,|w_{2,j}|)$. Similarly, the gradient fault indicator is defined as 
\begin{equation} \label{gindic}
I_{\Delta}(\pmb{x}_i,X_i)\defeq
\frac{| \widehat\Delta f (\pmb{x}_i)|}{ \sum_{ j \in J_i }  | w_j |\, \| \pmb{x}_j - \pmb{x}_i \|^2_2 }\,, 
\end{equation}
where 
\begin{equation} \label{gindicnum} 
\widehat\Delta f (\pmb{x}_i) \defeq \sum_{j \in J_i }  w_j f(\pmb{x}_j)
\end{equation}
is the  $\ell_2$-MNDF approximating $\Delta f(\pmb{x}_i)$ with polynomial exactness order $q=3$. It can be shown by analyzing its asymptotic behaviour that $I_{\nabla}(\pmb{x}_i,X_i)$ tends to take larger values for $\pmb{x}_i$ close to an ordinary fault and, similarly, $I_{\Delta}(\pmb{x}_i,X_i)$ tends to take larger values if $\pmb{x}_i$ is close to an ordinary or to a gradient fault \cite{bracco2018, bracco2019}.  

Considering these behaviours, the sets
\begin{equation*}
F^G(\alpha^G,X):=\{\pmb{x}_i\in X:\, I_{\Delta}(\pmb{x}_i,X_i) > \alpha^G\}, \qquad F^O(\alpha^O,X):=\{\pmb{x}_i\in X:\, I_{\nabla}(\pmb{x}_i,X_i) > \alpha^O\},
\end{equation*}
shall contain the points close to ordinary or gradient faults and close to only ordinary faults, respectively. The following detection algorithm was introduced in \cite{bracco2019}, using these two sets and two input parameters, denoted as $C^O$ and $C^1$.
\begin{enumerate}
\item Set $\alpha^G_1= \hbox{median}(\{I_{\Delta}(\pmb{x}_i,X_i):\, \pmb{x}_i\in X\})$,  and  compute $F^G(\alpha^G_1,X)$.
\item Set $\alpha^G_2=C^G\cdot \hbox{median}(\{I_{\Delta}(\pmb{x}_i,X_i):\, \pmb{x}_i\in F^G(\alpha^G_1,X)\})$, and  compute $F\defeq F^G(\alpha^G_2,F^G(\alpha^G_1,X))$.
\item Set $\alpha^O_1=  \hbox{median}(\{I_{\nabla}(\pmb{x}_i,X_i):\, \pmb{x}_i\in X\})$, and  compute $F^O(\alpha^O_1,X)$.
\item Set $\alpha^O_2=C^O \cdot \hbox{median}(\{I_{\nabla}(\pmb{x}_i,X_i):\, \pmb{x}_i\in F^O(\alpha^O_1)\})$, and  compute $F^O(\alpha^O_2,F^O(\alpha^O_1,X))$.
\item Return the sets:
\begin{equation}\label{class}
F,  \qquad 
F^O \defeq F\cap F^O(\alpha^O_2,F^O(\alpha^O_1,X)), \qquad 
F^G \defeq F  \setminus F^O,
\end{equation}
with points close to ordinary or gradient faults, $F$, points close to ordinary faults, $F^O$, and points close to gradient fault, $F^G$.
\end{enumerate}
Concerning the two input parameters, we assume $C^O\ge 1$ and $C^G \ge 1$, since for reasonable sampling strategies the majority of the points in $X$  are not close to faults. 

The points of $F$ are contained in strip shaped areas of certain thickness, along the true fault curves in $\Omega$. The exact location of the latter can be betted identified from $F$ if we apply a narrowing technique to each $\pmb{x} \in F$, which defines a corresponding point $\pmb{x}_N \in \Omega$ nearer to the true fault curve located inside the strip shaped area where $\pmb{x}$ lies in. 
We denote by $F_N$ the resulting narrowed set. As already done in \cite{bracco2019}, $F_N$ is obtained by using the technique based on the computation of local least squares approximations described in \cite{inkwon2000}. Essentially, for each point $\pmb{x} \in F$, we collect all the points of $X$ included in a ball of suitable radius centered at $\pmb{x}$, and compute the local linear least squares approximation for these subset of points. We then define a local system of coordinates where the $x$ axis has the direction, denoted by ${\bf d}({\pmb{x}_N})$, of the just obtained linear approximation. In this new system of coordinates, we compute a quadratic least squares approximation and then the projection of $\pmb{x}$ on it, which is $\pmb{x}_N$. Therefore, beside $\pmb{x}_N,$ such narrowing process naturally provides also an associated local direction ${\bf d}({\pmb{x}_N})$, which can be considered as an approximation of the tangent to the true fault curve at its point closest to $\pmb{x}_N$.
Note that $F_N$ is not contained in $X$ and hence we have no access to the values in $f(F_N)$ which will make the computation of the jump harder. 




\section{Jump Estimates and Classification}\label{sec:jstim}
In this section we present a new strategy to obtain directly from the scattered dataset also the information about  the  jump variation along each detected fault. This is clearly of interest on its own but it will be also exploited in Section \ref{sec:qi} to ensure a considerable reduction of the degrees of freedom in the adaptive reconstruction of $f$. We explain first the procedure to estimate the jump in ordinary faults. Namely, for each true fault curve $\cF^O,$ we aim to estimate the corresponding jump function $\cJ^O: \cF^O \rightarrow \RR$, where
\begin{equation}
\cJ^O(\widehat{\pmb{x}}) \defeq 
 \vert f_+(\widehat{\pmb{x}}) - f_-(\widehat{\pmb{x}}) \vert\,, \quad \forall\, \widehat{\pmb{x}} \in \cF^O\,,  
\label{jumpO}
\end{equation}
with $f_\pm(\widehat{\pmb{x}}) = \lim_{\pmb{x} \in \Omega^O_\pm, \pmb{x} \rightarrow \widehat {\pmb{x}}} f(\pmb{x}),$  and
$\bar \Omega^O_+ \cup \bar \Omega^O_- = \Omega, \,\, \Omega^O_+ \cap  \Omega^O_- = \varnothing\,,\,\,
\cF^O \subseteq  \bar \Omega^O_+ \cap \bar \Omega^O_-.$ Given $R > 0$, we define 
$$G_R = G_R(\widehat{\pmb{x}}) \defeq \sup\{\norm{\nabla f(\pmb{x})}_2 \,:\, \norm{\pmb{x} - \widehat{\pmb{x}}} \le R \land \pmb{x} \notin \cF^O\},$$
as first step, we prove the following proposition.
\begin{prop}
Let us assume $\widehat{\pmb{x}} \in \cF^O$ and $\pmb{x} \in X \cap \Omega^O_+$ with the segment $s(\pmb{x},\widehat{\pmb{x}})$ fully contained in $\Omega^O_+$ and let us define $\rho \defeq \norm{\pmb{x} - \widehat{\pmb{x}}}.$ Furthermore, let $R > \rho$ be a positive constant such that there exist $\pmb{x}_1 \in X \cap \Omega^O_+$ and $\pmb{x}_2 \in X \cap \Omega^O_-$ with $\rho \le r_i := \norm{\pmb{x}_i - \widehat{\pmb{x}}}_2 \le R, i=1,2.$
Then, if  $\cJ^O(  \widehat{\pmb{x}}) \ge 4R G_R$ and the segments $s(\pmb{x}_i, \widehat{\pmb{x}}), i=1,2$ are fully contained in $\Omega^O_+$ and $\Omega^O_-$, respectively, it holds
\begin{equation}\frac{ \vert f(\pmb{x}_2) - f(\pmb{x}) \vert}{r_2} \ge  \frac{ \vert f(\pmb{x}_1) - f(\pmb{x}) \vert}{r_1}\,.\end{equation}
\label{prop1}
\end{prop}  
\begin{proof}
Setting $\vv \defeq (\pmb{x} - \widehat{\pmb{x}}) / \rho,$ we can write
 $$f(\pmb{x}) = f_+(\widehat{\pmb{x}}) + \rho \left(\nabla f (\pmb{\xi}) \cdot \vv \right),$$
for some $\pmb{\xi} \in s(\pmb{x},\widehat{\pmb{x}})$, because of the assumption on $s(\pmb{x},  \widehat{\pmb{x}}).$ In a similar manner, defining $\vv_i :=  (\pmb{x}_i - \widehat{\pmb{x}}) / r_i,  i=1,2$, we have
$$f(\pmb{x}_1) = f_+(\widehat{\pmb{x}}) + r_1 \left(\nabla f (\pmb{\xi}_1) \cdot \vv_1 \right), \qquad f(\pmb{x}_2) = f_-(\widehat{\pmb{x}}) + r_2 \left(\nabla f (\pmb{\xi}_2) \cdot \vv_2 \right), $$
with $\pmb{\xi}_i \in s(\pmb{x}_i,\widehat{\pmb{x}})$.
Then, it follows that
$$\begin{array}{ll} 
\vert f(\pmb{x}_1) - f(\pmb{x}) \vert = \vert  r_1 \left(\nabla f (\pmb{\xi}_1) \cdot \vv_1 \right) -  \rho \left(\nabla f (\pmb{\xi}) \cdot \vv \right) \vert & \le G_R (r_1 + \rho)\\\\
\vert f(\pmb{x}_2) - f(\pmb{x}) \vert = \vert f_-(\widehat{\pmb{x}}) - f_+(\widehat{\pmb{x}}) + r_2 \left(\nabla f (\pmb{\xi}_2) \cdot \vv_2 \right) -  \rho \left(\nabla f (\pmb{\xi}) \cdot \vv \right) \vert & \ge \cJ^O(\widehat{\pmb{x}}) -G_R(r_2+\rho).
\end{array}
$$
Thus, we have
$$ \frac{\vert f(\pmb{x}_2) - f(\pmb{x}) \vert}{r_2} \ge \frac{\cJ^O(\widehat{\pmb{x}})}{R} - 2G_R \geq 2G_R \quad \mbox{and} \quad  \frac{ \vert f(\pmb{x}_1) - f(\pmb{x}) \vert}{r_1}\ \le 2G_R.\qedhere$$
\end{proof}

\begin{figure}
\centering
\subfloat[]{
\includegraphics[width = 0.3\textwidth]{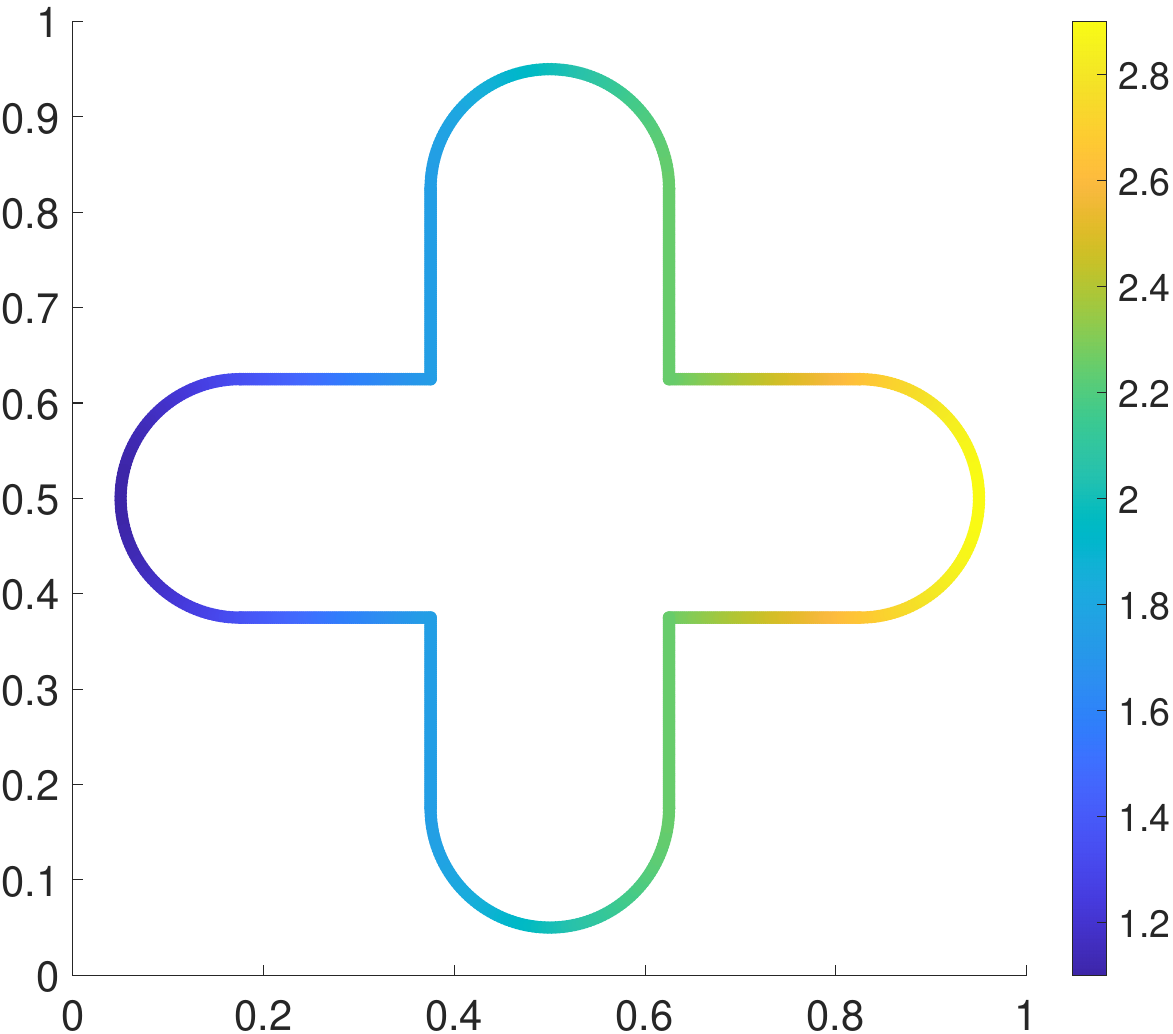}
}\quad
\subfloat[]{
\includegraphics[width = 0.3\textwidth]{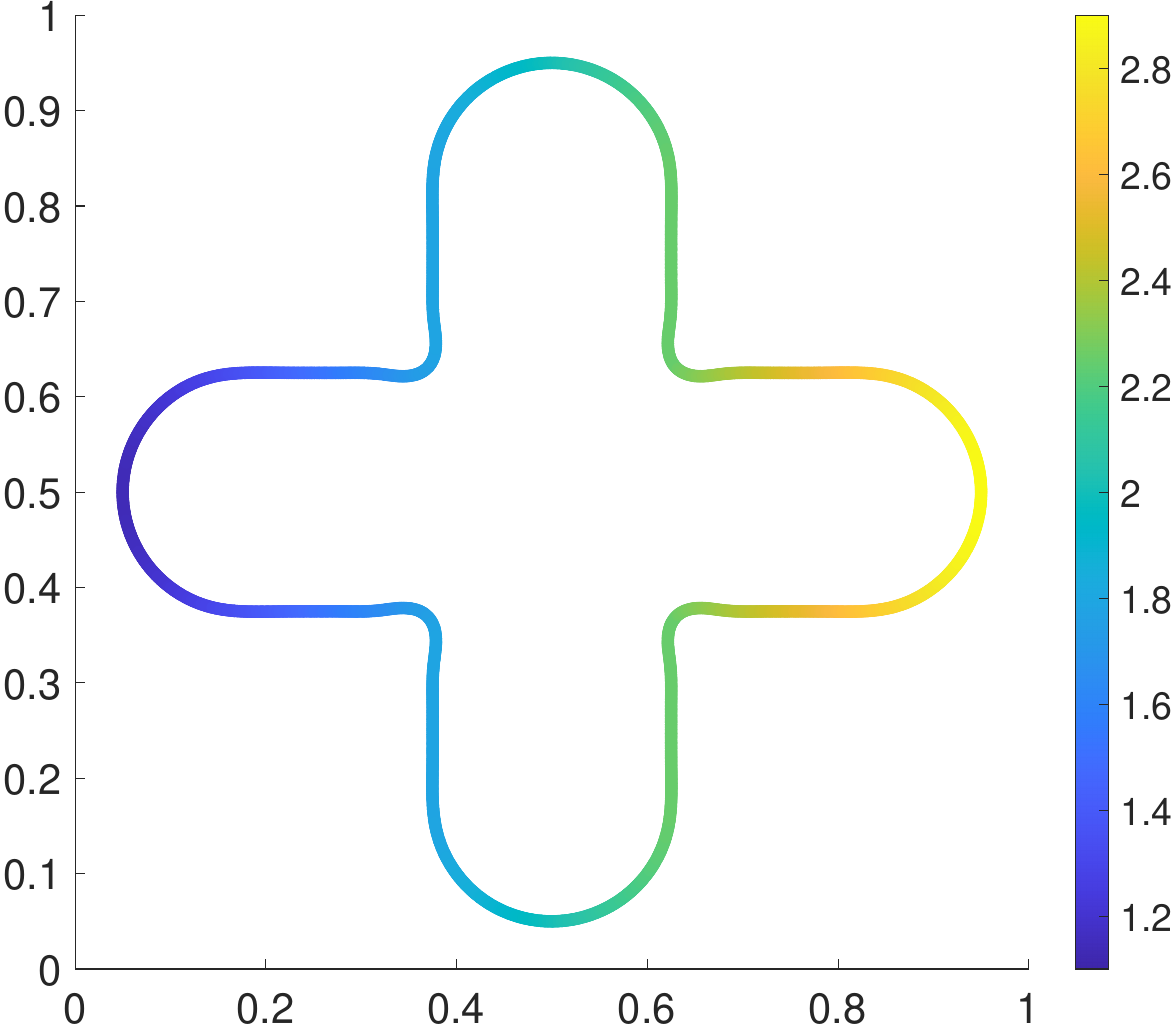}
}\quad
\subfloat[]{
\includegraphics[width = 0.3\textwidth]{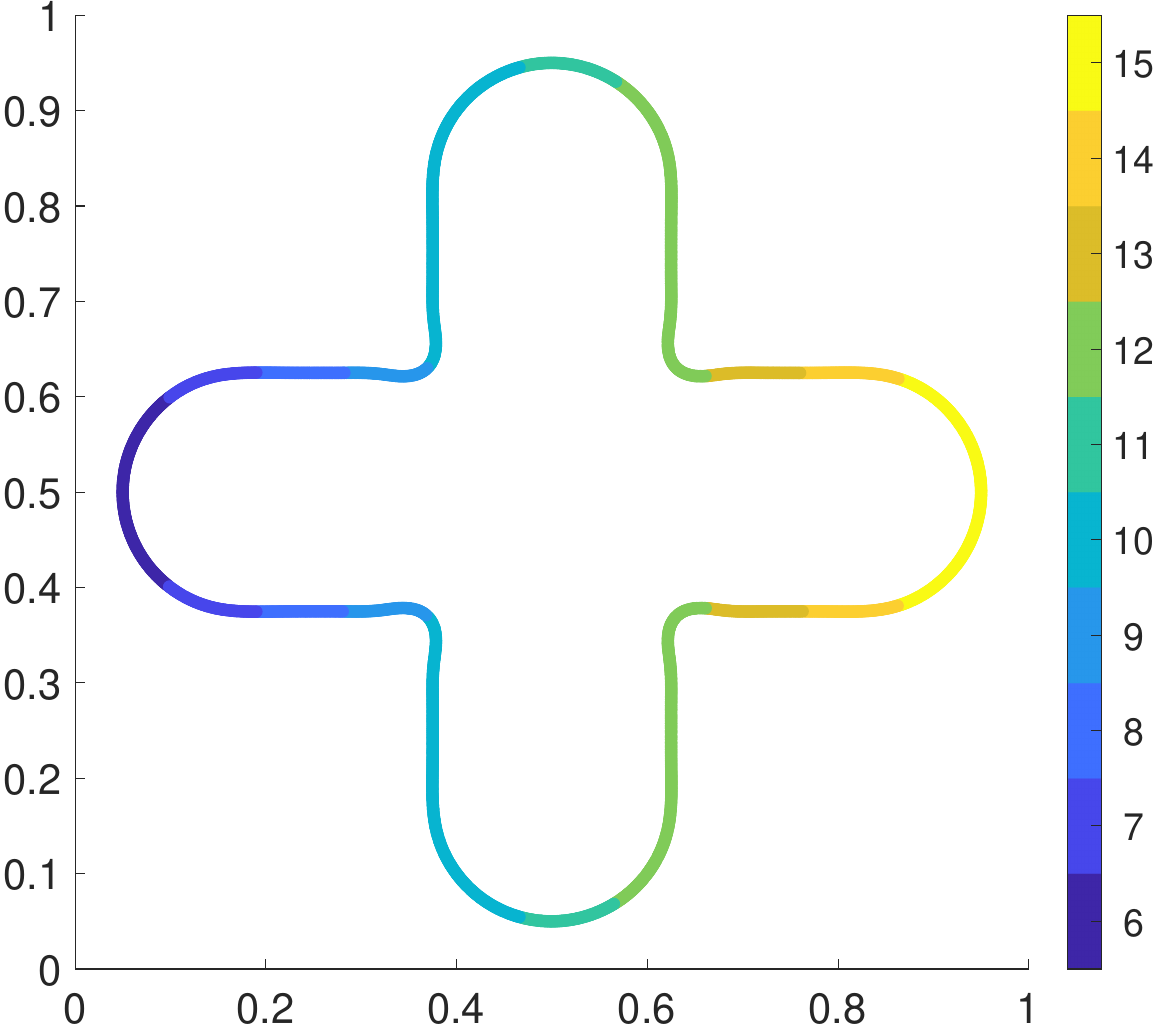}
}
\caption{Example of jump estimate and classification for an ordinary fault. We consider the sampling on one million Halton points of a piecewise surface in $[0, 1]^2$ which is zero outside of a region $R$, in the shape of a plus symbol, and it is the restriction of the plane $f(x, y) = 2x + 1$ inside $R$. In figure (a) we report the analytical jump. In figure (b) we show the fault and jump reconstructed. In figure (c) the classification of the jump when allowing $L = 15$ levels of refinements.}\label{fig:Plusjump}
\end{figure}

In principle, we do not have access to $\cF^O$. We have just available from the detection and narrowing phases, respectively, two discrete sets of points: $F^O \subset X$ and the corresponding set $F^O_N$, which is closer to $\cF^O.$ Furthermore, the value of $f$ is not assigned at the points in $F^O_N.$ Relying on the available information and the result of Proposition \ref{prop1}, we have formulated Algorithm \ref{alg:JSTIM} to assign a value $J^O(\pmb{x}_N)$ to each point $\pmb{x}_N \in F_N^O$. Since $\pmb{x}_N$ is obtained by the narrowing process \cite{inkwon2000}, and therefore represents an approximation of a point belonging to the true fault $\cF^O$, $J^O(\pmb{x}_N)$ provides an estimate of the jump $\cJ^O(\widehat{\pmb{x}})$ at the point $\widehat{\pmb{x}} \in \cF^O$  nearest to $\pmb{x}_N$. The algorithm works as follows. $J^O(\pmb{x}_N)$ is a weighted mean of differences $|f(\pmb{y}) - f(\pmb{x})|$,  with $\pmb{x} \in X$ denoting the point that generates $\pmb{x}_N$ after the narrowing precedure and $\pmb{y}$ varying in the set defined as follows. Assume, without loss of generality, that $\pmb{x}$ belongs to $\Omega_+^O$. Then, the points $\pmb{y}$ are the points of $X$ in $X_{\pmb{x}_N} \cap \Omega_-^O$ with $X_{\pmb{x}_N}$ a small annulus centered at $\pmb{x}_N$, that is, those points $\pmb{y}$ in $X$ verifying $\norm{\pmb{x} - \pmb{x}_N}_2 \eqcolon \rho \leq \norm{\pmb{y} - \pmb{x}_N}_2 \leq R$ for some small $R$. In order to select points in the annulus only on the other side of the fault, first we compute the quantities $|f(\pmb{y}) - f(\pmb{x})|/\norm{\pmb{y} - \pmb{x}_N}_2$ for all  the points $\pmb{y}$ of $X$ in the whole annulus $X_{\pmb{x}_N}$ and we sort them in a decreasing ordered vector $\m{q}$. Then, for Proposition \ref{prop1}, we know that the higher values in the entries of $\m{q}$ should correspond to points on the other side of the fault with respect to $\pmb{x}$, that is, in $\Omega_-^O$. Hence, we identify the index $k$ for which the compontent $q_{k + 1}$ is significatlly lower than the previous $q_k$. Such index is likely the index of the last point $\pmb{y}$ on the other side of the fault. The others should be on the same side of $\pmb{x}$, i.e., in $\Omega_+^O$. Thus, by doing the average of the weighted difference $|f(\pmb{y}) - f(\pmb{x})|$ only for those first $k$ entries of $\m{q}$, we make sure to consider only $\pmb{y}$ in $\Omega_-^O$. Concerning the weights, these are chosen such that they increase as the distance between $\pmb{y}$ and $\pmb{x}_N$ decreases, as closer points to $\pmb{x}_N$ should contribute more to the approximation of the jump. Namely, the weights chosen are  the normalizations of $w_{\pmb{y}} \defeq \norm{\pmb{y} - \pmb{x}_N}_2^{-1}$, in order to have all of them in $(0, 1)$ and summing to $1$.
As final remark, we stress that in principle, we should have considered $f(\pmb{x}_N)$ in the differences defining $J^O(\pmb{x}_N)$. However, we do not have access to this information as $\pmb{x}_N$ does not belong to $X$. Approximating $f(\pmb{x}_N)$ as an average is also not easy as the surronding points in $X$ could belong to either sides of the fault. On the other hand, $f(\pmb{x})$ is a good approximation of $f(\pmb{x}_N)$ for $\rho$ sufficiently small.
Figure \ref{fig:Plusjump} (b) shows the outcome of such Algorithm \ref{alg:JSTIM} in a test case which can be compared with the corresponding exact jump distribution presented in Figure  \ref{fig:Plusjump} (a).
 
\begin{algorithm}
{\footnotesize
\textbf{Input:}\begin{itemize}
\item $X, f(X)$: scattered dataset of distinct points and associated function values,
\item $\pmb{x}$: point of $X$ detected as ordinary fault, i.e., in $F^O$, 
\item $\pmb{x}_N$: point of $F_N^O$ corresponding to $\pmb{x}$, 
\item $R$: radius for a neighbourhood of $\pmb{x}_N$, with $R > \norm{\pmb{x} - \pmb{x}_N}_2 \eqqcolon \rho$. 
\end{itemize}
\textbf{Output:}\begin{itemize}
\item $J^O(\pmb{x}_N)$: estimate of the jump value $\cJ^O(\widehat{\pmb{x}})$  with $\widehat{\pmb{x}}$ being the point of $\cF^O$ closest to $\pmb{x}_N$.
\end{itemize}
\linespread{1.35}\selectfont
Let $X_{\pmb{x}_N}$ be the set of points $\pmb{y} \in X$ in the annulus $\rho \leq \norm{\pmb{y} - \pmb{x}_N}_2 \leq R$\; 
\For{$\pmb{y} \in X_{\pmb{x}_N}$}{
Compute $q(\pmb{y}) \defeq \frac{\vert f(\pmb{y}) - f(\pmb{x}) \vert}{\norm{\pmb{y} - \pmb{x}_N}_2}$\;
}
Sort the points in $X_{\pmb{x}_N}$ so that the values $q_j \defeq q(\pmb{y}_j)$ are decreasingly ordered\; 
Compute the index $k \defeq\displaystyle{ \argmax_{1\le j \le |\m{q}|-1} }(q_j - q_{j+1})$\;  
Set $J^O(\pmb{x}_N) \defeq \displaystyle\sum_{j=1}^k q_j\Big/\sum_{j=1}^k \norm{\pmb{y}_j - \pmb{x}_N}_2^{-1}$\;
\caption{Jump estimate, $J = \texttt{jump\_estimate}(X, f(X), \pmb{x},\pmb{x}_N, R)$} \label{alg:JSTIM}
}
\end{algorithm}

For gradient faults, we are interested in approximating the jump of the normal derivative of $f$ across any fault $\cF^G$, that is, the function $\cJ^G: \cF^G \rightarrow \RR$
\begin{equation}
\cJ^G(\widehat{\pmb{x}}) \defeq  
\vert \partial_{\widehat{\m{n}}} f_+(\widehat{\pmb{x}}) - \partial_{\widehat{\m{n}}} f_-(\widehat{\pmb{x}}) \vert\,, \quad \forall\, \widehat{\pmb{x}} \in \cF^G \,, 
 \label{jumpG}
\end{equation}
where, extending the notation introduced for ordinary faults, it is $\partial_{\widehat{\m{n}}} f_\pm(\widehat{\pmb{x}}) = \lim_{\pmb{x} \in \Omega^G_\pm, \pmb{x} \rightarrow \widehat{\pmb{x}}} \partial_{\widehat{\m{n}}} f(\pmb{x}),$ with $\widehat{\m{n}} = \widehat{\m{n}}(\widehat{\pmb{x}})$ denoting the unit normal to the fault at $ \widehat{\pmb{x}}$. If we had a good approximation of $\nabla f(\pmb{x})$ at all the points $\pmb{x}$ of $X$, we could just reuse Algorithm \ref{alg:JSTIM} simply by replacing the values of $f$ with the estimated values of its normal derivative. In fact, since we have a local direction $\m{d}(\pmb{x}_N)$ of the fault at every $\pmb{x}_N \in F_N^G$, we could easily estimate the jump of the normal derivative across the point $\widehat{\pmb{x}} \in \cF^G$ nearest to $\pmb{x}_N$  with respect to the direction $\m{n}_N = \m{n}_N(\pmb{x}_N) \defeq \m{d}(\pmb{x}_N)^\perp$, which constitutes an  approximation of the normal direction $\widehat{\m{n}} = \widehat{\m{n}}(\widehat{\pmb{x}})$. Unfortunately, getting good approximations of the gradient near a gradient fault is a delicate task. Indeed, at any $\pmb{x} \in X,$ the gradient estimate given by \eqref{indicnum} is not reliable if the adopted proximity set includes points on both sides of the fault. 
We then preliminary need a new way to approximate safely the gradient at least at some points of $X$ in a neighbourhood of $\pmb{x}_N$, for any $\pmb{x}_N \in F_N^G$. The procedure we propose works as follows. 

\begin{figure}
\centering
\includegraphics[width = 0.5\textwidth]{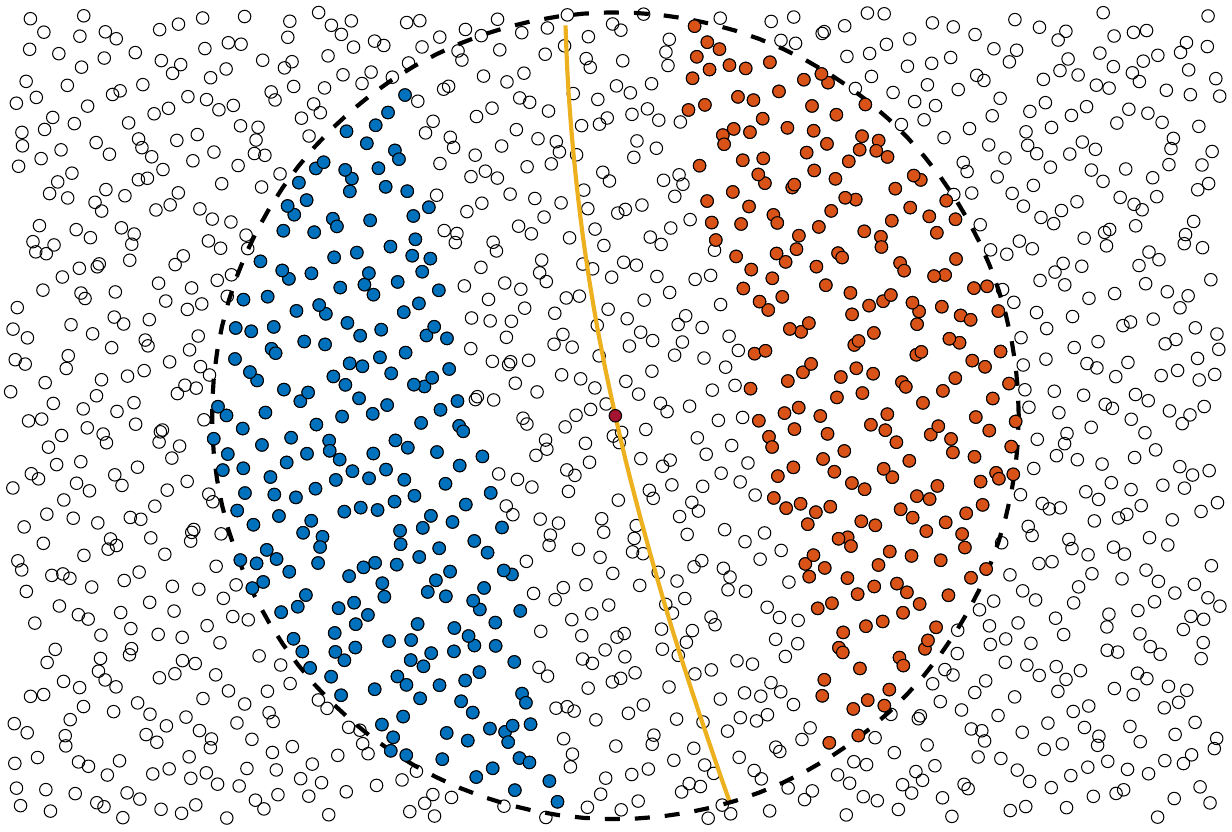}
\caption{Setting for the gradient estimate in an offset of the fault curve. All the dots belong to $X$ except the central red dot. The red central marker is the narrowed point $\pmb{x}_N \in F_N^G$ defined from a detected fault point $\pmb{x} \in F^G$. The yellow curve is the restriction of the fault curve in the considered neighbourhood (bounded by the dashed circumference). The two coloured groups of dot markers are the two sets, composing $X_{\pmb{x}_N}$, considered for the gradient estimate on the two sides of the fault curve.}\label{fig:grad4jump}
\end{figure}

For every fault point $\pmb{x}_N \in F_N^G$ and $R>0$, let us define $O(\pmb{x}_N)$ as the offset of the straight line $r_{\pmb{x}_N}(t):\pmb{x}_N + t\m{d}(\pmb{x}_N)$ of width $\nicefrac{2R}{3}$, i.e., the strip with axis $r_{\pmb{x}_N}$ and size $\nicefrac{2R}{3}$, and let $X_{\pmb{x}_N} \defeq X \cap B_R(\pmb{x}_N) \setminus O(\pmb{x}_N)$, see Figure \ref{fig:grad4jump}.
We add to $X_{\pmb{x}_N}$ the original detected point $\pmb{x} \in X$ (i.e., before narrowing replacement) if not already contained. In $X_{\pmb{x}_N}$ we have two distinct groups of points separated by the strip $O(\pmb{x}_N)$, which can be easily distinguished with geometric considerations. We choose $R$ such that there are enough points to apply the $\ell_2$-MNDF formula to estimate the gradient in each point $\pmb{y}$ of $X_{\pmb{x}_N}$ using only points of the group to which $\pmb{y}$ belongs to. Of course, for the sake of efficiency, we apply this procedure only to those points that have not a gradient value assigned yet. After applying this algorithm to all $\pmb{x}_N \in F_N^G$, a new, smaller, dataset $X^G$, and corresponding set of estimated gradient values $\widehat{\nabla}f(X^G)$, are produced, with
\begin{equation*}
X^G \defeq \bigcup_{\pmb{x}_N \in F_N^G} X_{\pmb{x}_N} \subseteq X.
\end{equation*}
Finally, we compute while running Algorithm \ref{alg:JSTIM} the normal derivatives with respect to $\m{n}_N(\pmb{x}_N)$ for every point of $X^G$ in the neighbourhood of $\pmb{x}_N \in F_N^G$. More precisely, provided the collection of local directions $\{\m{d}(\pmb{x}_N)\}_{\pmb{x}_N \in F^G}$, $X^G$ and $\widehat{\nabla}f(X^G)$ as inputs instead of $X$ and $f(X)$, we change the for cycle of Algorithm \ref{alg:JSTIM} as follows:
\begin{equation*}
\begin{minipage}{.75\textwidth}
{\footnotesize
\begin{algorithm}[H]
\linespread{1.35}\selectfont
Estimate the normal derivative at $\pmb{x}$ as $\widehat{\partial_{\m{n}_N}} f(\pmb{x}) \defeq \scalar{\widehat{\nabla}f(\pmb{x})}{\m{n}_N(\pmb{x}_N)}$\; 
\For{$\pmb{y} \in X_{\pmb{x}_N}$}{
Estimate the normal derivative at $\pmb{y}$ as $\widehat{\partial_{\m{n}_N}} f(\pmb{y}) \defeq \scalar{\widehat{\nabla}f(\pmb{y})}{\m{n}_N(\pmb{x}_N)}$\;
Compute $q(\pmb{y}) \defeq \frac{\vert \widehat{\partial_{\m{n}_N}}f(\pmb{y}) - \widehat{\partial_{\m{n}_N}}f(\pmb{x}) \vert}{\norm{\pmb{y} - \pmb{x}_N}_2}$\;
}
\end{algorithm}}
\end{minipage}
\end{equation*}
 Figure \ref{fig:esse} (b) shows the estimate of the jump in the normal derivative along a gradient fault. The analytic jump is reported in Figure \ref{fig:esse} (a).
\begin{oss}
The proposed approximation of the gradient relies on the fact that the restriction of the fault curve within the neighbourhood $B_R(\pmb{x}_N)$ is inside the stripe $O(\pmb{x}_N)$. Thereby, the two identified groups of points used for the gradient estimates are one on one side and one on the other side of the fault. 
This means that $R$ has to be small enough, especially close to points of larger curvature in the true fault, to verify this assumption. On the other hand, a smaller $R$ requires a dense enough point cloud in order to have enough points to apply the $\ell_2$-MNDF. Therefore, we assume that the point cloud is sufficiently dense, especially along the gradient faults.
\end{oss}
\begin{oss}
As already mentioned, we may need to force the inclusion of $\pmb{x}$ in $X_{\pmb{x}_N}$, if $\pmb{x}$ lies within $O(\pmb{x}_N)$. The side of the line $r_{\pmb{x}_N}$ where $\pmb{x}$ lies in decides which sub-group of points in $X_{\pmb{x}_N}$ we shall use for the $\ell_2$-MNDF. In the worst case scenario, if $\pmb{x}$ is sufficiently close to the true fault curve $\cF^G$, it would have belonged to the opposite group of points in $X_{\pmb{x}_N}$ if the sides were defined with respect to $\cF^G$ rather than $r_{\pmb{x}_N}$. Hence, the gradient estimate for $\pmb{x}$ is not reliable in these cases. However, if we use the set $X_{\pmb{x}_N} \setminus \{\pmb{x}\}$ for the $\ell_2$-MNDF, the effect is basically to extend to $\pmb{x}$ the restriction of the function to the opposite side of the gradient fault. In the interest of estimating the jump in the normal derivative using Algorithm \ref{alg:JSTIM}, we only need a coherent computation of a gradient approximation, that is, using only points on one side of the true fault. Thus, by excluding $\pmb{x}$, the value $\widehat{\nabla}f(\pmb{x})$ becomes an estimate of the gradient expression on the opposite side of the fault with respect to $\pmb{x}$. This is enough for having a reliable jump estimate from Algorithm 1, despite the original wrong grouping of $\pmb{x}$ in $X_{\pmb{x}_N}$.
\end{oss}
\subsection{Jump Classification}
\begin{figure}
\centering
\subfloat[]{
\includegraphics[width = 0.3\textwidth]{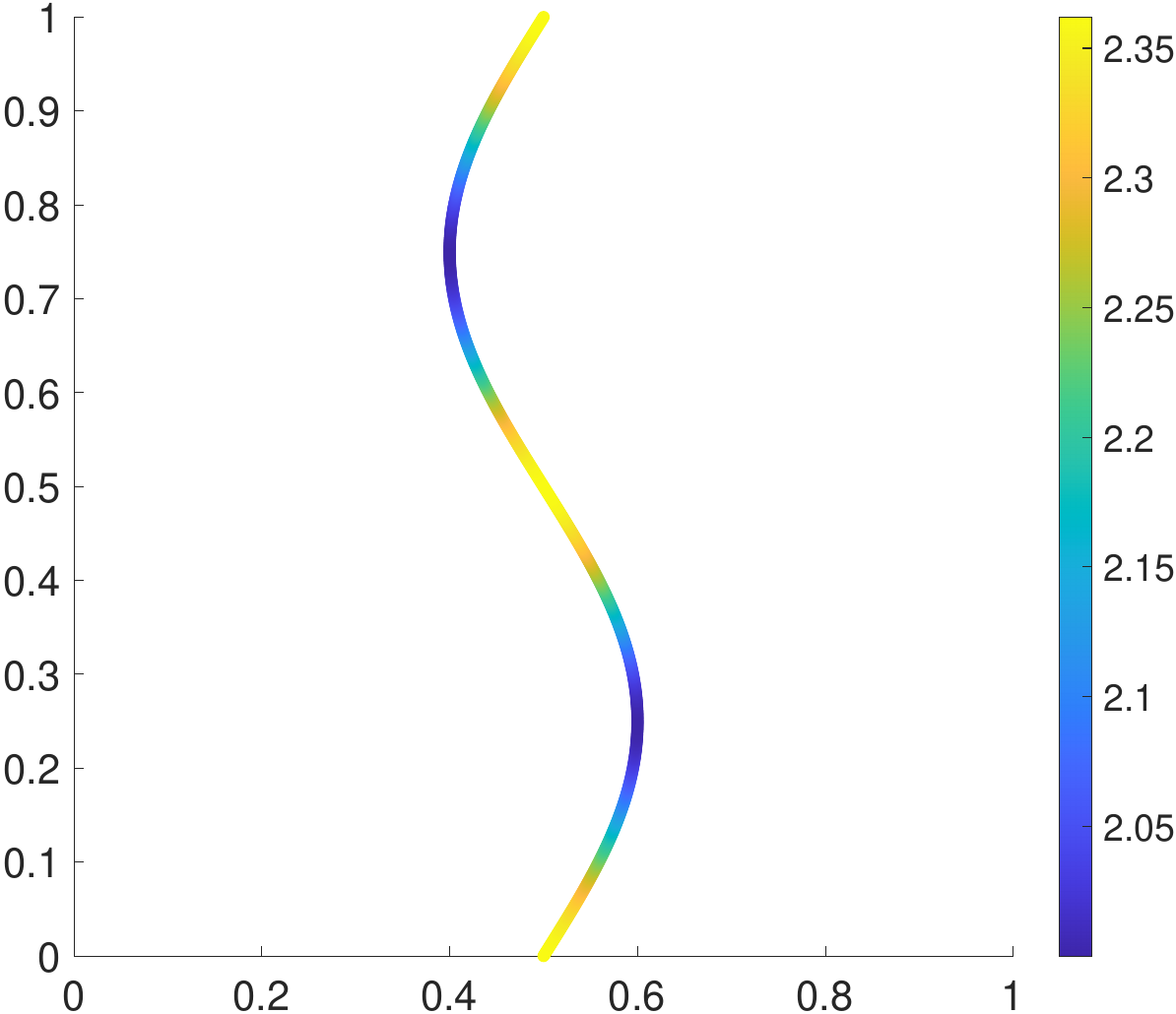}
}\quad
\subfloat[]{
\includegraphics[width = 0.3\textwidth]{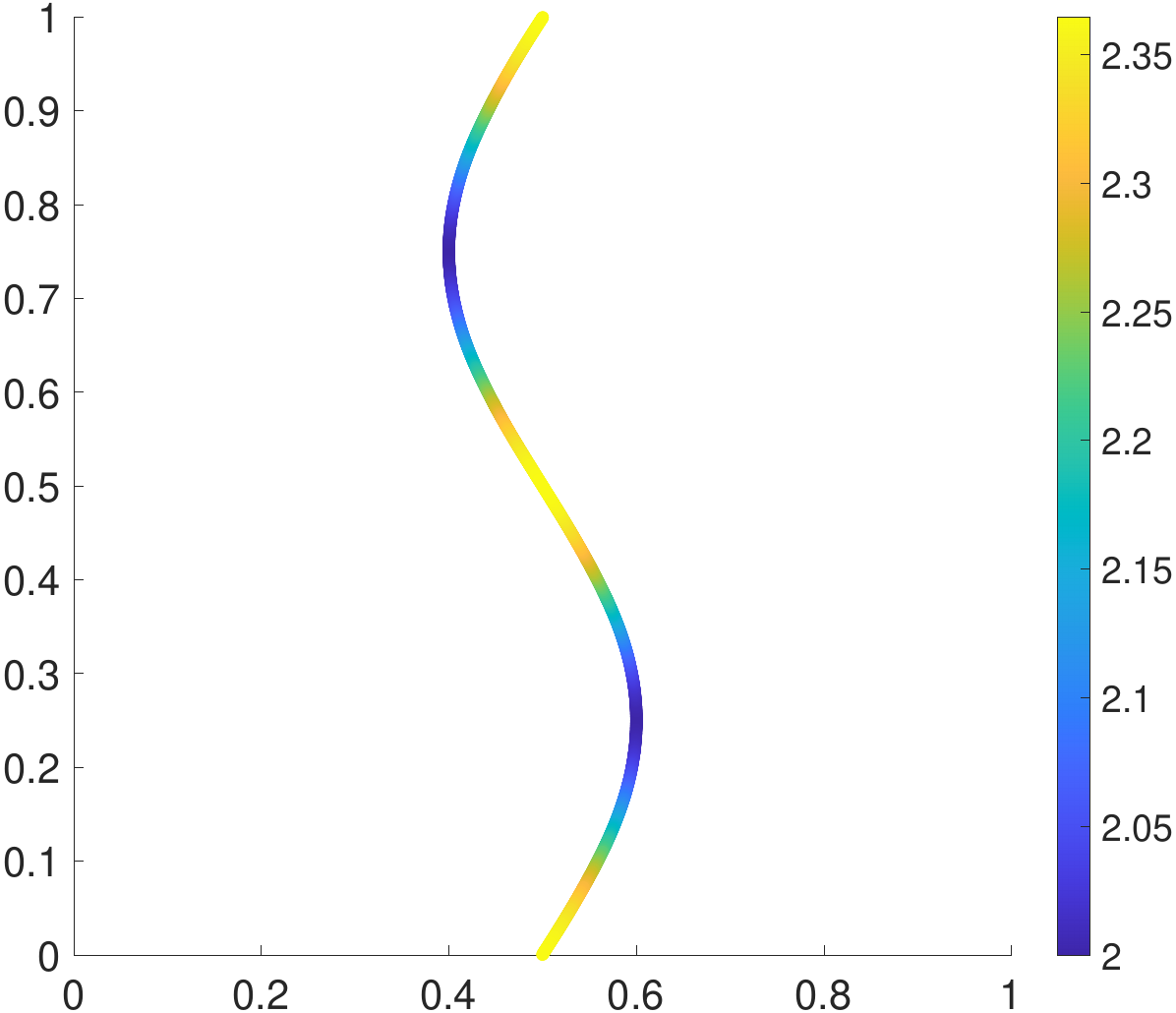}
}\quad
\subfloat[]{
\includegraphics[width = 0.3\textwidth]{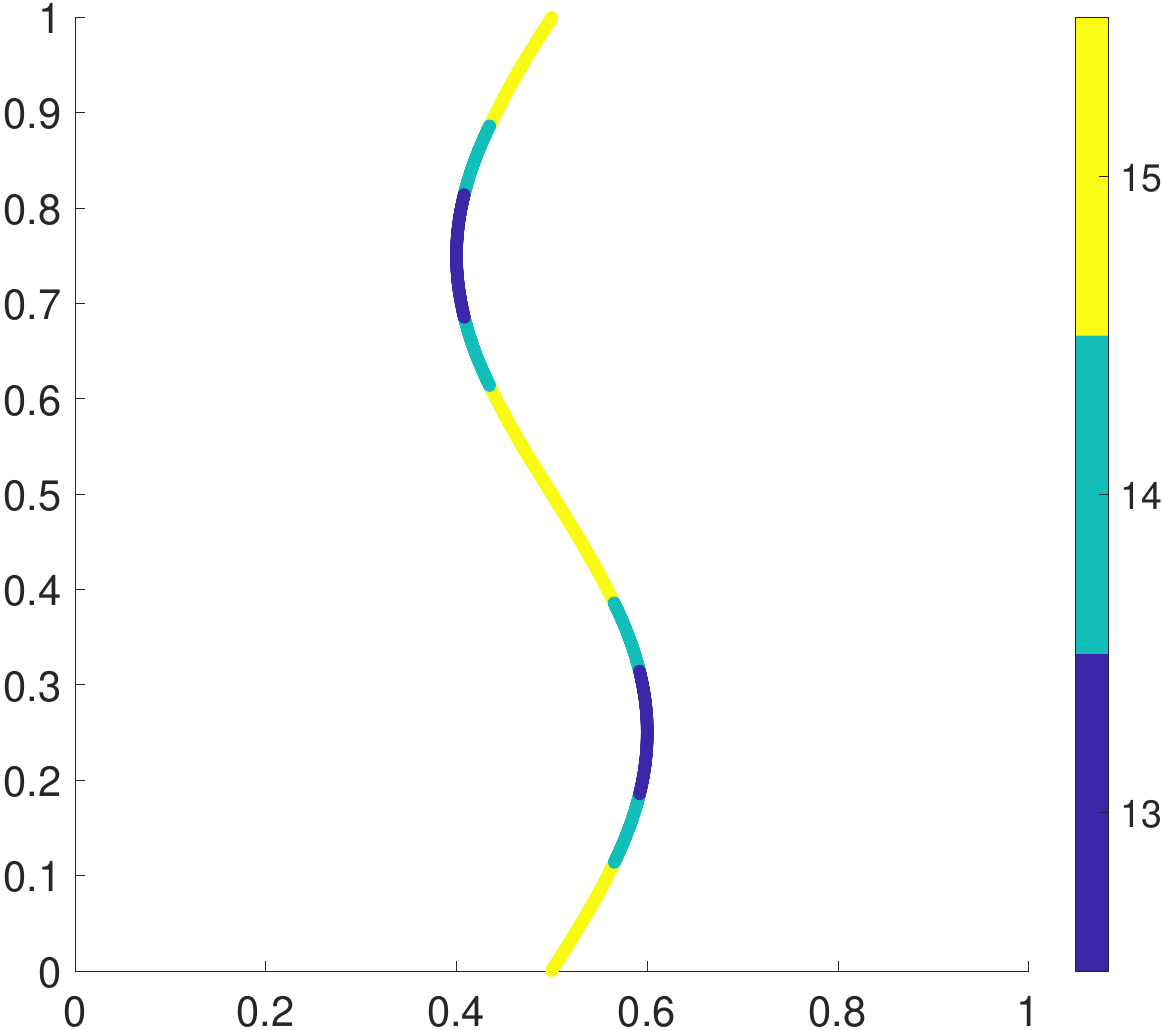}
}
\caption{Example of jump estimate and classification for a gradient fault. We consider the sampling on one million Halton points in $[0, 1]^2$ of the function $f(x, y) = |x - \frac{1}{2} - \frac{1}{10}\sin(2\pi y)|$, which is only $C^0$ along the gradient fault $x = \frac{1}{2} + \frac{1}{10}\sin(2\pi y)$. The analytical jump in the normal derivative across the fault is $g(y) = 2\sqrt{1 + \frac{\pi^2}{25}\cos^2(2\pi y)}$ and it is reported in figure (a). In figure (b) we show the reconstructed fault and jump. In figure (c) the classification of the jump when allowing $L = 15$ levels of refinements. Despite the availability of 15 classes, the detected points are gathered only in the last three by the jump classification algorithm. This is due to the fact that the range of the jump estimate is relatively small with respect to its magnitude and therefore it is reckon as fairly ``close to constant'' by the classification algorithm.}\label{fig:esse}
\end{figure}

In order to consider the computed estimates into an adaptive refinement scheme, the jump values assigned to the fault points have to be clustered into classes. Each class corresponds to a maximal level that should be applied in the refinement process. We highlight that the straightforward strategy of uniformly distributing the jump values is not an option as, in case of constant jumps, the estimates would hardly be exactly constant as well and points that should belong to the same class would rather be spread over all the available classes. Furthermore, if the range of values of the jump is not very large compared to the overall magnitude, only the last few classes should be invoked, as the fineness of the mesh should be close to the maximal everywhere along the fault. With these considerations in mind, we propose the procedure schematized in Algorithm \ref{alg:jumpclass}, which has proved to be reliable in these aspects, see Figures \ref{fig:Plusjump},\ref{fig:esse} and Figures \ref{fig:SS}--\ref{fig:FB} in the numerical tests.

\begin{algorithm}
{\footnotesize
\textbf{Input:}\begin{itemize}
\item $F_N, J$: set of narrowed points and related jump estimates,
\item $L, \bar{L}$: Maximal and minimal class indices.
\end{itemize}
\textbf{Output:}\begin{itemize}
\item $C^J$: set of jump classes of indices $\ell \in \{\bar{L}, \ldots, L\}$ related to $F_N$.
\end{itemize}
\linespread{1.35}\selectfont
Initialize $\widehat{L} \defeq \bar{L}$, $J_m \defeq \min_{\pmb{x} \in F_N} J(\pmb{x}_N)$ and $J_M \defeq \max_{\pmb{x} \in F_N} J(\pmb{x}_N)$\;
Define $\widehat{J}_{\text{old}} \defeq 0$ and 
\begin{equation}\label{eq:Jnew}\widehat{J}_{\text{new}} \defeq \frac{J_M - J_m}{J_M} \left(\frac{1}{|F_N|} \sum_{\pmb{x}_N \in F_N} J(\pmb{x}_N)\right);
\end{equation}
\While{$F_N \neq \varnothing$ and $\widehat{L} \geq \bar{L}$}{
Update $L \leftarrow \max\left\{L - \floor{\nicefrac{\widehat{J}_{\text{old}}}{\widehat{J}_{\text{new}}}}, \widehat{L}\right\}$\;
Let $D \defeq L-\widehat{L}$\;
\uIf{$D > 0$}{
Define $v_D \defeq J_M - \nicefrac{(J_M - \widehat{J}_{\text{new}})}{D}$\;}
\Else{
Define $v_D \defeq 0$;
}
Set $C^J(\pmb{x}_N) = L$ for all $\pmb{x}_N \in F_N$ with $J(\pmb{x}_N) \geq v_D$ and erase $\pmb{x}_N$ from $F_N$, $F_N \leftarrow F_N\setminus\{\pmb{x}_N\}$\;
\If{$F_N \neq \varnothing$}{
Update $\widehat{J}_{\text{old}} \leftarrow \widehat{J}_{\text{new}}$\;
Update minimal and maximal jumps, $J_m$ and $J_M$, and $\widehat{J}_{\text{new}}$ as in Equation \eqref{eq:Jnew}\;
Update $\widehat{L} \leftarrow \floor{\nicefrac{\widehat{J}_{\text{new}}}{J_M}L}$\;
}}
\If{still $F_N \neq \varnothing$}{
Define $v_0 \defeq 0$ and $v_k \defeq \widehat{J}_{\text{old}} + \frac{k - 1}{D}(J_M - \widehat{J}_{\text{old}})$ for $k = 1, \ldots, D$\; 
Set $C^J(\pmb{x}_N) \defeq \widehat{L} + k - 1$ with $\pmb{x}_N \,:\, v_{k-1} \leq J(\pmb{x}_N) < v_k$ for $k = 1, \ldots, D$\;
}}
\caption{Jump classification, $C^J = \texttt{jump\_classification}(F_N, J, L, \bar{L})$}\label{alg:jumpclass}
\end{algorithm}

The classification can be summarized as follows. Let $F_N$ be either $F_N^O$ or $F_N^G$ and $J$ be the corresponding set of estimated jump values, i.e., $J^O$ or $J^G$. We identify the points whose jump estimates should be at the maximal level $L$. Once they have been clustered, we remove them from $F_N$. We update by suitably lowering the maximal level $L$ and we determine the points, in the remaining collection, that should be in the new maximal level. We iterate this routine until nothing is left in the narrowed points collection or when a minimal level $\bar{L}$ has been reached for the classification. 
More precisely, both the current maximal level $L$ and the points whose jump values are assigned to class $L$ are identified by a quantity $\widehat{J}_{\text{new}}$, defined in Equation \eqref{eq:Jnew}, that is linear in the relative range size and in average of the estimated jump values in the remaining points of $F_N$ throughout the algorithm. In particular, the smaller the range size and/or average is, the lower will be $\widehat{J}_{\text{new}}$. As reported in Algorithm \ref{alg:jumpclass}, after each assignment of class and extraction of such just assigned points from $F_N$, $L$ is updated as $L \leftarrow \max\left\{L - \floor{\nicefrac{\widehat{J}_{\text{old}}}{\widehat{J}_{\text{new}}}}, \widehat{L}\right\}$ where $\widehat{J}_{\text{old}}$ is the $\widehat{J}_{\text{new}}$ that was computed before the reduction of $F_N$, in the previous class assignment. Whereas, the points whose assigned class will be the current $L$ are those $\pmb{x}_N \in F_N$ whose jump verifies $J(\pmb{x}_N) \geq v_D \defeq J_M - \nicefrac{(J_M - \widehat{J}_{\text{new}})}{D}$, with $D \defeq L - \widehat{L}$. $\widehat{L}$ is also a varying quantity that depends on the values of $\widehat{J}_{\text{new}}$ and the maximum $J_M$ and provides, through $D$, the range of classes that should be considered after each reduction of $F_N$. At first it is initialized as $\bar{L}$, the minimal class allowed. During the process it changes according to (i) how large the relative jump range is in the remaining collection and (ii) to how much the majority of fault points in the current $F_N$ have a jump value far from the jump maximum. If any of these two factors is small, then basically all of the remaining points should be classified at the new maximal level $L$. $\widehat{L}$ and its relation to $L$ are defined in such a way to achieve this result. 

\begin{oss}
We stress that the classification provided by Algorithm \ref{alg:jumpclass} does not necessarily involve the classes from $L$ down to $\bar{L}$. It may very well spread over less classes, depending on $L$, $\bar{L}$, the range of values of the jump and its magnitude. 
\end{oss}
\label{sec:jclass}

\section{Quasi-Interpolation of Scattered Data}\label{sec:qi}
In this section we describe a quasi-interpolation method exploiting the LR B-splines to allow adaptive refinements. The classical adaptivity cycle for numerical approximations is based on the achievement of a given tolerance by the error, computed in some norm. Assuming to allow up to $L$ refinement iterations, the structure of the procedure is the following:
\begin{equation*}
\scalebox{0.85}{
\begin{tikzpicture}
\tikzstyle{decision} = [diamond, draw, text badly centered, inner sep = -7pt, rounded corners, ultra thick]
\matrix (m)[matrix of math nodes,column sep=5em,row sep=2em]{
\text{Start} \pgfmatrixnextcell \dblock{SOLVE} \pgfmatrixnextcell \dblock{ESTIMATE} \pgfmatrixnextcell \raisebox{-0.85cm}{\begin{tikzpicture}\node[decision] {{\scriptsize $\begin{array}{c}\vspace{-.15cm}\\\ell \leq L\\ \text{and}\\ \text{err $\geq$ \texttt{tol}}\\\vspace{.05cm}\end{array}$}};\end{tikzpicture}} \pgfmatrixnextcell\text{End} \\  \pgfmatrixnextcell\pgfmatrixnextcell \dblock{REFINE}  \pgfmatrixnextcell \dblock{\begin{minipage}{\widthof{ERROR BASED}}\begin{center}ERROR BASED MARK\end{center}\end{minipage}}\\
};
\draw[ultra thick, -stealth] (m-1-1) --node[above]{$\ell = 1$} (m-1-2);
\draw[ultra thick, -stealth] (m-1-2) -- (m-1-3);
\draw[ultra thick, -stealth] (m-1-3) --node[above]{$\ell \leftarrow \ell + 1$} (m-1-4);
\draw[ultra thick, -stealth] (m-1-4) --node[right]{Yes} (m-2-4);
\draw[ultra thick, -stealth] (m-1-4) --node[above]{No} (m-1-5);
\draw[ultra thick, -stealth] (m-2-4) -- (m-2-3);
\draw[ultra thick, -stealth, rounded corners] (m-2-3) -- (-4.45, -1.4) -- (m-1-2);
\end{tikzpicture}}
\end{equation*}
More precisely, for each iteration $\ell$, a numerical approximation in the module SOLVE is computed, then the error is estimated and $\ell$ is updated. As long as $\ell$ is less than the maximal allowed loop iteration $L$ and the error is still above the input tolerance, one keeps marking those elements (or basis functions) contributing the most to the error. Consequently, these cells (or functions) are refined and then the cycle starts over by computing the new approximation.

However, the SOLVE module as well as the estimation of the error, for large size problems, are severely time-consuming. For this reason, the adaptivity cycle might be terminated prematurely, preventing it to reach the given tolerance for the error or the maximal loop iteration.

On the contrary, in the fault and jump driven strategies, the SOLVE module and error estimation are done only once at the end of the process. In fact, the scheme of the procedure is the following
\begin{equation*}
\scalebox{0.85}{
\begin{tikzpicture}
\tikzstyle{decision} = [diamond, draw, text badly centered, inner sep=3pt, rounded corners, ultra thick]
\matrix (m)[matrix of math nodes,column sep=4em,row sep=2em]{
\dblock{\begin{minipage}{\widthof{JUMP ESTIMATES AND}}\begin{center}JUMP ESTIMATES AND CLASSIFICATION\end{center}\end{minipage}}  \pgfmatrixnextcell \dblock{\begin{minipage}{\widthof{DETECTION}}\begin{center}FAULT DETECTION\end{center}\end{minipage}} \pgfmatrixnextcell \text{Start}\\
\dblock{\begin{minipage}{\widthof{FAULT/JUMP DRIVEN}}\begin{center}FAULT/JUMP DRIVEN MARK\end{center}\end{minipage}}  \pgfmatrixnextcell \dblock{REFINE} \pgfmatrixnextcell \raisebox{-.7cm}{\begin{tikzpicture}\node[decision] {$\ell \leq L$};\end{tikzpicture}} \pgfmatrixnextcell \dblock{SOLVE} \pgfmatrixnextcell \dblock{ESTIMATE}\\ \pgfmatrixnextcell \pgfmatrixnextcell \pgfmatrixnextcell \pgfmatrixnextcell \text{End}\\
};
\draw[ultra thick, -stealth] (m-1-3) -- (m-1-2);
\draw[ultra thick, -stealth] (m-1-2) -- (m-1-1);
\draw[ultra thick, -stealth] (m-1-1) --node[right]{$\ell = 1$} (m-2-1);
\draw[ultra thick, -stealth] (m-2-1) -- (m-2-2);
\draw[ultra thick, -stealth] (m-2-2) --node[above]{$\ell \leftarrow \ell + 1$} (m-2-3);
\draw[ultra thick, -stealth] (m-2-3) --node[above]{No} (m-2-4);
\draw[ultra thick, -stealth] (m-2-4) -- (m-2-5);
\draw[ultra thick, -stealth] (m-2-5) -- (m-3-5);
\draw[ultra thick, -stealth, rounded corners] (m-2-3) -- (1.625, -2) --node[above]{Yes} (-7.15, -2) -- (m-2-1);
\end{tikzpicture}}
\end{equation*}

The tremendous lightening of the approximation construction constitutes the main advantage and motivation for our approach. Note that the fault or jump driven adaptive cycle could also be used in combination with the standard error based cycle. We could in fact skip $L$ intermediate approximations and error evaluations with the former, and then check the quality of the surface constructed. Based on that, we further refine where needed, possibly by lowering the tolerance and provided a sufficiently high density of the point cloud. If required, an higher accuracy could also be attained by raising the maximum level $L$ and/or clustering differently the jump estimates in refinement levels by raising the parameter $\bar{L}$. On the contrary, if all the dataset points are within the tolerance at the end of the procedure one could consider to lower such $\bar{L}$ and $L$ to save more degrees of freedom while still achieving an acceptable approximation.

\subsection{Jump Estimates Guided Marking}
Once we have divided the jump estimates into classes, as explained in Section \ref{sec:jclass}, we have a set $C^J$, associated to the points in $F_N$, whose elements are indices $\ell$ in the range $\{1, \ldots, L\}$, with $L$ the maximum number of refinements allowed on the mesh.

Provided an LR B-spline basis $\cL$, obtained after $\bar{\ell}$ refinements of the mesh with the N$_2$S$_2$ strategy from an initial tensor mesh, we now describe the marking procedure to refine further such basis, according to the jump classification $C^J$. Basically, an LR B-spline in $\cL$ has to verify two conditions in order to be marked. First, we ask that the majority of the points of $F_N$, in the support of the LR B-spline being analized, have associated classes larger than $\bar{\ell}$. Then, we check if there are ``enough'' data points of the global dataset $X$ in its support. This second requirement is needed to assemble a local least squares system and construct the  approximation of the point cloud, as it will be explained in Section \ref{sec:solve}. It is clear that the quality of the local problem depends on the number of points of $X$ in the support. The minimal number $m \geq 3$ of data considered acceptable for this purpose is provided as input by the user. 

More precisely, the marking procedure, schematized in Algorithm \ref{alg:mark}, works as follows. For each LR B-spline $B \in \cL$, let $F_B \defeq F_N \cap \supp B\subseteq F_N$. Note that possibly $F_B$ is the empty set, in this case $B$ will not be marked for refinement. Otherwise, we split $F_B$ into groups $F_B^{\ell_1}, \ldots, F_B^{\ell_{n_B}}$ with $F_B^{\ell_k}\defeq \{ \pmb{x}_N \in F_B \,:\, C^J(\pmb{x}_N) = \ell_k\text{ for }k = 1, \ldots, n_B\}$, with $n_B$ the number of different classes to which the points in $F_B$ belong to. Then, we define the sets of those points in $F_B$ with class $\ell$ over and under $\bar{\ell}$, respectively, as
\begin{equation}\label{eq:classesoverunder}
F_B^O \defeq \bigcup_{\ell_j \geq \bar{\ell}} F_B^{\ell_j}, \qquad F_B^U \defeq \bigcup_{\ell_j < \bar{\ell}} F_B^{\ell_j}.
\end{equation}

Then, if the cardinality $|F_B^U|$ is strictly greater than $|F_B^O|$, that is, the majority of points in $F_B$ has class strictly less than $\bar{\ell}$, the LR B-spline $B$ is not refined. Otherwise, we look at the subset $X_B \defeq X \cap \supp B \subseteq X$. If $X_B$ contains at least $m$ points, then the LR B-spline $B$ is marked for refinement.

\begin{algorithm}{\footnotesize
\textbf{Input:}\begin{itemize}
\item $\cL$: LR B-spline set,
\item $X$: Dataset of scattered data,
\item $F_N, C^J$: set of narrowed points and corresponding set of jump classes,
\item $m$: Minimal number of points of $X$ to have in an LR B-spline support,
\item $\bar{\ell}$: Current refinement level.
\end{itemize}
\textbf{Output:}
\begin{itemize}
\item $\cL^M$: subset of $\cL$ of those LR B-splines marked for refinement.
\end{itemize}
\linespread{1.35}\selectfont
\For{$B \in \cL$}{
Define $F_B$ as the set of points in $F_N$ lying in $\supp B$\;
\If{$F_B = \varnothing$}{
$B$ is not marked, continue\;
}
Split the points of $F_B$ into groups, according to jump class: $F_B^{\ell_1}, \ldots, F_B^{\ell_{n_B}}$\;
Define the sets of points in $F_B$ with class over and under $\bar{\ell}$, i.e., $F_B^O$ and $F_B^U$ as in \eqref{eq:classesoverunder}\;
\If{$|F_B^O| \geq |F_B^U|$}{
Define $X_B$ as the set of points in $X$ lying in $\supp B$\;
\If{$|X_B| \geq m$}{
Mark $B$ for refinement, i.e., add it to $\cL^M$\;
}}
}}
\caption{Jump estimates driven marking, $\cL^M = \texttt{marking}(\cL, X, F_N, C^J, m, \bar{\ell})$}\label{alg:mark}
\end{algorithm}

\begin{oss}
The density of the dataset $X$ establishes the maximal number $L$ of refinement iterations that can be considered. Indeed, as the supports of the LR B-splines become smaller and smaller as we refine, the number of points in them decreases. 
Hence, if we set an $L$ too large compared to the density of the points, it will not be reached. 
\end{oss}
\begin{oss}
Although in principle the marking algorithm ensures that we refine only LR B-splines having enough data for the local least squares problem, the \NS~property demands some propagation of the refinement at each step, possibly out of the marked region. Therefore, the recovery of the local linear independence of the basis could lead to an LR B-spline set in which some LR B-splines lack of data for the approximation of the point cloud afterwards. Despite this situation is rarely experienced, as the propagation is still limited in practice, a workaround will be proposed in Section \ref{sec:solve} to deal with this inconvenience.
\end{oss}

Once we have created the subset $\cL^M \subseteq \cL$ of the marked LR B-splines, we have to decide wheater or not they should be refined anisotropically, that is, only in one of the two directions, according to the local displacement of the fault points in their supports.

As a pre-processing step, given the set $F_N$ and the associated set of local directions, we make a characteristic vector $\m{a}$ of length $|F_N|$ with values in $\{1, 2, 3\}$. Given $\pmb{x}_N \in F_N$ and a ``small'' angle $\theta$, let $\m{d}(\pmb{x}_N)$ be the related local direction. We have that 
\begin{equation}\label{eq:allinedpointvector}
a(\pmb{x}_N) \defeq\left\{ \begin{array}{ll}
1 & \text{if $\m{d}(\pmb{x}_N)$ makes an angle less than $\theta$ with $\m{e}_1 \defeq (1, 0)$},\\\\
2 & \text{if $\m{d}(\pmb{x}_N)$ makes an angle less than $\theta$ with $\m{e}_2 \defeq (0, 1)$},\\\\
3 & \text{otherwise.}
\end{array}\right.
\end{equation}
Such vector $\m{a}$ indicates which fault points have a local direction aligned, except for a small angle $\theta$, with one of the two axes and which have not. Such $\theta$ is arbitrary. In the numerical tests we performed, including those of Section \ref{sec:exm}, we have set $\theta = \nicefrac{\pi}{12}$ and we reckon that the fault points classified as axis-aligned are somehow reasonable.

\begin{oss}
The use of the angle $\theta$ is necessary for the identification of the axis-aligned fault points. In fact, even when the real fault is a line parallel to one of the axis, the fault detected could be not completely straight and the local directions associated to the fault points will not be all parallel to such line. Hence, it is necessary to compensate such deflections by allowing some room in the angles that the local directions can make with the axes.
\end{oss}

As another free parameter, we can decide from which level $L^A\in \{1, \ldots, \bar{L}, \ldots, L\}$ one should perform anisotropic refinements. Here $\bar{L}\in \{1, \ldots, L\}$ is the minimal class in the jump classification $C^J$ and so the level from which local refinements is driven by jump intensity. In fact, not necessarily we should introduce anisotropy in the refinement as soon as we start placing local insertions in the mesh where the jump classes indicate. One may require anisotropy only at the very last steps. This choice could be considered in order to capture details in the direction of the fault as well, at least for some further steps after $\bar{L}$.

With vector $\m{a}$ and level $L^A$ at hand, we decide if we should refine an LR B-spline $B \in \cL^M$ anisotropically, according to the following procedure. First of all, we should have that the refinement step we are performing is greater than $L^A$, otherwise the refinement is isotropic. If this is the case, for every marked LR B-spline $B \in \cL^M$, we consider again the subset $F_B \subseteq F_N$ as well as $\m{a}_B\subseteq \m{a}$ the restriction of $\m{a}$ to such sub-collection. Let then $F_B^k$ for $k = 1, 2, 3$ be the set of points in $F_B$ with $k$ as associated value in $\m{a}_B$. We define $\bar{k} \defeq \argmax_{k \in \{1, 2, 3\}} |F_B^k|$. If $\bar{k}$ is not unique, then it is set as $\bar{k} \defeq 3$. If $\bar{k} = 1$, $B$ will be refined only horizontally, if $\bar{k} = 2$, only vertically and if $\bar{k} = 3$, in both directions, i.e., isotropically.

\begin{oss}
A trivial observation, which however has ultimately a strong impact on the number of degrees of freedom, is that if most of the fault points makes an angle in the range $\alpha \pm \theta$ with one of the axes, we could always apply a rotation by $-\alpha$ to the point cloud, in order then to consider axis-aligned most of them and therefore perform anisotropic refinements along a larger portion of the faults.
\end{oss}

\subsection{Solving}
In this section we explain the SOLVE module, schematized in Algorithm \ref{alg:SOLVE}, which provides the approximation in the LR spline space, given as input, by means of a quasi-interpolation (QI) technique. The quasi-interpolant here described is a local method, that is, the QI coefficient assigned to each LR B-spline is computed by solving a local problem, i.e., a problem that can be set up by using only the small sub-collection of points of the dataset lying in (or close to) the support of the LR B-spline at hand.

Before starting the description of the method, we point out a couple of remarks. First, note that as opposed to the problem of approximating analytical functions, which can be evaluated at as many points as demanded to assemble the local least squares problems, when facing scattered data approximation, the least squares systems might be not always constructed, due to the possible lack of points in the dataset and/or to the bad displacement of them, e.g., in the case of almost collinear points, which leads to ill-conditioned systems \cite{thbscattered}. A special treatment is therefore required in these situations and it is detailed later in this section. As an overview and an insight of what we propose, when there is a lack of points, we enlarge to some extent the region where we look for data. If still there are not enough points, the assignment of the QI coefficient for that LR B-spline is not done by solving a least squares system, but rather by inheritance, that is, a weighted average of the QI coefficients computed for some of the LR B-splines defined on the previous, coarser, LR meshes. 
Furthermore, in order to avoid singular or nearly singular local least squares systems a smoothing correction is adopted \cite{thbscattered}. Nevertheless, the inheritance routine might still be necessary and invoked, despite the use of such smoothing correction if the system is still badly conditioned.
 
Let $\cL^A$ be the collection of all the LR B-spline sets computed so far in the adaptive cycle. Let then $\cL \in \cL^A$ be the latest LR B-spline set computed, that is, the set of LR B-splines on the current LR mesh and for which we want to compute the coefficients of the quasi-interpolant for a better approximation of the dataset. Let also $X$ be the collection of scattered points forming such dataset and $f(X)$ the corresponding set of function values. For any LR B-spline $B \in \cL$, we define $X_B$ as the subset of points in $X$ lying in $\supp B \defeq [a, b] \times [c, d]$. Let us assume, at first, that $X_B$ has at least $m \geq 3$ points, with $m$ a threshold set by the user, that is, we assume the amount of points in $X_B$ acceptable to set the local least squares system. This is the simplest case, we shall see later what to do when this condition is not verified. 
Let $\bbol{C}_B$ be the collocation matrix at the points in $X_B$ with respect to the polynomial power basis $\Pi^{\m{p}}_B$ of bi-degree $\m{p} = (p_1, p_2)$, centered in $\supp B$, that is,
\begin{equation*}
\Pi^{\m{p}}_B \defeq \left\{\left(\frac{x-a}{b-a}\right)^i\left(\frac{y - c}{d-c}\right)^j \,:\, i = 0, \ldots, p_1; j = 0, \ldots, p_2 \right\}.
\end{equation*}
Provided a so-called smoothing matrix $\bbol{S}_B$, that we shall define at the end of this section, we then assemble the least squares with smoothing matrix $\bbol{L}_B \defeq \bbol{C}_B^T\bbol{C}_B + \mu_B \bbol{S}_B$, with $\mu_B \geq 0$ a parameter which scales the smoothing effect yield by $\bbol{S}_B$ in the approximation. At this stage we have two cases, depending on the displacement of the points in $\supp B$: $\bbol{L}_B$ is ill-conditioned or not. In the former case, we will proceed as in the case of $X_B$ of cardinality less than $m$. Therefore we postpone the treatment of this instance for the time being. Let us hence assume that $\bbol{L}_B$ is not ill-conditioned. We solve the least squares problem
\begin{equation}\label{eq:localleastsquares}
\bbol{L}_B\m{u} =\bbol{C}_B^T\m{z}
\end{equation}  
with $\m{z} \defeq f(X_B)$. Let $\cB_B$ be the tensor product B-spline basis defined on the open tensor mesh set up from the knots of $B$. Note that $B \in \cB_B$. If we denote by $\bbol{R}_B$ the representation matrix of the polynomial basis $\Pi^{\m{p}}_B$ in terms of the tensor product B-splines in $\cB_B$ (see, e.g., \cite[Corollary 3.6]{splinemethods}), then we define $\pmb{\lambda} = \m{u}^T\bbol{R}_B$. If $\lambda_B$ is the coefficient of $B \in \cB_B$ and $\m{q} \defeq \{q_B \in \RR \,:\, B \in \cL\}$ is the set of QI coefficients we are looking for the LR B-spline basis $\cL$, we assign $q_B = \lambda_B$, that is, we use the same coefficient $\lambda_B$ computed as solution of the local least squares problem in terms of the tensor product B-spline basis $\cB_B$ also for the global quasi-interpolation problem in the LR B-spline basis $\cL$.

If $X_B$ has less than $m$ points, let $\widehat{\cL} \in \cL^A$ be the second last set of LR B-splines provided by the jump driven adaptivity cycle, just before $\cL$. Let then $\widehat{\cL}_B\subseteq \widehat{\cL}$ be the set of LR B-splines in $\widehat{\cL}$ whose supports contain the support of $B$. We will refer to these functions as the mothers of $B$ (one could actually create a direct graph with the LR B-splines created during the refinement procedure of a mesh as nodes and oriented edges based on the support inclusions when the mesh has the N$_2$S property, see \cite[Appendix A]{tor}). Let us define $\cR\subseteq \Omega$ as the region given by the union of the supports of the LR B-splines in $\widehat{\cL}_B$, i.e., 
\begin{equation*}
\cR \defeq \bigcup_{\widehat{B} \in \widehat{\cL}_B} \supp \widehat{B}.
\end{equation*} 
Let $X_\cR\supseteq X_B$ be the set of points in $X$ lying in $\cR$. If the cardinality of $X_\cR$ is greater or equal to $m$, we add to $X_B$ the closest $m - |X_B|$ points of $X_\cR \setminus X_B$ to the center of $\supp B$. Thereby, we ensure that $m$ points are now in $X_B$. 

However, it could happen that even by enlarging the region where we look for data to $\cR$ we still cannot find an acceptable number of points to set the local least squares system. This happens more often as we reach higher levels of fineness in the LR mesh. In fact, as the LR B-splines support become smaller and smaller, the number of data points in $\cR$ decreases. 

As possible workaround to this occurrence, we propose the following inheritance algorithm. As first step we recursively call the SOLVE module in order to have a set of coefficients $\widehat{\m{q}}_B \subseteq \widehat{\m{q}}$ for the mothers of $B$, obtained either by solving local least squares systems or by inheritance recursively, when needed. Let then $c_{\widehat{B}}$ be the coefficient of $B$ in the expression of $\widehat{B}$ in terms of a set of tensor product B-splines with $B$ among them. The QI coefficient of $B$ is assigned as the following weighted average
\begin{equation*}
q_B \defeq \sum_{\widehat{B} \in \widehat{\cL}_B} c_{\widehat{B}} \widehat{q}_{\widehat{B}}.
\end{equation*}
Note that the sum of $\{c_{\widehat{B}}\}_{\widehat{B} \in \widehat{\cL}_B}$ is the weight to assign to $B$ for the partition of unity in $\cL$, according to \cite[Section 7]{tor}. By the N$_2$S property, such weight is equal to one, that is, the $\{c_{\widehat{B}}\}_{\widehat{B} \in \widehat{\cL}_B}$ sum to one and therefore $q_B$ is a convex combination of the coefficients of the mothers of $B$. As previously mentioned, the inheritance algorithm is also invoked whenever the local least squares system \eqref{eq:localleastsquares} is (nearly) singular.

It remains to define the smoothing matrix $\bbol{S}_B$. Such matrix has been introduced in \cite{thbscattered} to ensure existence and uniqueness of a solution to the regularized least squares system, in case of non collinear points. Beyond this, the smoothing matrix yields other benefits to the approximation, such as, lowering the risk of overfitting which results in a more reasonable behaviour of the computed surface, especially in presence of outliers in the dataset. Given a pair of indices $i \in \{0, \ldots, p_1\}$ and $j \in \{0, \ldots, p_2\}$ assume to order the basis $\Pi_B^\m{p}$ by flattening the double index in some way and let $(ij) \defeq r \in \{1, \ldots, (p_1 + 1)(p_2 + 1)\}$ be the flattened index corresponding to $i,j$. If we denote by 
$$
\pi_{(ij)} \defeq \left(\frac{x - a}{b-a}\right)^i\left(\frac{y - c}{d - c}\right)^j
$$
such $(ij)$th basis function, the general $(r, s)$th entry of matrix $\bbol{S}_B$ is given by
\begin{equation}\label{eq:smoothingmatrixentries}
(\bbol{S}_B)_{r, s} \defeq \int_{\supp B} (\partial_{xx} \pi_{(i_1j_1)})(\partial_{xx} \pi_{(i_2j_2)}) +2 (\partial_{xy}\pi_{(i_1j_1)})(\partial_{xy} \pi_{(i_2j_2)}) + (\partial_{yy}\pi_{(i_1j_1)})(\partial_{yy} \pi_{(i_2j_2)})\text{d}x\text{d}y,
\end{equation}
with $r = (i_1j_1)$ and $s = (i_2j_2)$. By linearity, we can therefore write
$$
\bbol{S}_B = \bbol{S}_B^1 + \bbol{S}_B^2 + \bbol{S}_B^3
$$
where $\bbol{S}_B^k$ has $(r, s)$th entry given by the $k$th term in \eqref{eq:smoothingmatrixentries}. After a simple change of variables, we rewrite this latter sum of matrices as
\begin{equation}\label{eq:smoothingmatrix}
\bbol{S}_B = |\supp B| \left(\frac{1}{(b - a)^4} \bbol{S}^1 + \frac{2}{|\supp B|^2} \bbol{S}^2 + \frac{1}{(d - c)^4}\bbol{S}^3\right)
\end{equation}
with $\bbol{S}^k$ obtained by changing the domain of integration from $\supp B$ to the unit square $[0,1]^2$ and by replacing the basis $\Pi_B^\m{p}$, used in the expression of $\bbol{S}_B^k$, with the classical power polynomial basis $\Pi^\m{p}$, whose $(ij)$th element is
$$
\bar{\pi}_{(ij)} \defeq x^iy^j.
$$

The pseudo-code of the SOLVE module is presented in Algorithm \ref{alg:SOLVE}. Note that we give the smoothing parameter as input, that is, we set $\mu_B$ equal to $\mu$ for all the LR B-splines $B \in \cL$. Obviously, this is not the only possible choice. The smoothing effect could be dynamically changed according to the size of the supports and/or number of data points available, in order, e.g., to cope with outliers and noise.

\begin{algorithm}{\footnotesize
\vspace{8.825cm}
\hfill\begin{minipage}{0.92\textwidth}\vspace{-4cm}\begin{tcolorbox}[colback=gray!30!white, boxrule = 0pt,right=0cm,left=0cm,opacityback=1] 
\hfill \textbf{Enlargement}\vspace{1.45cm}
\end{tcolorbox}\end{minipage}\\
\hfill\begin{minipage}{0.92\textwidth}\vspace{-0.65cm}\begin{tcolorbox}[colback=gray!30!white, boxrule = 0pt,right=0cm,left=0cm,opacityback=1] 
\hfill \textbf{Inheritance}\vspace{2.2cm}
\end{tcolorbox}\end{minipage}\\
\hfill\begin{minipage}{0.92\textwidth}\vspace{2cm}\begin{tcolorbox}[colback=gray!30!white, boxrule = 0pt,right=0cm,left=0cm,opacityback=1] 
\hfill \textbf{Inheritance}\vspace{2.2cm}
\end{tcolorbox}\end{minipage}\\
\vspace{-16.125cm}
\textbf{Input:}\begin{itemize}
\item $\cL$: LR B-splines for which we want to compute the QI coefficients,
\item $\cL^A$: Archive of all the LR B-spline sets computed so far,
\item $X, f(X)$: Dataset of scattered points and corresponding set of function values,
\item $m$: Minimal number of points in $X$ considered sufficient to set a least squares system,
\item $\mu$: Smoothing parameter to scale the smoothing effect,
\item $\kappa_M$: Condition number threshold above which we consider a matrix ill-conditioned.
\end{itemize}
\textbf{Output:}\begin{itemize}
\item $\m{q}$: vector of the QI coefficients for the LR B-spline set $\cL$.
\end{itemize}
\linespread{1.35}\selectfont
\For{$B \in \cL$}{
Define $X_B$ as the set of points of $X$ lying in $\supp B$\;
\If{$|X_B| < m$}{
Extract from $\cL^A$ the set of mothers of $B$, $\widehat{\cL}_B$. Define $\displaystyle\cR \defeq \bigcup_{\widehat{B} \in \widehat{\cL}_B} \supp \widehat{B}$\;
Define $X_{\cR}\supseteq X_B$ as the set of points of $X$ contained in $\cR$\;
\uIf{$|X_{\cR}| \geq m$}{
Add to $X_B$ the closest $m - |X_B|$ points of $X_\cR \setminus X_B$ to the center of $\supp B$\;
}\Else{
Define $\widehat{\m{q}}^B = \texttt{solve}(\widehat{\cL}_B, \cL^A, X, f(X), m, \mu)$\;
\For{$\widehat{B} \in \widehat{\cL}_B$}{
Express $\widehat{B}$ in terms of a set of tensor product B-splines with $B$ among them\;
Let $c_{\widehat{B}}$ be the coefficient of $B$ in such expression\;
}
Set $q_B = \sum_{\widehat{B} \in \widehat{\cL}_B} c_{\widehat{B}} \widehat{q}_{\widehat{B}}^B$\;
Continue\;
}}
Set $\bbol{C}_B$ as the collocation matrix at the points in $X_B$ with respect to the polynomial basis $\Pi_B^{\m{p}}$\;
Compute the smoothing matrix $\bbol{S}_B$ for $B$, as in Equation \eqref{eq:smoothingmatrix}\;
Assemble the least squares with smoothing matrix $\bbol{L}_B \defeq \bbol{C}_B^T\bbol{C}_B + \mu \bbol{S}_B$\;
\uIf{$\bbol{L}_B$ is ill-conditioned, that is, $\kappa(\bbol{L}_B) > \kappa_M$}{
Define $\widehat{\cL}_B$ as the set of mothers of $B$, extracted from $\cL^A$\;
Set $\widehat{\m{q}}^B = \texttt{solve}(\widehat{\cL}_B, \cL^A, X, f(X), m, \mu)$\;
\For{$\widehat{B} \in \widehat{\cL}_B$}{
Express $\widehat{B}$ in terms of a set of tensor product B-splines with $B$ among them\;
Let $c_{\widehat{B}}$ be the coefficient of $B$ in such expression\;
}
Set $q_B = \sum_{\widehat{B} \in \widehat{\cL}_B} c_{\widehat{B}} \widehat{q}_{\widehat{B}}^B$\;
}\Else{
Let $\m{z} \defeq f(X_B)$. Then, solve $\bbol{L}_B\m{u} = \bbol{C}_B^T\m{z}$\;
Let $\bbol{R}_B$ be the representation matrix of the basis $\Pi_B^\m{p}$ in terms of tensor product B-splines  $\cB_B$\;
Compute $\pmb{\lambda} \defeq \m{u}^T\bbol{R}_B$\;
Extract $\lambda_B$ from $\pmb{\lambda}$\;
Set $q_B = \lambda_B$\;
}}}
\caption{SOLVE module, $\m{q} = \texttt{solve}(\cL, \cL^A, X, f(X), m, \mu, \kappa_M)$}\label{alg:SOLVE}
\end{algorithm}
We highlight that if $|X| > m$ there always exist a sufficiently large condition number threshold $k_M$ and a coarse initial tensor mesh for which Algorithm \ref{alg:SOLVE} is able to produce an LR QI approximant, provided that the scattered data distribution is not degenerate (e.g., all aligned data-sites). In fact, under these assumptions the algorithm is able to assign a first set of QI coefficients for the B-splines on such coarse tensor mesh, without the need of enlargements or inheritances, that can be instead used for such routines when requested. 


\begin{oss}
There are two main reasons why we solve the local least squares problems in the power basis of polynomials and then we write the solution in terms of a local tensor B-spline basis, instead of solving them directly in this latter basis, as done for instance in \cite{thbscattered}. First, the assemble of the smoothing matrix is direct and very fast with the power polynomial basis, as one can see from Equation \eqref{eq:smoothingmatrix}. we just compute once the matrices $\bbol{S}^1, \bbol{S}^2, \bbol{S}^3$ and then, for any $B \in \cL$, scale them accordingly to the sizes of $\supp B$. As second reason, there would be not a straightforward way to set the tensor mesh out of the support of the LR B-spline at hand to define the local tensor B-spline basis, when enlarging the region where to seek for data points to $\cR$.
\end{oss}

\begin{figure}
\centering
\subfloat[SM, max err = 3.24e-02, RMSE = 2.86e-03]{\includegraphics[width = .25\textwidth]{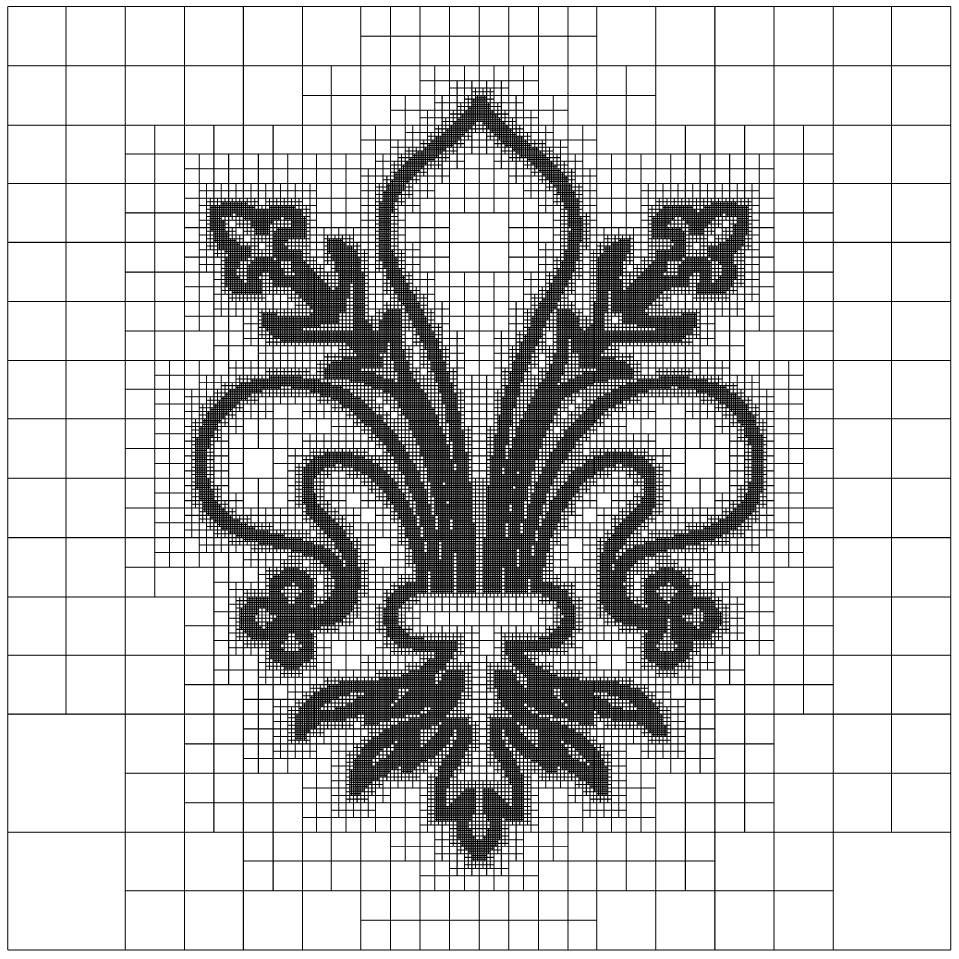}$\quad$
\includegraphics[width = .2\textwidth]{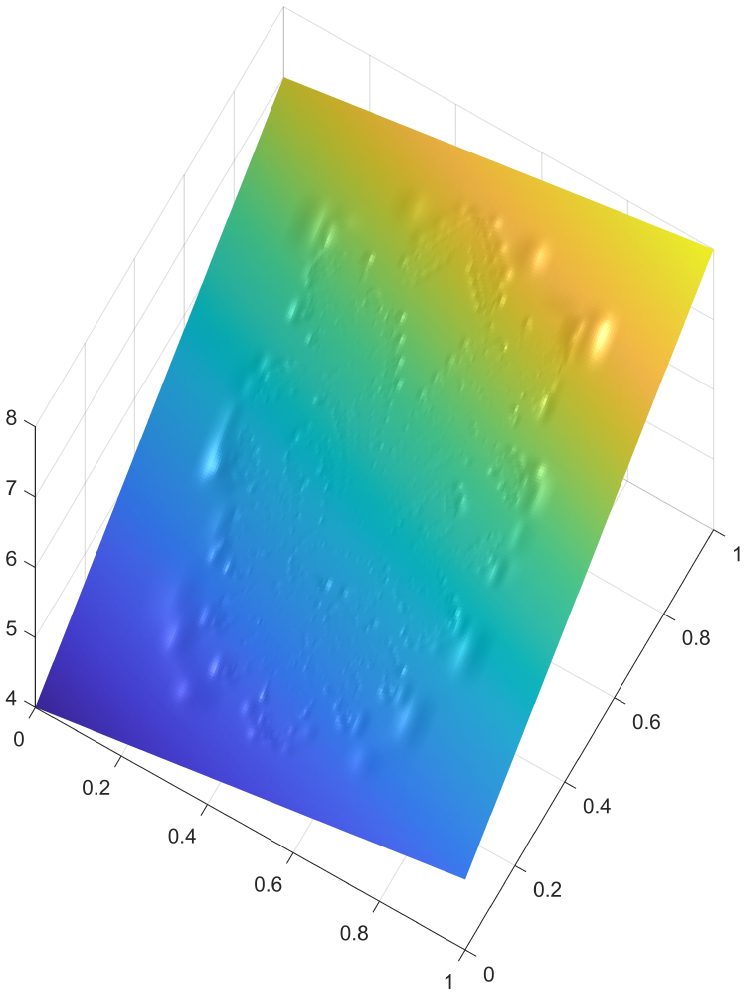}}\quad
\subfloat[N$_2$S$_2$, max err = 4.78e-13, RMSE = 5.70e-14]{\includegraphics[width = .25\textwidth]{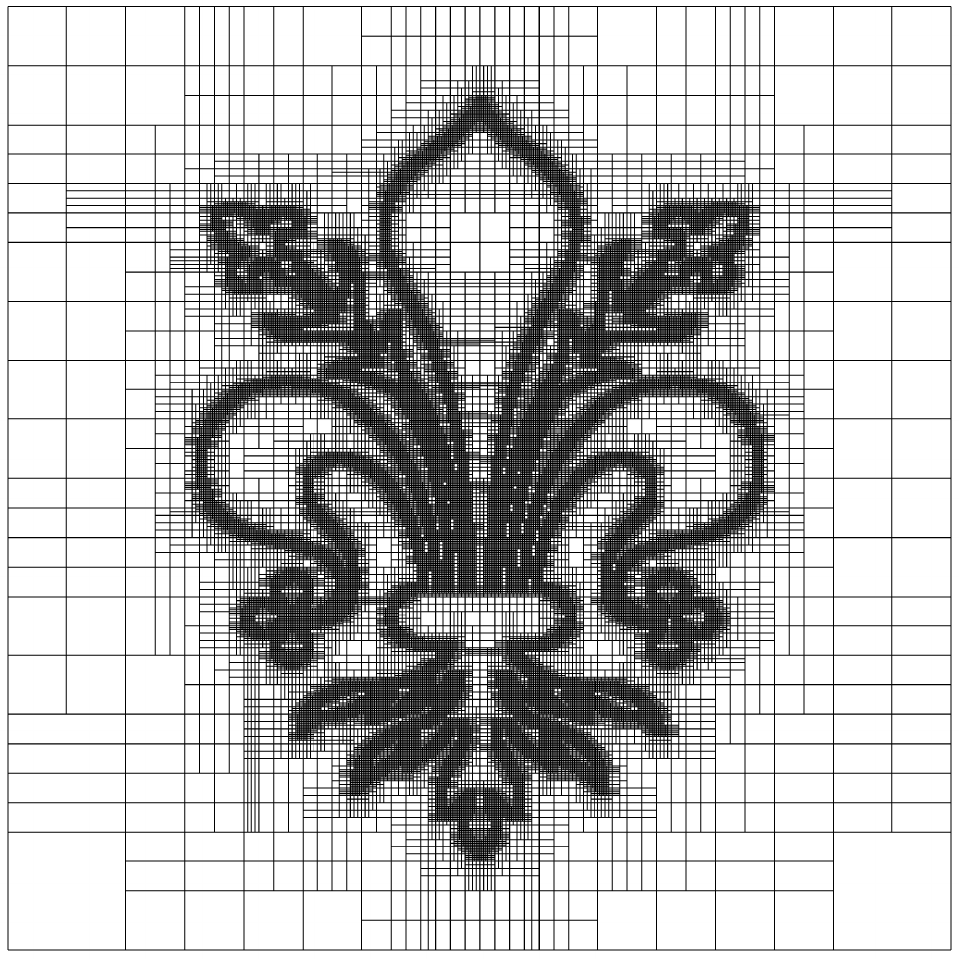}$\quad$
\includegraphics[width = .2\textwidth]{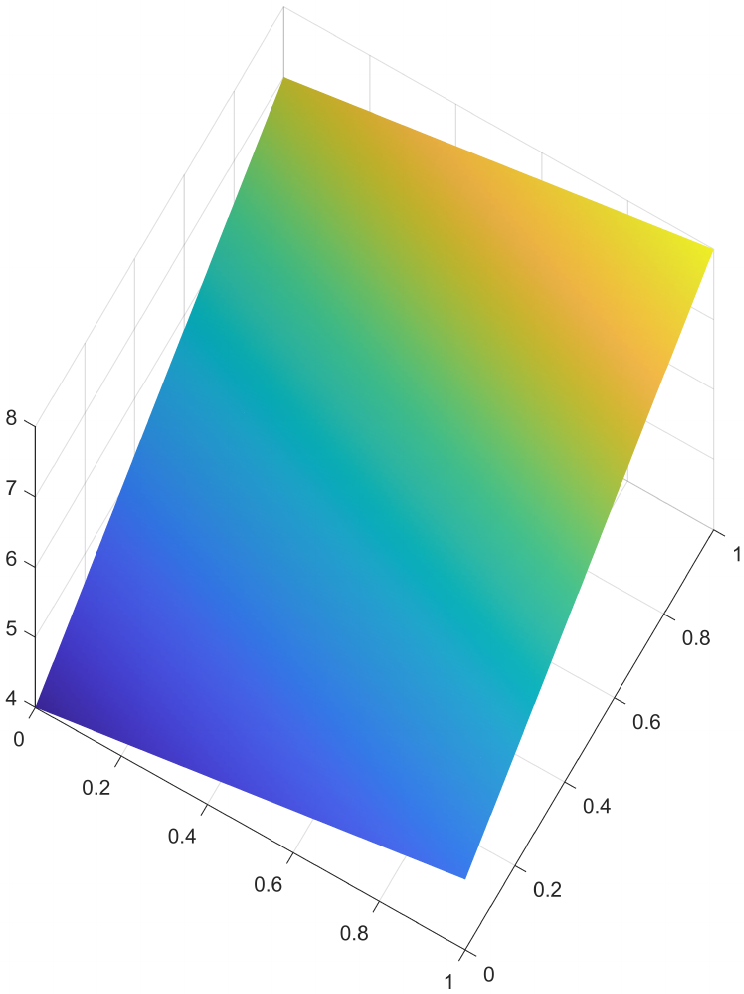}}
\caption{Linear reproduction without and with the N$_2$S property. Let $\m{p} = (2, 2)$. In figure (a) we use 9 levels of structured mesh (SM) refinement localized along the ``giglio of Florence'' curve to build an LR mesh in $[0, 1]^2$. In figure (b) we do the same with the N$_2$S$_2$ strategy. We test the reproduction of the linear polynomial $f(x, y) = x + 3y$ using the quasi-interpolation technique described in this section with the LR B-spline sets on such meshes. In order to do so, we generate 1 million Halton points and create a scattered dataset by evaluating $f$ at these points. As the structured mesh refinement lacks of the N$_2$S property, $f$ is not reproduced, or even well approximated, and artefacts appear on the computed surface. Instead, the N$_2$S$_2$ strategy ensures the reproduction of $f$, despite the smoothing parameter that here has been set as $\mu = 10^{-6}$. In the captions  we have reported the maximum (pointwise) error in the Halton points and the RMSE for the two approximations.}\label{fig:polireproduction}
\end{figure}

We conclude this section by investigating the polynomial reproduction capabilities of our method. The theoretical result in \cite[Proposition 4.4]{N2S2} ensures that, under the N$_2$S property of the mesh, for LR QI splines of bidegree $\p$, if the local approximation scheme used to define the coefficient associated with each LR basis function has polynomial reproduction up to bidegree $\p,$ the same property is inherited by the final LR QI spline. More in general, under the N$_2$S property, using the same argument one shows that if the local scheme ensures polynomial reproduction of bidegree up to $\m{d} \le \p$ or, even less, up to a total degree $d \le \max\{p_1,p_2\}$, then such reduced polynomial reproduction capability is shared with the final LR QI spline of bidegree $\p$. For this reason we are able to prove the next Proposition, in which the linear polynomial space shall be denoted by $\Pi_1$. Note that $\Pi_1 \subsetneq \Pi^{\m{1}}$ the space of bilinear polynomials.
\begin{prop}\label{prop:linearreproduction}
Assume $p_i \ge 1, i=1,2,$ and that Algorithm \ref{alg:SOLVE} is able to produce a QI of bidegree $\p$ in an LR spline space defined on an input LR mesh with the N$_2$S property. Then,  if $\mu >0$ and the function to be approximated  is a polynomial $P \in \Pi_1$,  then the created QI reproduces exactly $P$.
\end{prop}
\begin{proof}
We observe that the space of linear bivariate polynomials $\Pi_1$ is included in ${\cal L},$ since it is $p_i \ge 1, i=1,2.$ Thus, the restriction to $\Omega$ of any $P \in \Pi_1$ has a unique representation in ${\cal L}$, that is, $P_{|\Omega} = \sum_{B\in {\cal L}}P_B B$. Furthermore, the restriction $P_{|\supp B}$ can be also represented in any tensor product B-spline basis ${\cal B}_B$ in $\supp B$. In the Algorithm, we choose ${\cal B}_B$ such that it includes $B$ and, as already noticed in \cite[Proposition 4.4]{N2S2}, this implies that the the coefficient assigned to $B$ for representing $P|_{\supp B}$ is still $P_B.$ \\
Denoting with  $\mathfrak{Q}(P) = \sum_{B\in {\cal L}}q_B B$ the QI spline produced by the algorithm when applied to $P$, what we need then to prove is that $\mathfrak{Q}(P) = P$, that is, $q_B = P_B, \forall B \in {\cal L}.$ Now, considering the control flow in the algorithm, essentially two different strategies can be adopted for defining $q_B$, depending on the size of the local set $X_B \subset X$ and also matrix conditioning considerations. The first consists in solving a local approximation problem which relies on functional information at the subset $X_B \subset X$ (or at a suitable enlargement of $X_B$, i.e., $X_\cR$).  More precisely, referring for simplicity to $X_B$, in this case $q_B$ is defined as the coefficient assigned to $B$ in the basis ${\cal B}_B$ of the polynomial $g = g(P) \in \Pi_\p,$ such that $g = g(P) \defeq \argmin_{Q \in \Pi_\p} A_B(P,Q)  + \mu I_B(Q) \,$, with
$$\scalebox{0.85}{$
A_B(P,Q) \defeq \displaystyle\sum_{(x, y) \in X_B} \left( Q(x, y) - P(x, y) \right)^2\,, \quad
I_B(Q) \defeq \displaystyle\int_{\supp B} \left(\partial_{xx} Q (x,y)\right)^2 + 2 \left(\partial_{xy} Q (x,y)\right)^2 + \left(\partial_{yy} Q (x,y)\right)^2\,\, \text{dxdy}\,.$}
$$
Now, being $\mu>0$, such polynomial $g$ is unique, for details we refer to \cite[Appendix]{thbscattered}, and it necessarily coincides with $P$, since $A_B(P,P) = 0$ and also $I_B(P) = 0$, being $P\in \Pi_1.$ Hence the coefficient $q_B$ coincides with $P_B$ in this case.
In the second alternative, $q_B$ is inherited from coefficients associated with mothers of $B$ belonging to $\widehat{\cL}$ which is a coarsening of ${\cal L}$ and it is extracted from the archive $\cL^A$ of all the LR B-spline sets in input to the algorithm. In this case the proof can be obtained by using an induction argument and considering that $P_B$ can be also defined exactly through the same inheritance formula, starting from the expression of $P$ in the coarsest LR space in ${\cal L}^A.$
\end{proof}
We consider such linear reproduction capability of the scheme highly demanded and indispensable in the context of shape approximation, as a faithful reconstruction of the flat areas of a surface can be accomplished only through it. 
This is highlighted by Figure \ref{fig:polireproduction} where two different approximations of a plane are compared. They are both generated by our SOLVE module but rely on two different LR meshes. Even if both are refined in the same (arbitrary) regions, one is obtained with the N$_2$S$_2$ refinement strategy and the other with the structured mesh approach  \cite{johannessen} which does not guarantee the N$_2$S property. One can clearly see the artefacts (bumps) in the approximation when using the LR B-splines on the structured mesh.    

\label{sec:solve}

\section{Numerical Tests}\label{sec:exm}
In this section we show the performance of the quasi-interpolation scheme for LR B-splines of Section \ref{sec:qi} in some numerical tests. The location and number of iterations of the refinements on the LR meshes will be provided by the fault detections and jump estimates of Section \ref{sec:jstim}.

In order to verify the quality of the jump driven approach, we shall compare the results with the purely fault guided method, in which the refinements are performed everywhere along the fault curves detected. This procedure has been proved to be reliable \cite{BGC23}, providing approximation precision close to the error driven algorithms.

In the tests of this section we always consider scattered datasets composed of one million points, either generated by evaluation of functions at pseudo-random Halton points or acquired by LiDAR or similar technologies from real world landscapes and openly distributed on the internet. We quantify the quality of the approximation comparing the RMSE of the reconstructed surface using our jump based approach and the RMSE of the surface obtained with the purely fault guided method, both computed at all the scattered data-sites. In addition to saving a remarkable number of degrees of freedom while achieving an error of the same magnitude, as we shall see, it also turns out that the jump based method is faster than the fault based procedure both for assembling the quasi-interpolant and for evaluating it, as also shown in Table \ref{tab:speedup} for the two most challenging datasets considered, due to the smaller number of basis functions to handle. The marking and refining steps are significantly faster as well, as most of the fault points are dismissed after applying $\bar{L}$ refinement steps in the jump based algorithm. 
For what concerns the settings, we fix an initial $2 \times 2$ tensor mesh in all the examples. We describe the results and show the approximations in the figures for bidegree $\m{p} =(3, 3)$. However, similar considerations hold for  $\m{p}=(2, 2)$, as one shall see from Table \ref{tab:bidegree22}.
The maximal jump class will always be $L = 10$, while different minimal classes $\bar{L}$ will be adopted for the same experiment as reported in related tables. However, $\bar{L}$ will be the same both for ordinary and gradient faults. The smoothing parameter is fixed to $\mu = 10^{-6}$ and the minimum number of points $m$ for setting the local least squares problem is fixed to $m = 5$. Furthermore, when considering anisotropic refinements, we start performing one directional insertions as soon as the refinement iteration reaches the minimal jump class, $L^A = \bar{L}$. Finally, the least squares matrix $\bbol{L}$ shall be considered ill-conditioned if its condition number is above $\kappa_M = $1e05. However, this occurence has never happened in the reported tests.
Of course, these are not the only possible, and perhaps not optimal, settings: for instance, one may choose different minimal and maximal jump classes for ordinary and gradient faults, as well as an independent class from which to refine anisotropically. 
Nevertheless, the results show the potential of this kind of technique.
\subsection{Approximation of synthetic data from functional sampling at scattered points}\label{sec:funtests}
We test the algorithm for surface reconstruction from scattered data provided by the sampling of two functions, namely
\begin{equation}\label{eq:SS}
f_1(x, y) \defeq \left\{\begin{array}{ll}
\left|x - \frac{2}{5} - \frac{1}{10}\sin(2\pi y) \right| & \text{if } x \leq \frac{7}{10} + \frac{1}{5} \sin(2\pi y),\\\\
\left|x - \frac{2}{5} - \frac{1}{10}\sin(2\pi y) - \frac{1}{5}\right| & \text{otherwise,}
\end{array}\right.
\end{equation}
sampled in $[0, 1]^2$, and
\begin{equation} \label{eq:BranchCuts}
f_2(x, y) \defeq \text{Im}(\arcsin(z) + \arcsin(i\bar{z}))
\end{equation}
with $z \defeq x + i y$, sampled in $[-4, 4]^2$ ($f_2$ is a so-called \emph{branch cut}, details can be found in, e.g., \cite{arcoseno}). 

\begin{figure}
\centering
\subfloat[]{\includegraphics[height = 0.2\textwidth]{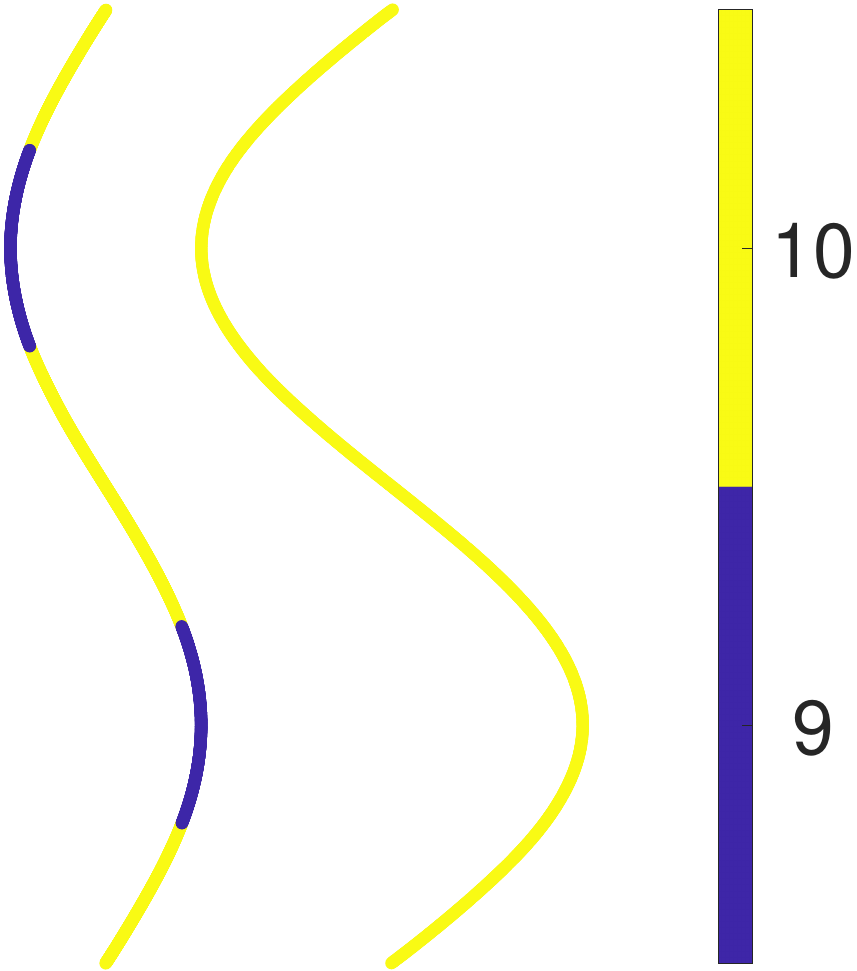}}\hspace{2cm}
\subfloat[ndofs 66634, RMSE 3.56e-03]{
\includegraphics[width = 0.2\textwidth]{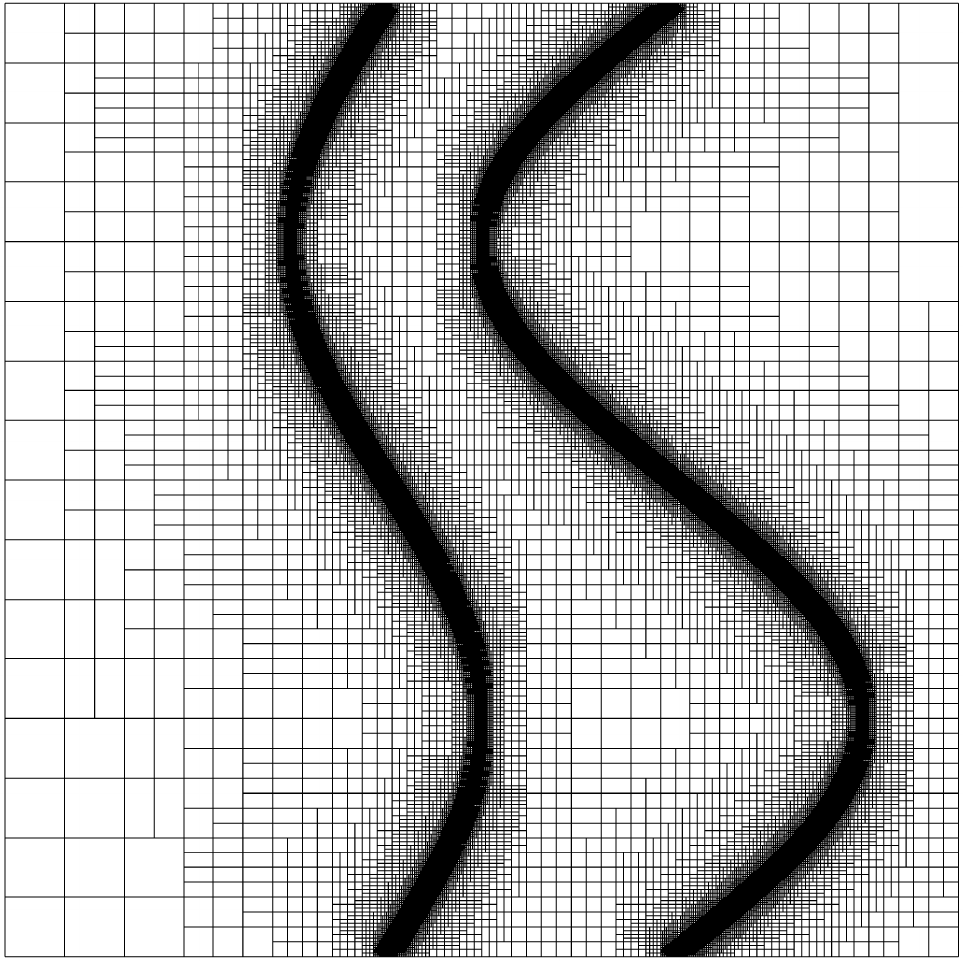}$\quad$\raisebox{.05\textwidth}{\includegraphics[width = 0.25\textwidth]{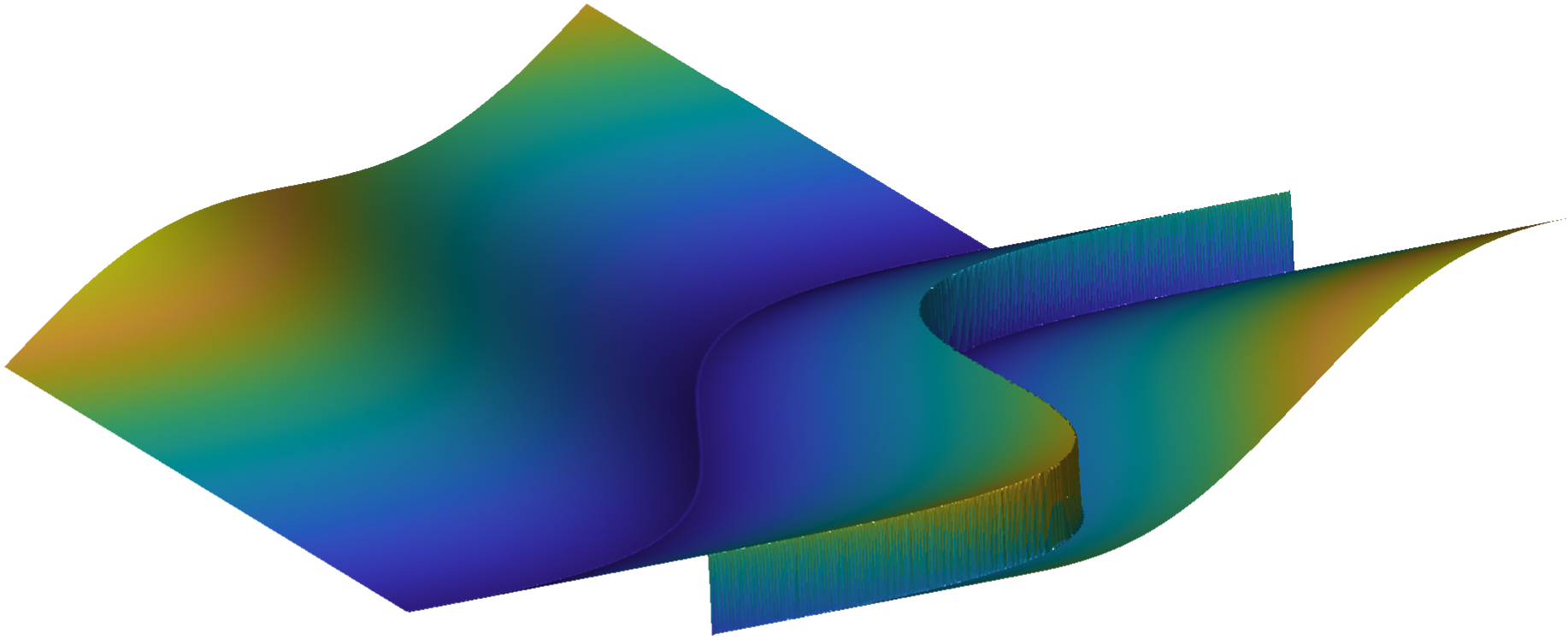}}
}\\
\subfloat[ndofs 61393, RMSE 3.56e-03]{
\includegraphics[width = 0.2\textwidth]{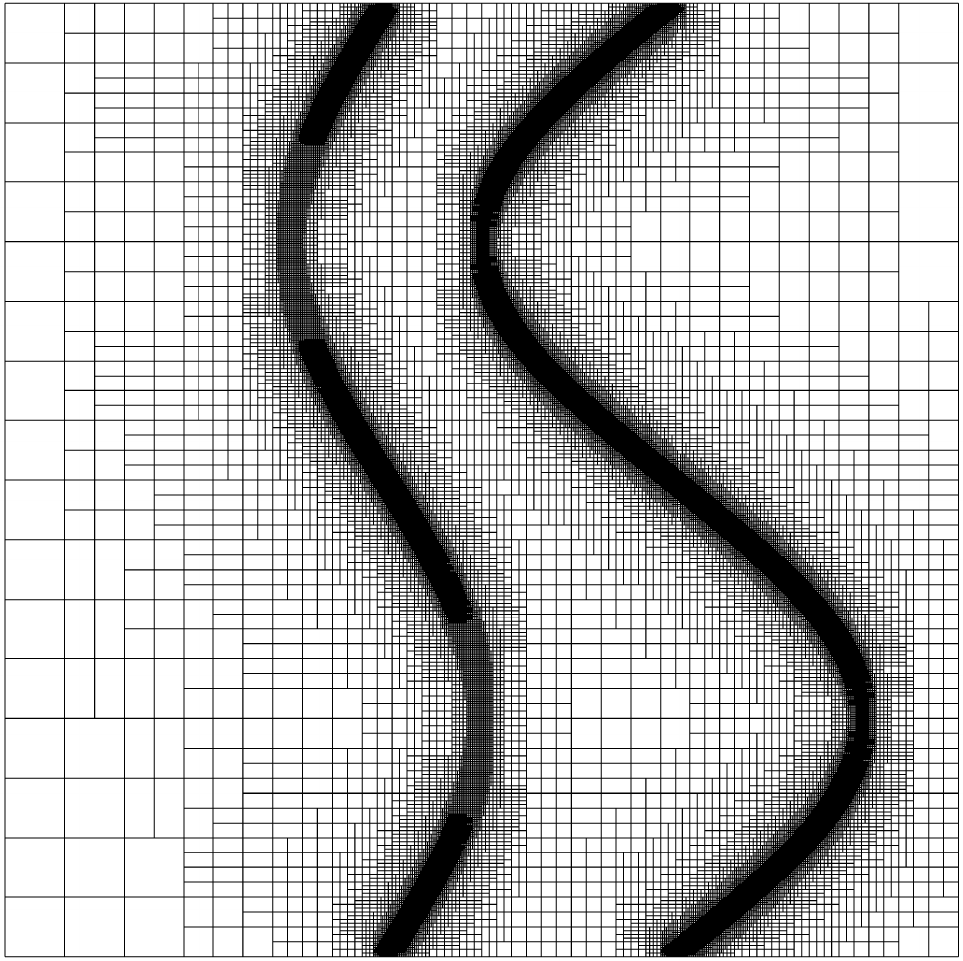}$\quad$\raisebox{.05\textwidth}{\includegraphics[width = 0.25\textwidth]{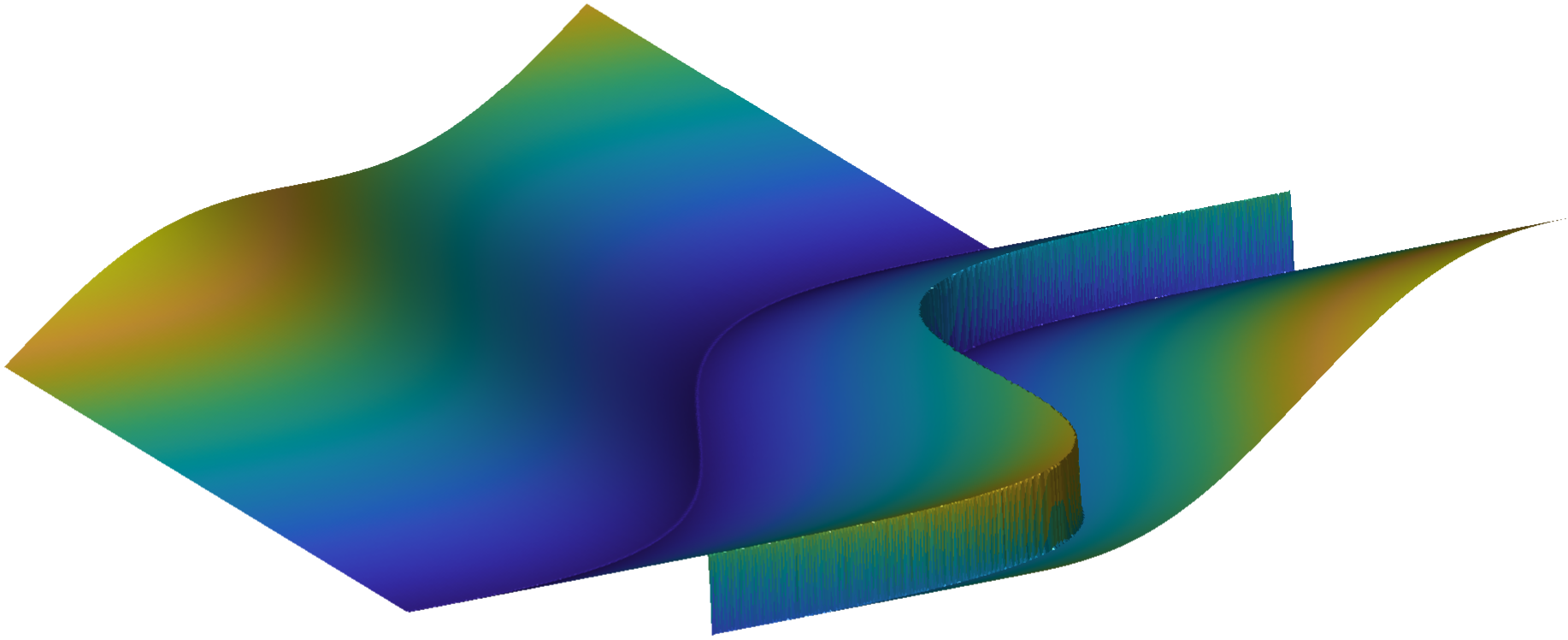}}
}\quad
\subfloat[ndofs 58907, RMSE 3.57e-03]{
\includegraphics[width = .2\textwidth]{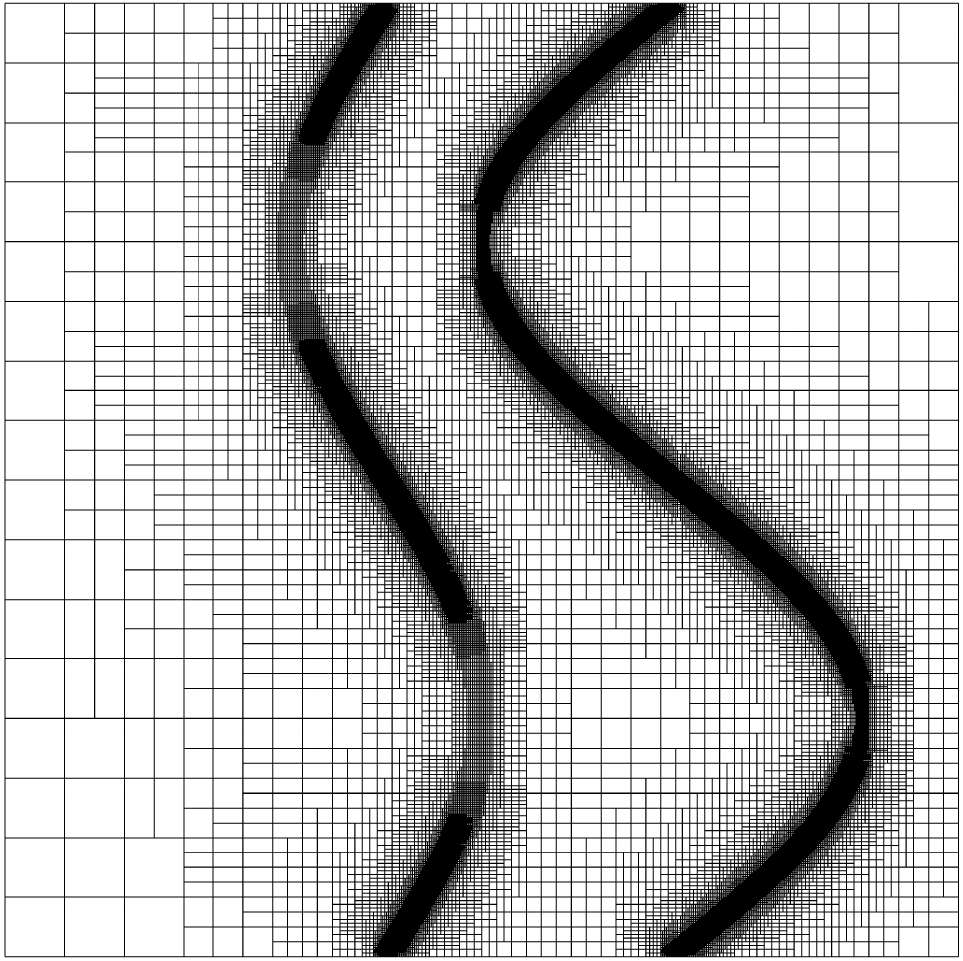}$\quad$\raisebox{.05\textwidth}{\includegraphics[width = 0.25\textwidth]{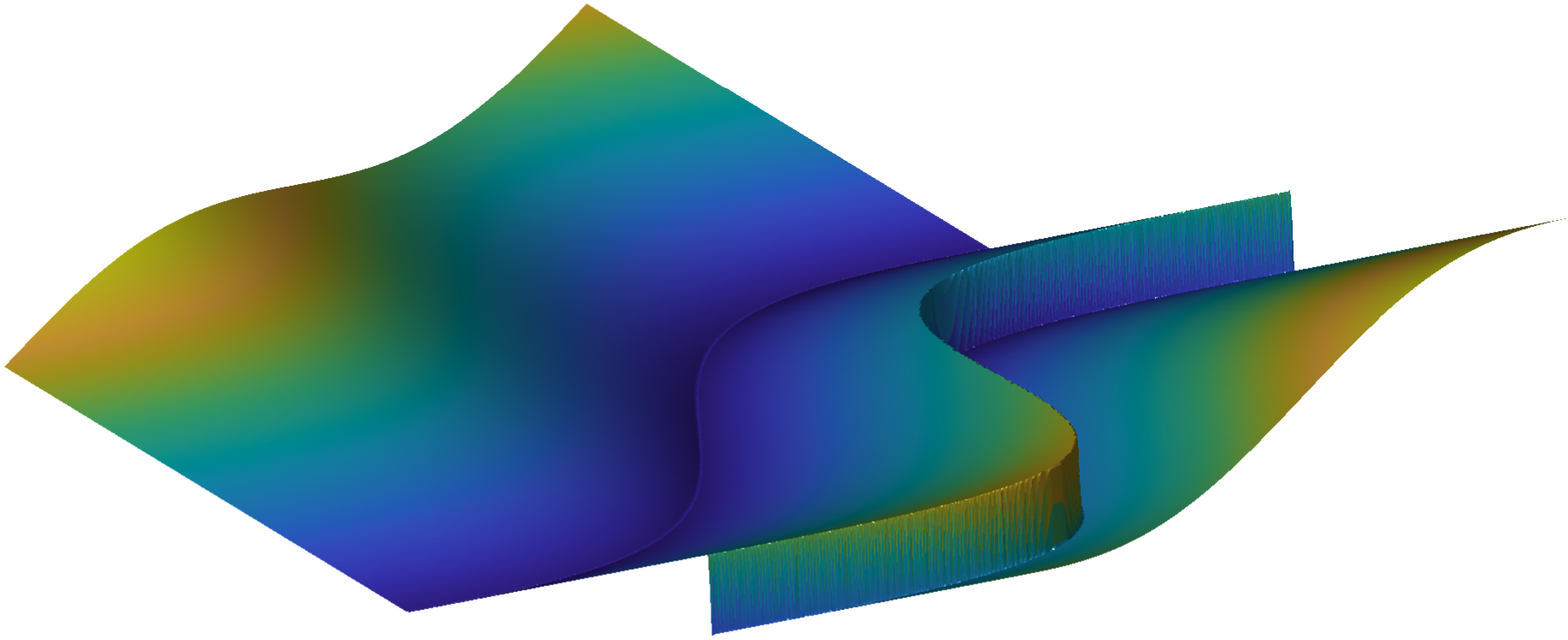}}}
\caption{QI approximation of scattered points from function \eqref{eq:SS}. Such function has two sinusoidal fault curves. The one on the right is a constant ordinary fault and the one on the left is a varying gradient fault. The points detected for the ordinary fault have all been classified at the maximum level of refinement $\ell = 10$ by the jump classification algorithm. The points detected for the gradient fault have instead been split into two classes corresponding to levels $\ell = 9$ and $\ell = 10$, see figure (a). In figures (b) we see mesh and approximation provided by the purely fault based refinement. The number of dofs and RMSE are reported in the caption of the figure. In figures (c) and (d) we show the results of the isotropic and anisotropic jump driven refinements, respectively.}\label{fig:SS}
\end{figure}

$f_1$ restricted to $[0, 1]^2$ has range $[0, 0.5]$ and it is characterized by two faults along sinusoidal curves: one ordinary fault, with a constant jump of $\frac{1}{5}$, and one varying gradient fault. The results of the jump driven quasi-interpolations are illustrated in Figure \ref{fig:SS} and compared with the purely fault guided strategy. Here, the advantages of the jump driven strategies are mitigated by the fact that the jump range is relatively small compared to its overall magnitude. This leads to a small variety in the jump class distribution of the detected fault points. Moreover, the axis-aligned pieces of the fault curve are quite small and, for the gradient fault curve, related to the lower jump class, which reduces the impact of the anisotropic insertions. Nevertheless, as one can see from the figure, essentially the same RMSE of the fault driven algorithm is reached both with the isotropic and anisotropic refinement approaches, despite saving approximately 5200 and 7700 degrees of freedom, respectively, i.e., a reduction of about $\nicefrac{1}{12}$ and $\nicefrac{1}{9}$ of the total.

\begin{figure}
\centering
\subfloat[]{\includegraphics[height = 0.2\textwidth]{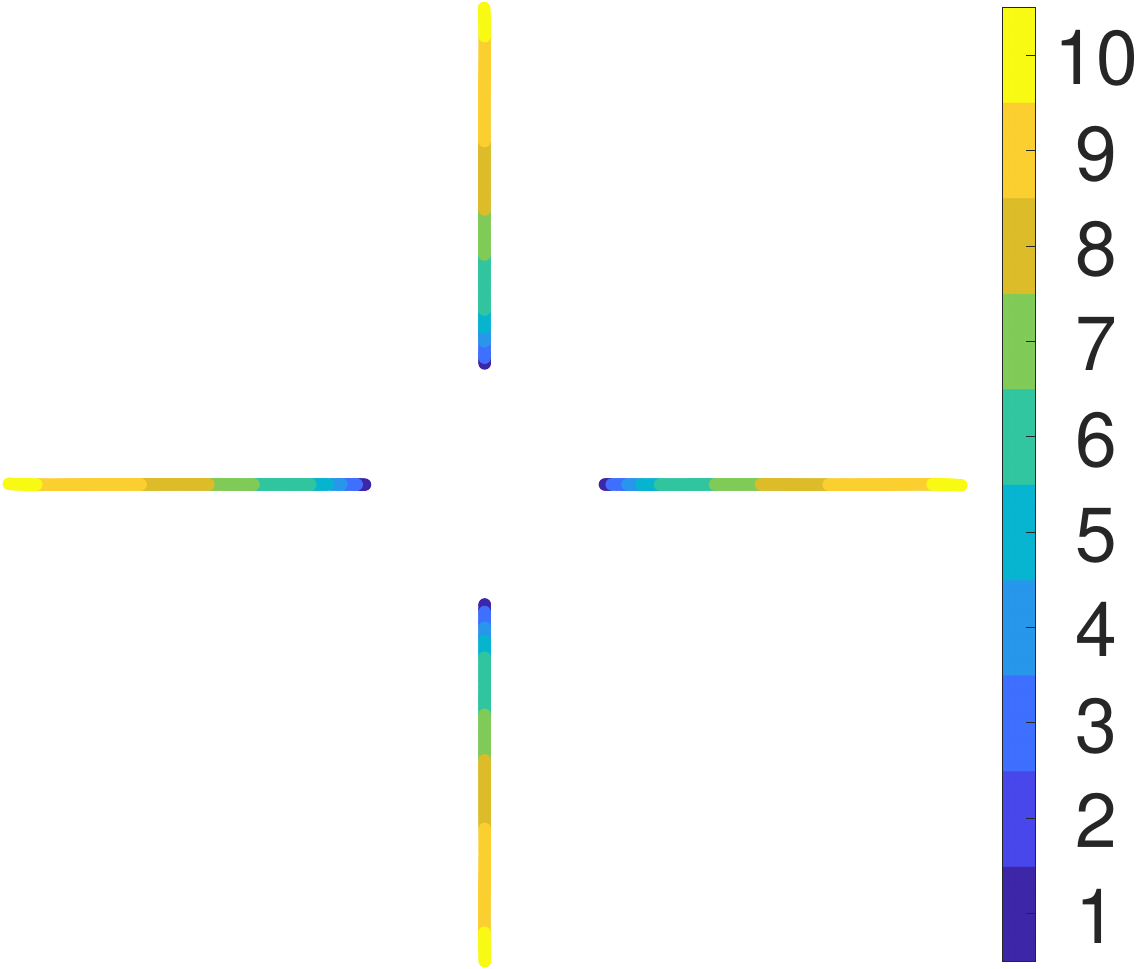}}\qquad\qquad
\subfloat[ndofs 30935, RMSE 6.06e-02]{
\includegraphics[width = .2\textwidth]{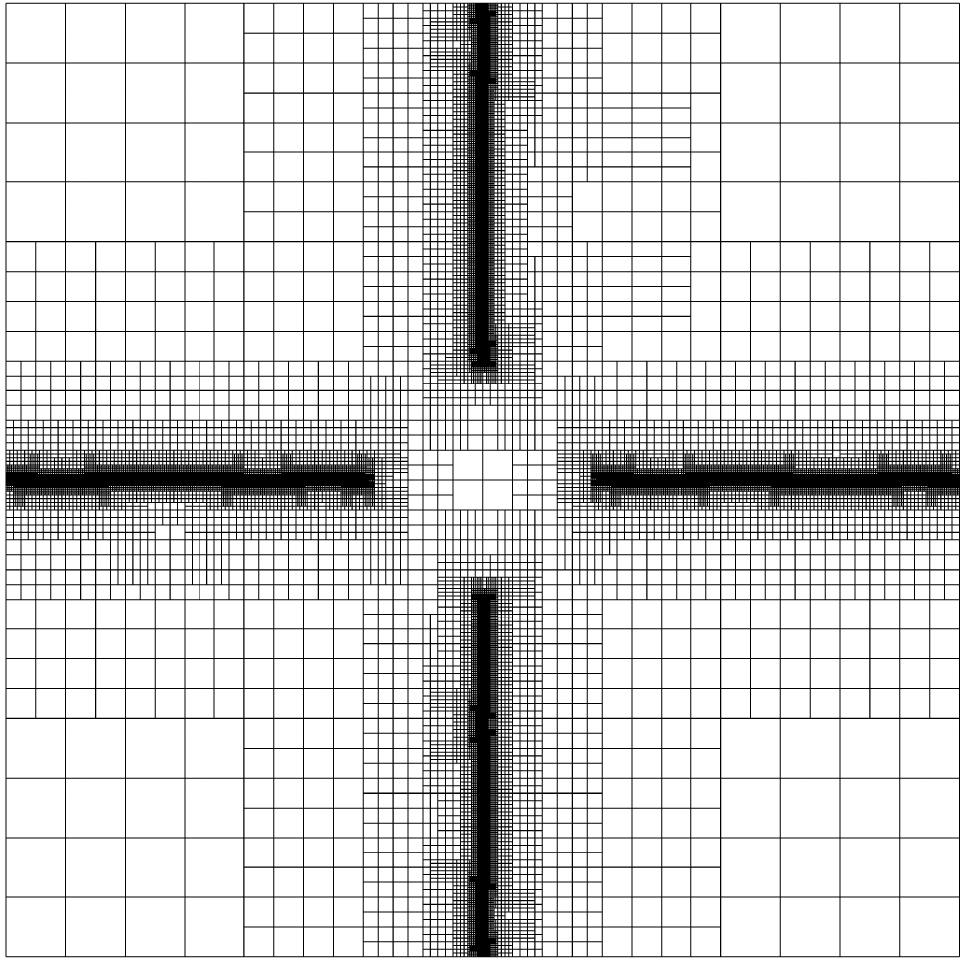}$\quad$\includegraphics[width = 0.2\textwidth]{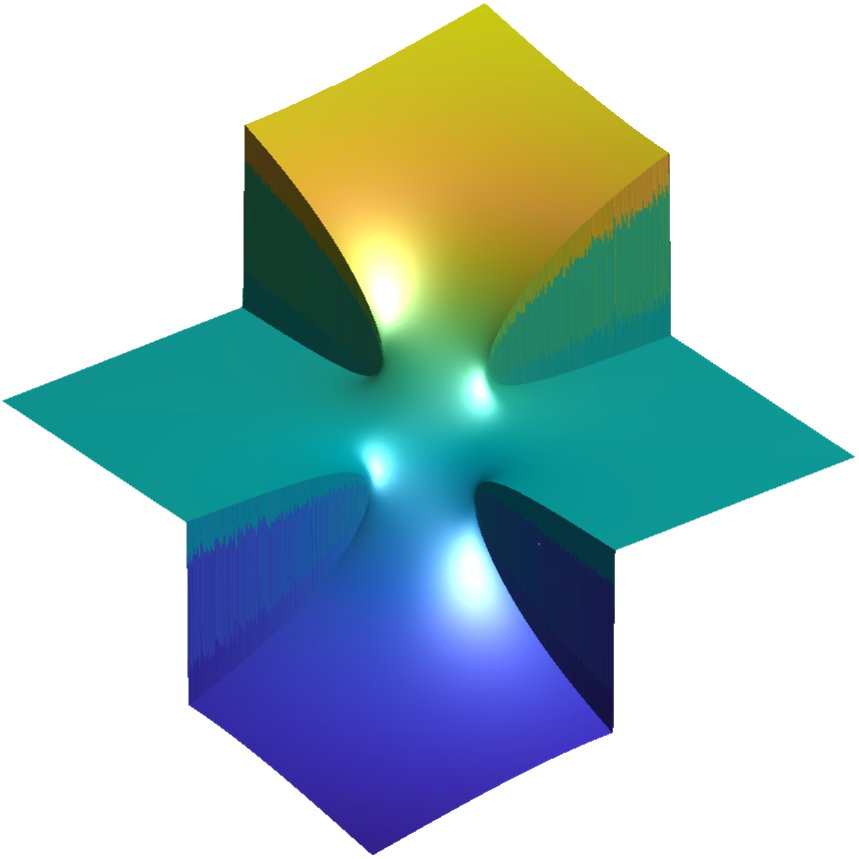}}\vspace{.25cm}\\
\subfloat[]{\scalebox{.85}{\renewcommand{\arraystretch}{2}\begin{tabular}{|c|c|c|}
\hline
$\bar{L}$ & ndofs & RMSE \\ \hline
1 			& 1558 & 8.37e-02 \\ \hline
4			& 1696 & 8.37e-02 \\ \hline
8			& 7310 & 7.28e-02 \\ \hline
\end{tabular}}}\qquad\qquad\quad\raisebox{-1.4cm}{
\subfloat[$\bar{L} = 8$, ndofs 7310, RMSE 7.28e-02]{
\includegraphics[width = 0.2\textwidth]{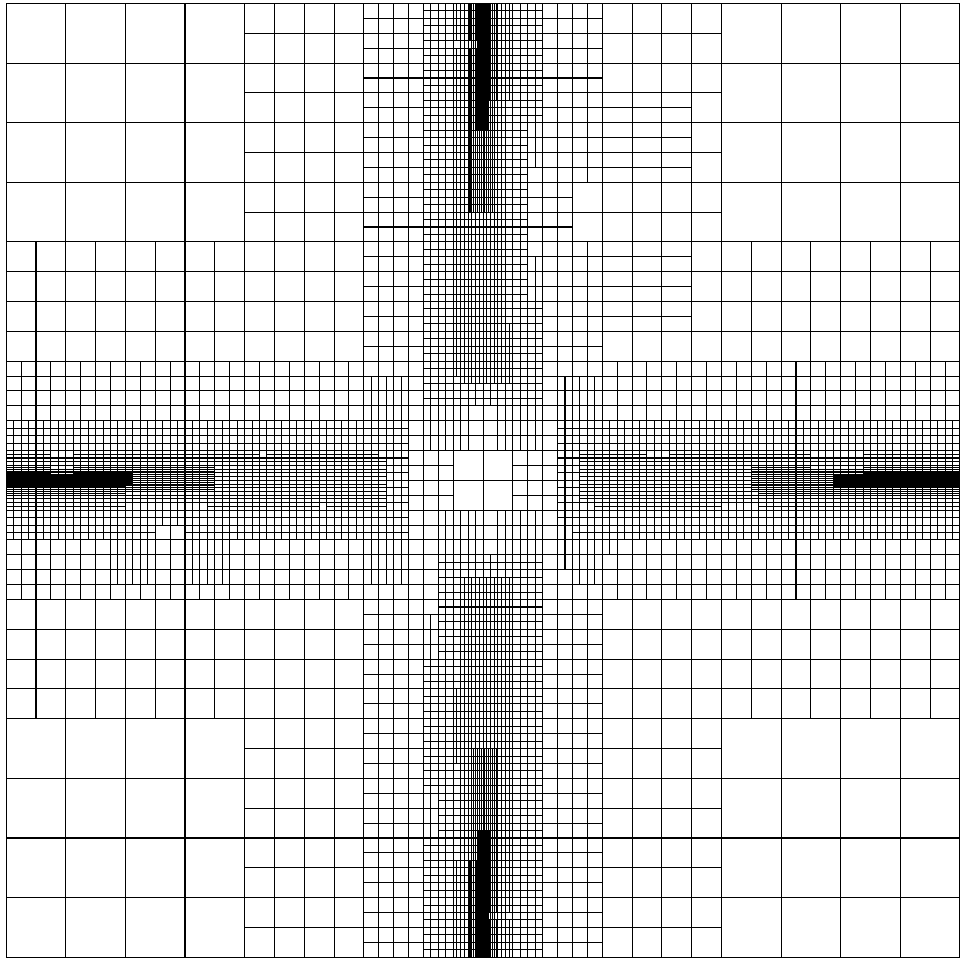}$\quad$\includegraphics[width = 0.2\textwidth]{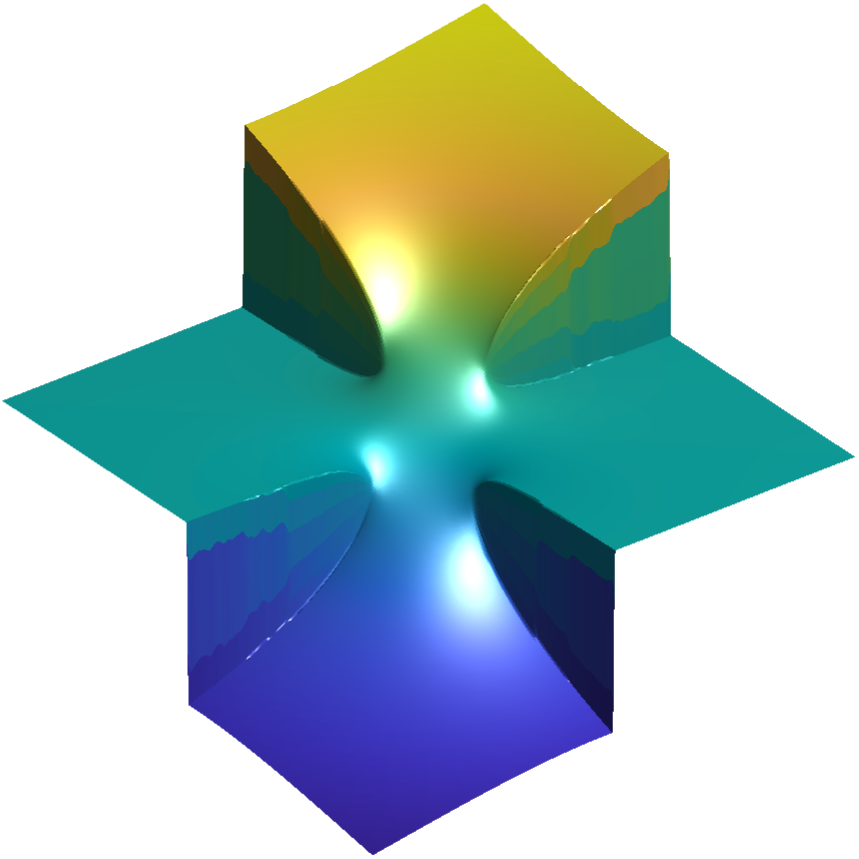}
}}
\caption{QI approximation of scattered data from function \eqref{eq:BranchCuts}.
The surface is characterized by four symmetric axis-aligned ordinary faults of decreasing intensity towards the center. The points detected as fault points can be distributed on jump classes $\ell \in \{\bar{L}, \ldots, L\}$, from any $\bar{L} \in \{1, \ldots, 10\}$, using the jump classification algorithm. In figure (a), we see the classification for $\bar{L} = 1$. In figure (b) we show mesh and approximation when using the purely fault driven refinement, as it will be used for comparison with the jump based approach. In the caption of the figure we have recorded the number of degrees of freedom and RMSE. In table (c) we have instead the results for the anisotropic jump driven refinement for different jump classifications, obtained by picking different $\bar{L}$. 
In figure (d) we present mesh and approximation for the case $\bar{L} = 8$.}\label{fig:BranchCuts}
\end{figure}

$f_2$ restricted to $[-4, 4]^2$ has (approximate) range $[-4.85, 4.85]$ and it shows four symmetric axis-aligned ordinary faults of varying jumps, with a decreasing intensity towards the center, where it is actually continuous. The results of the jump driven quasi-interpolation with anisotropic refinements are reported in Figure \ref{fig:BranchCuts} and compared with the fault guided approach. This function has been chosen in order to highlight and stress the possible terrific reduction in the degrees of freedom, while retaining the same magnitude of RMSE, by allowing anisotropic insertions when (most of) the detected fault points can be labelled as axis-aligned, possibly after a suitable rotation of the point cloud, in combination with jump driven marking. Indeed, in this test we can reach essentially the same RMSE of the purely fault guided refinement with approximately $\nicefrac{1}{4}$ (for $\bar{L} = 8$), $\nicefrac{1}{18}$ (for $\bar{L} = 4$) and $\nicefrac{1}{20}$ (for $\bar{L} = 1$), respectively, of the degrees of freedom.

\begin{figure}
\centering
\subfloat[]{
\includegraphics[width = 0.25\textwidth]{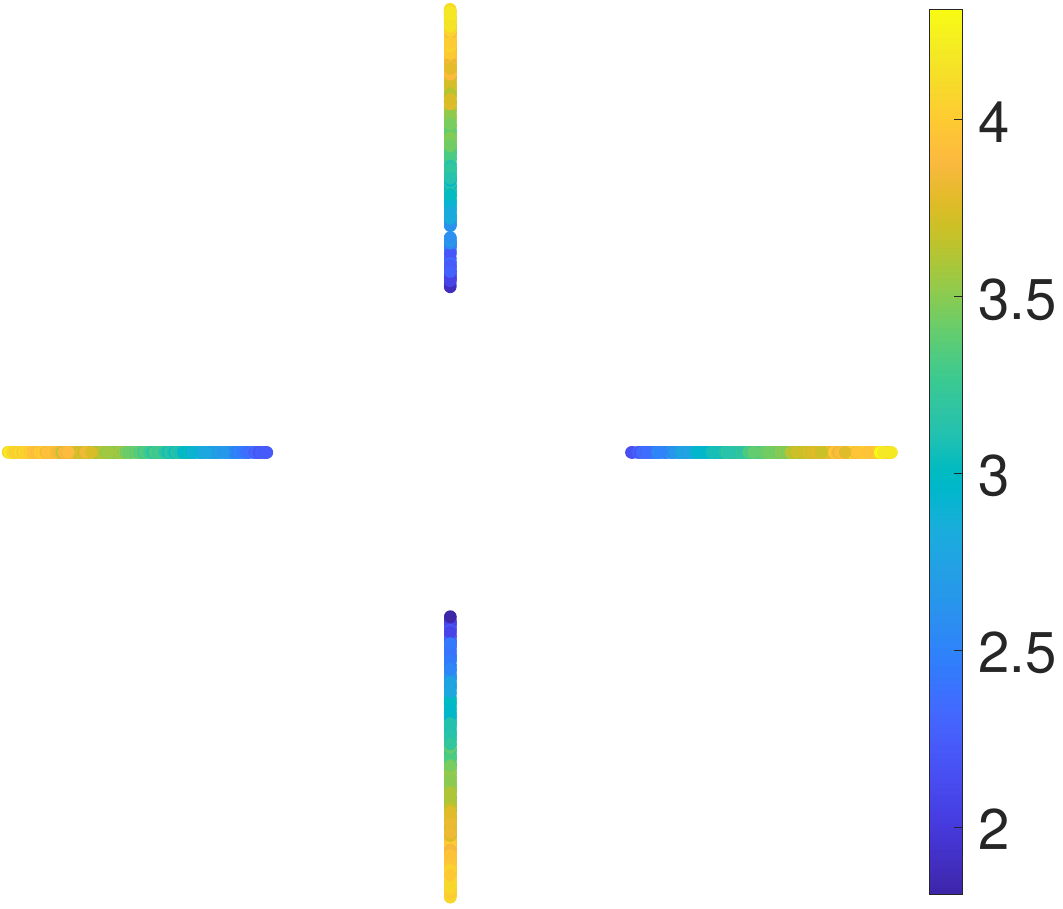}
}\qquad
\subfloat[]{
\includegraphics[width = 0.25\textwidth]{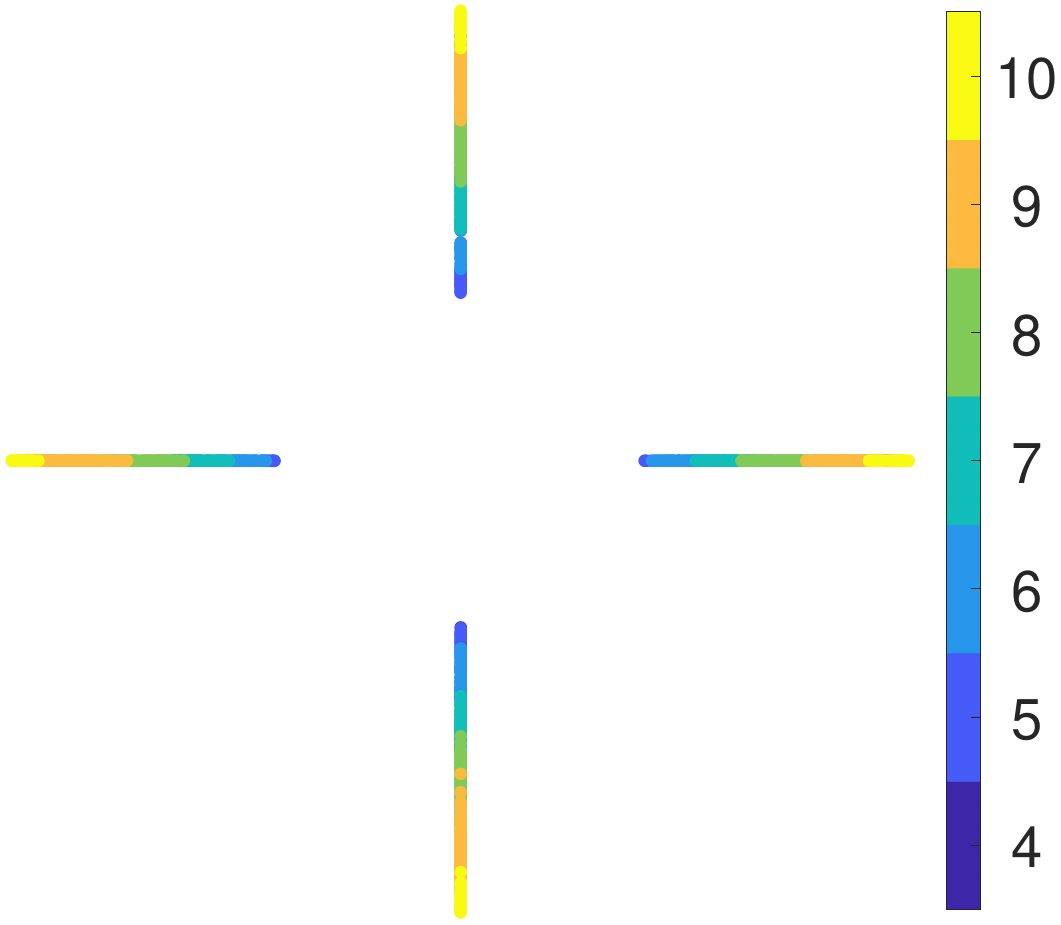}
}\qquad\raisebox{1.3cm}{
\subfloat[]{
\scalebox{.75}{\renewcommand{\arraystretch}{2}\begin{tabular}{|c|c|c|c|c|}
\hline
 & \multicolumn{2}{|c|}{$\m{p} = (2, 2)$} & \multicolumn{2}{|c|}{$\m{p} = (3, 3)$} \\ \hline
$\bar{L}$ & ndofs & RMSE & ndofs & RMSE\\ \hline
4 			 & 912 & 9.60e-02 & 1463 & 1.05e-01\\ \hline
6			 & 1847 & 8.71e-02 & 2592 & 9.01e-02\\ \hline
8			 & 5289 & 8.06e-02 & 7377 & 8.39e-02\\ \hline
\end{tabular}}}}
\caption{A preliminary test for noisy datasets. We have perturbed the values of the dataset created for function \eqref{eq:BranchCuts} with a Gaussian noise of magnitude 0.25, i.e., about $5.70\%$ of the maximal jump value. In figure (a) the jump estimates. In figure (b) the corresponding jump classification. In (c) the table providing the results for degrees $\m{p} = (2, 2)$ and $\m{p} = (3, 3)$ when choosing $\bar{L} = 4, 6, 8$ respectively.}\label{fig:noise}
\end{figure}

This latter function has been adopted also for carrying a first, preliminary, test in the presence of  data affected by noise. Namely, we have perturbed the dataset sampled from function \eqref{eq:BranchCuts} with a Gaussian distributed noise of magnitude 0.25, that is, about the 5.70\% of the maximum jump value of the unperturbed function. Figure \ref{fig:noise} summarizes the results. As the jump becomes lower toward the center of the surface, the detection struggles to distinguish it from the noise. Furthermore the noise is handed down to the jump estimates to some extent. This latter drawback has however little-to-none effect on the jump guided refinement ultimately. For what concerns the former, the loss of the part of fault where the jump is lower can be compensated by considering a finer initial tensor mesh. This is advisable also for mitigating the effect of the noise. Such finer mesh shall be reasonable for having a fair approximation of the fault where the jump is low, as the local refinements there would stop pretty soon anyway as lower jump values correspond to lower jump classes.

\subsection{Approximation of landscape datasets}
\begin{figure}
\centering
\subfloat[]{\includegraphics[height = 0.19\textwidth]{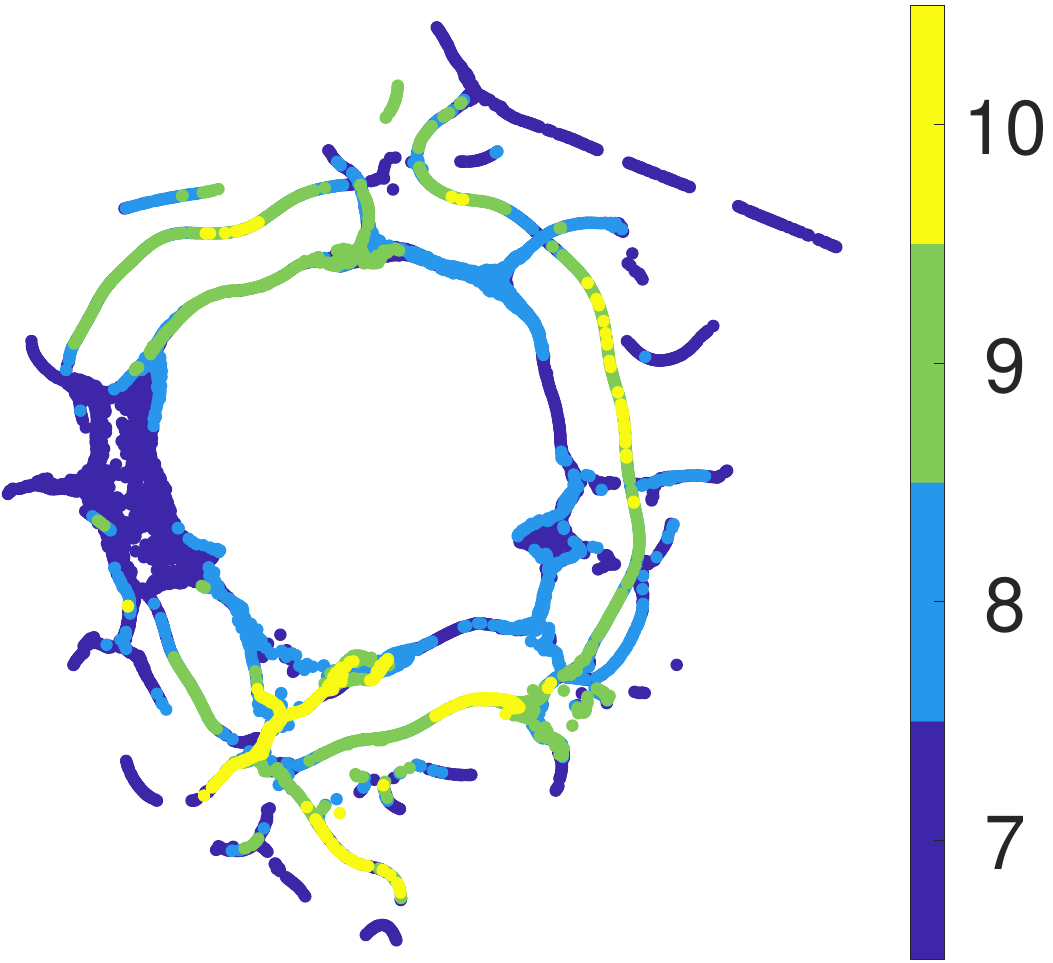}}\quad
\subfloat[ndofs 183802, RMSE 2.22e-04]{
\includegraphics[width = .2\textwidth]{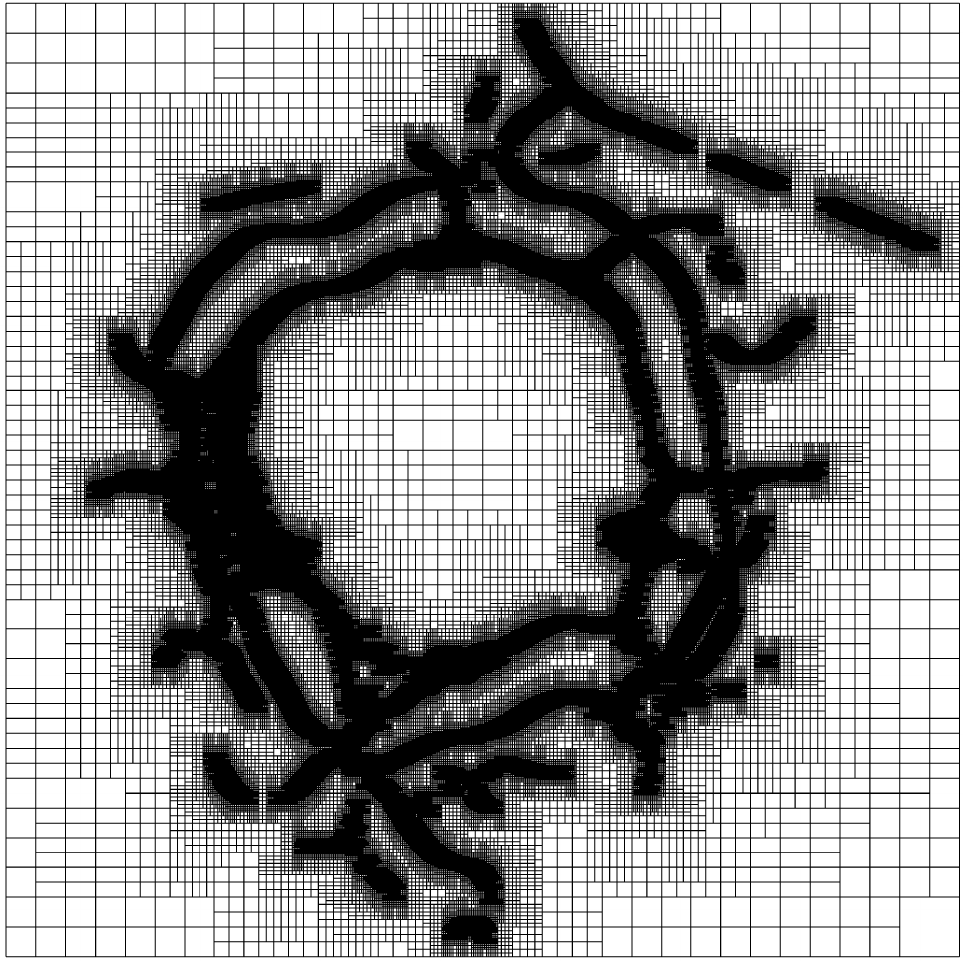}$\quad$\includegraphics[width = 0.25\textwidth]{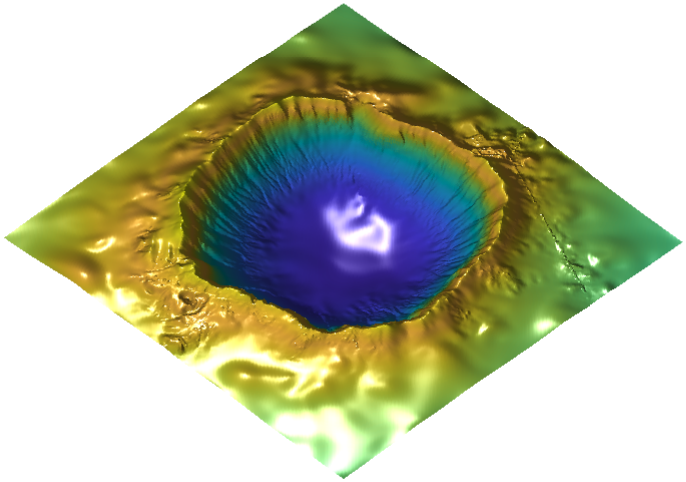}}\quad\raisebox{1.25cm}{
\subfloat[]{\scalebox{.65}{\renewcommand{\arraystretch}{2}\begin{tabular}{|c|c|c|c|c|}
\hline
 & \multicolumn{2}{|c|}{isotropic} & \multicolumn{2}{|c|}{anisotropic} \\ \hline
$\bar{L}$ & ndofs & RMSE & ndofs & RMSE\\ \hline
5 			 & 8174 & 7.70e-04& 6725 & 7.99e-04\\ \hline
6			 & 15425 & 4.39e-04& 13396 & 4.55e-04\\ \hline
7			 & 29910 & 2.93e-04& 26897 & 3.06e-04\\ \hline
\end{tabular}}}}\\
\subfloat[$\bar{L} = 7$, isotropic, ndofs 29910, RMSE 2.93e-04]{
\includegraphics[width = .2\textwidth]{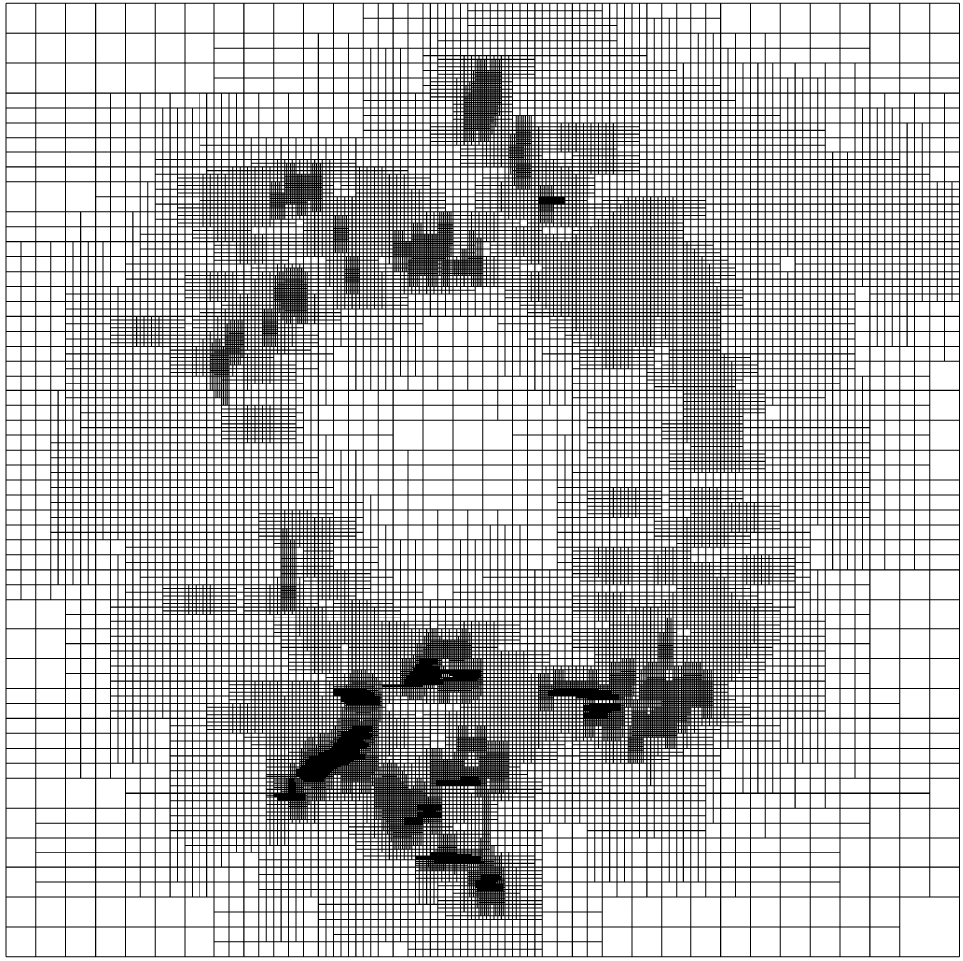}$\quad$\includegraphics[width = 0.25\textwidth]{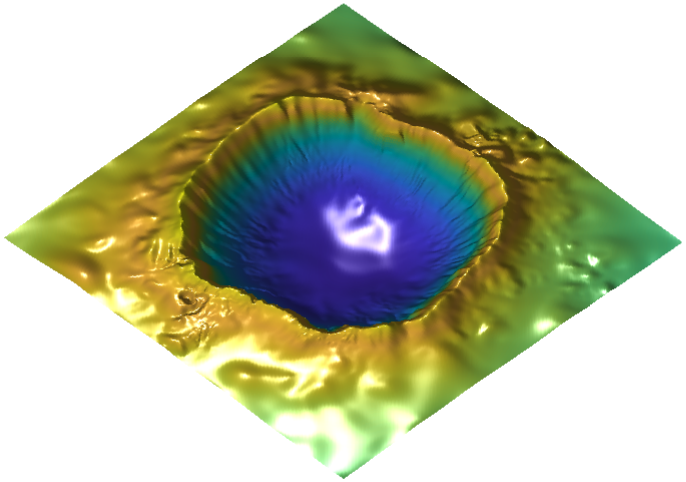}}
\quad
\subfloat[$\bar{L} = 7$, anisotropic, ndofs 26897, RMSE 3.06e-04]{
\includegraphics[width = .2\textwidth]{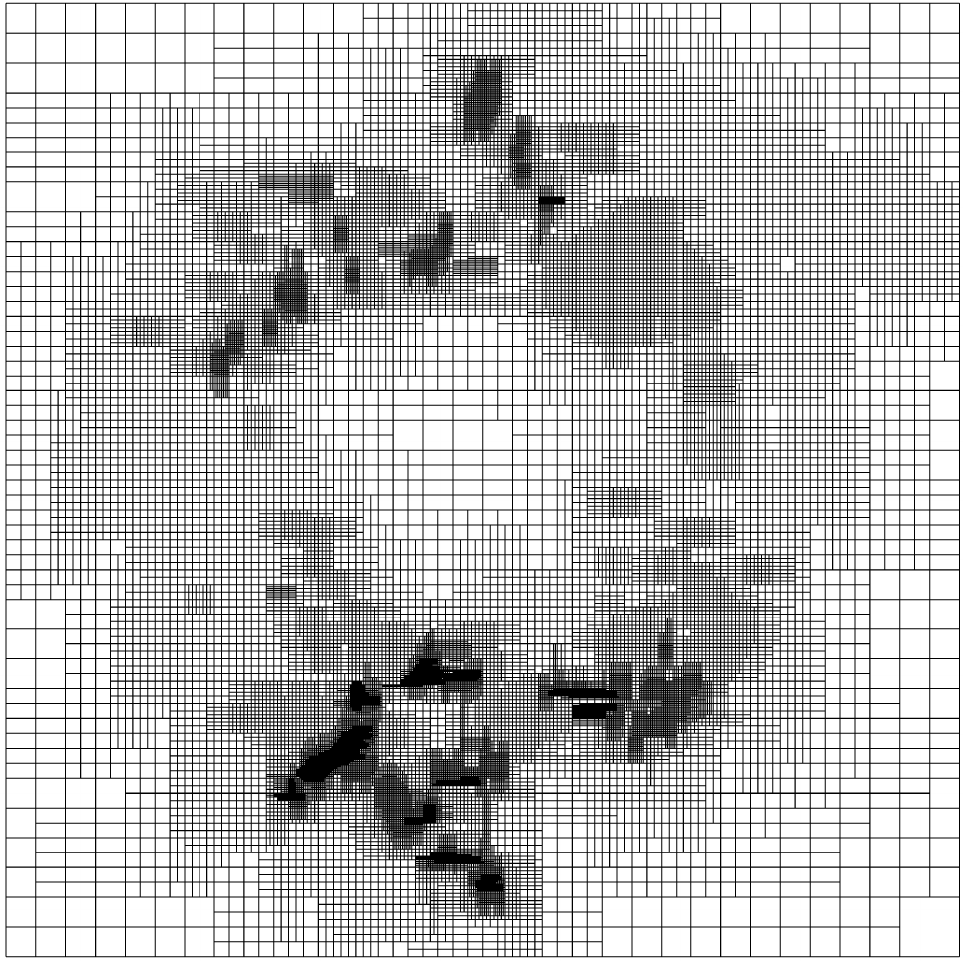}$\quad$\includegraphics[width = 0.25\textwidth]{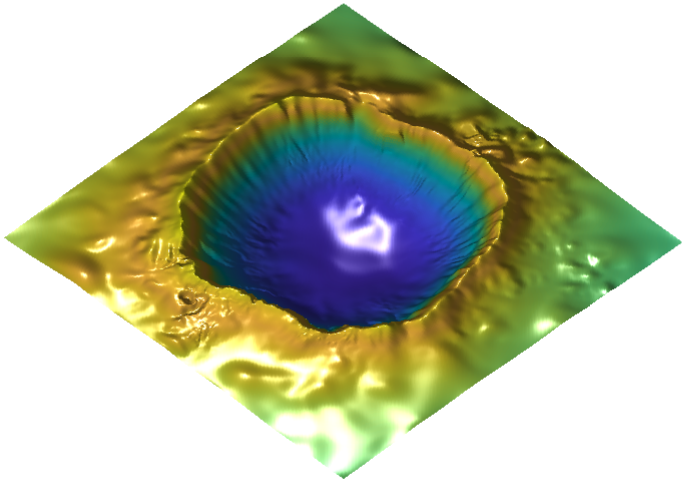}}
\caption{QI approximation of a one million scattered point cloud of the Barringer Crater in Arizona, USA \cite{MC}. In figure (a) the jump classification for $\bar{L} = 7$. In figure (b) we have mesh and approximation using the purely fault driven approach to localize the refinement, when setting the maximal level to $L = 10$. The number of degrees of freedom and RMSE are listed in the caption of the figure. In table (c) we report the results for the isotropic and anisotropic jump driven refinements, for different choices of $\bar{L}$ in the jump classification. In figures (d)--(e) we present meshes and approximations of the two approaches in the instance of $\bar{L} = 7$.}\label{fig:MC}
\end{figure}

In this section we test our jump driven quasi-interpolation on some datasets acquired by the scanning of real world landscapes. The first point cloud we consider has been captured at the Barringer Crater in Arizona \cite{MC}. The sampling has been rescaled to be in the range $[0.5, 0.6]$. Figure \ref{fig:MC} shows the results both for isotropic and anisotropic refinements, compared to the fault driven approach in figure (b). In figure (a) we show the jump classified when $\bar{L} = 7$. The RMSE remains of the same magnitude for any choice of $\bar{L} \geq 5$ and the amount of details is also qualitatively comparable, as one can see from figures (d) and (e) for $\bar{L} = 7$. On the other hand, we achieve a remarkable reduction in the degrees of freedom as illustrated by the table in figure (c).

\begin{figure}
\centering
\subfloat[]{\includegraphics[height = 0.19\textwidth]{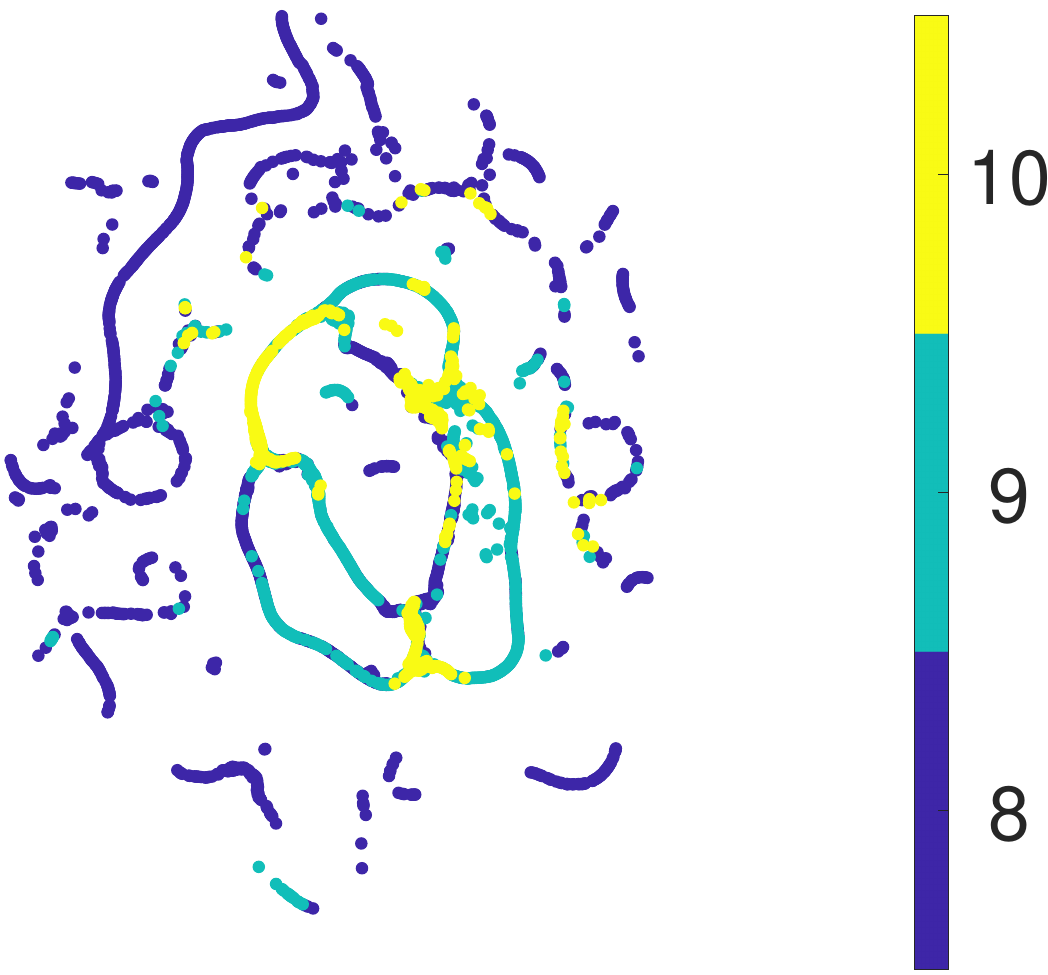}}\quad
\subfloat[ndofs 148700, RMSE 2.31e-03]{
\includegraphics[width = .2\textwidth]{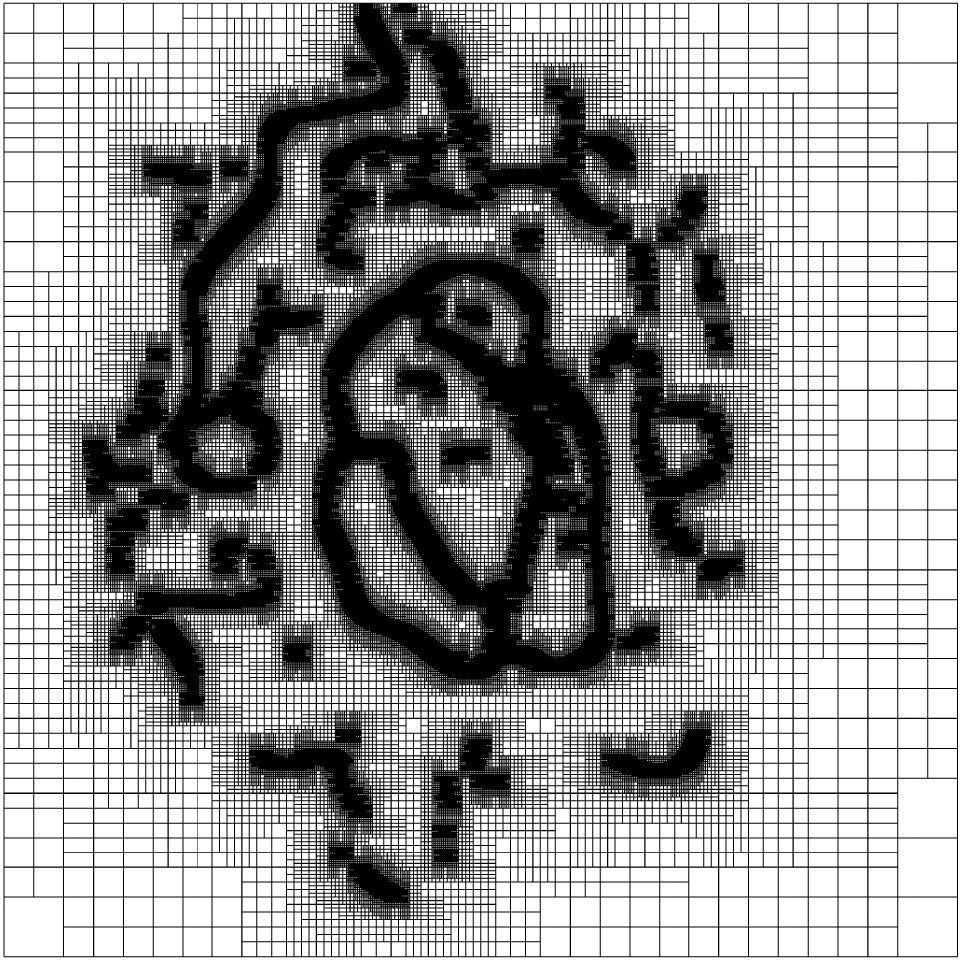}$\quad$\raisebox{1.25cm}{\includegraphics[width = 0.25\textwidth]{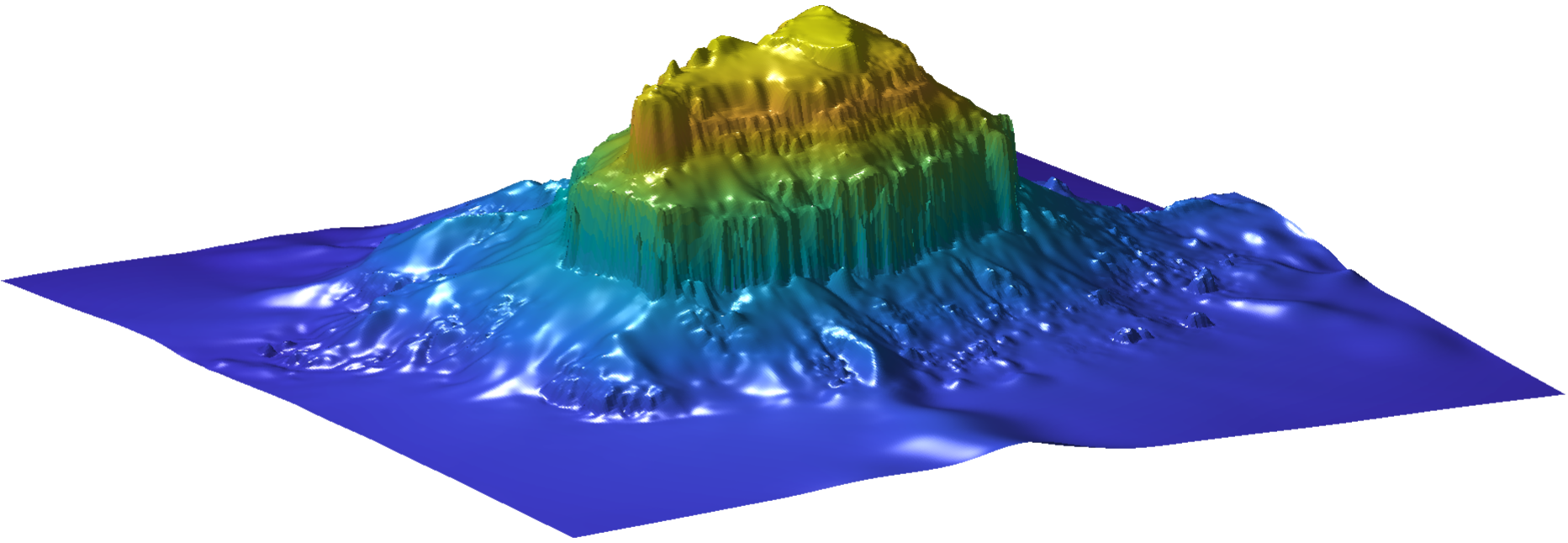}}}\quad\raisebox{1.25cm}{
\subfloat[]{\scalebox{.65}{\renewcommand{\arraystretch}{2}\begin{tabular}{|c|c|c|c|c|}
\hline
 & \multicolumn{2}{|c|}{isotropic} & \multicolumn{2}{|c|}{anisotropic} \\ \hline
$\bar{L}$ & ndofs & RMSE & ndofs & RMSE\\ \hline
6 			 & 16820 & 2.87e-03& 14536 & 2.93e-03\\ \hline
7			 & 27734 & 2.59e-03& 24237 & 2.63e-03\\ \hline
8			 & 51607 & 2.42e-03& 45324 & 2.47e-03\\ \hline
\end{tabular}}}}\\
\subfloat[$\bar{L} = 8$, isotropic, ndofs 51607, RMSE 2.42e-03]{
\includegraphics[width = .2\textwidth]{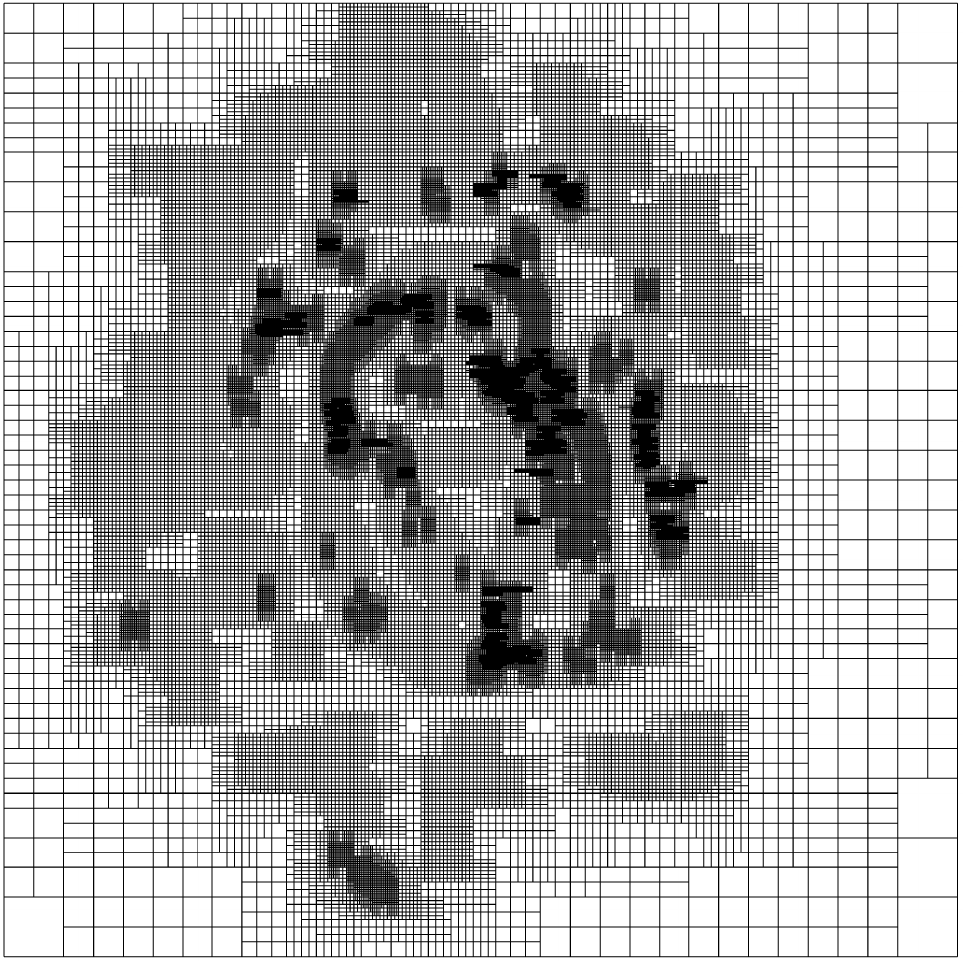}$\quad$\raisebox{1.25cm}{\includegraphics[width = 0.25\textwidth]{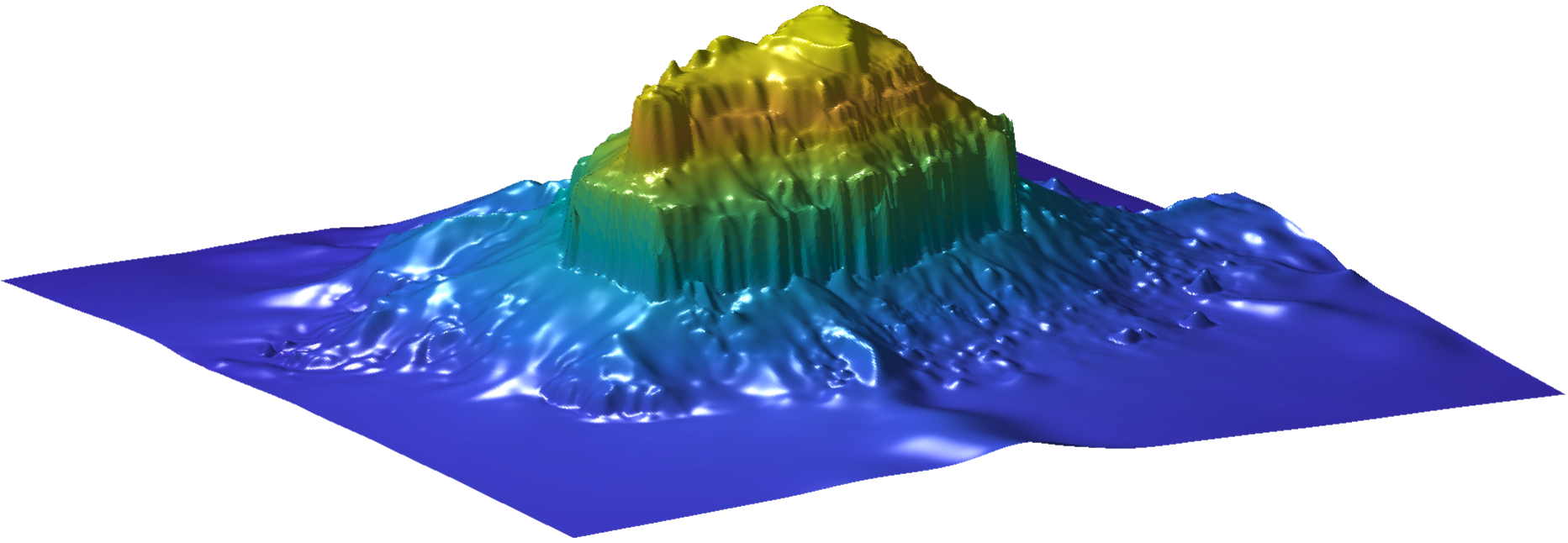}}}
\quad
\subfloat[$\bar{L} = 8$, anisotropic, ndofs 45324, RMSE 2.47e-03]{
\includegraphics[width = .2\textwidth]{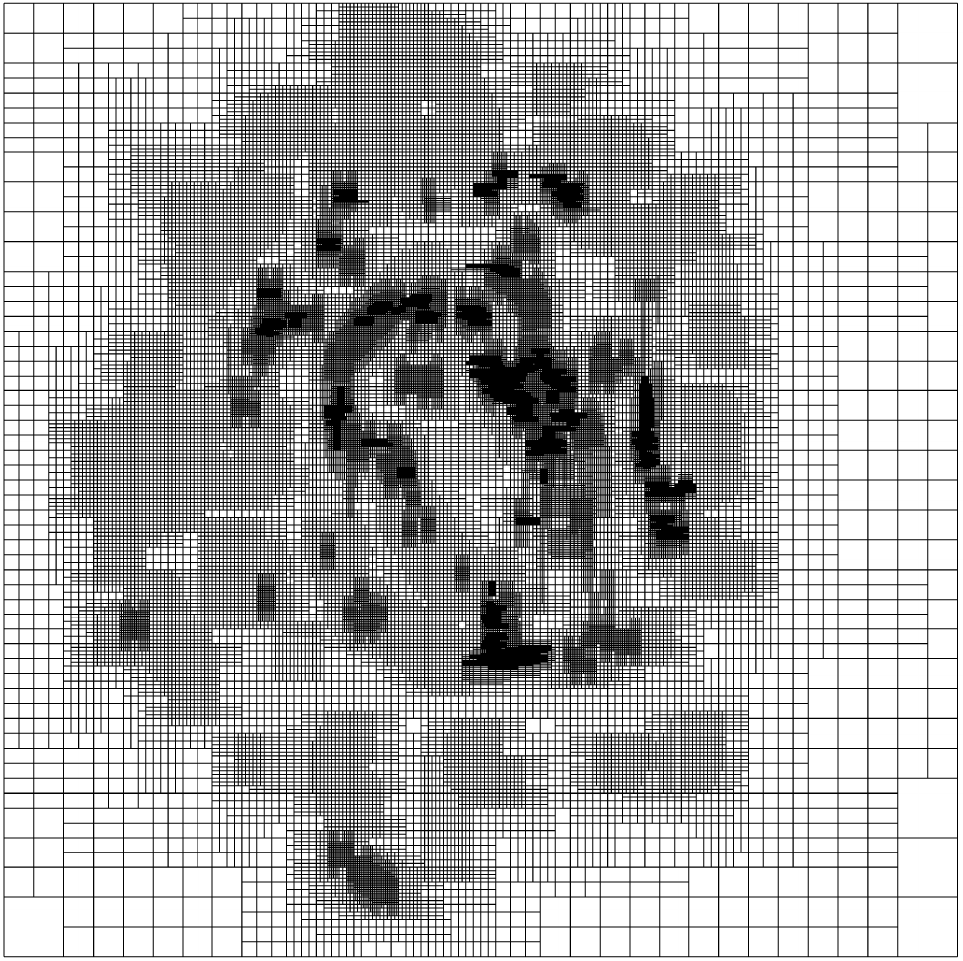}$\quad$\raisebox{1.25cm}{\includegraphics[width = 0.25\textwidth]{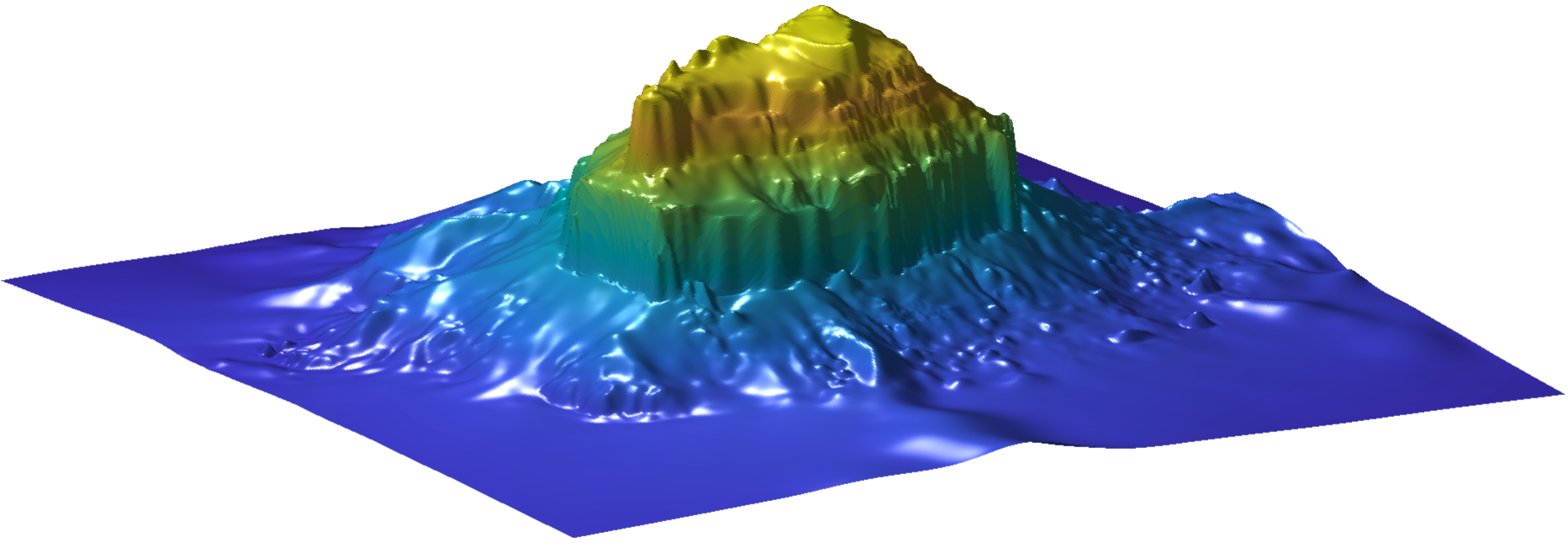}}}
\caption{QI approximation of a one million scattered point cloud of the Fajada Butte in New Mexico, USA \cite{FB}. In figure (a) the jump classification for $\bar{L} = 8$. In figure (b) we have mesh and approximation using the purely fault driven approach to localize the refinement, when setting the maximal level to $L = 10$. The number of degrees of freedom and RMSE are listed in the caption of the figure. In table (c) we report the results for the isotropic and anisotropic jump driven refinements, for different choices of $\bar{L}$ in the jump classification. In figures (d)--(e) we present meshes and approximations of the two approaches in the instance of $\bar{L} = 8$.}\label{fig:FB}
\end{figure}

As second test, we present the approximation of a dataset of the Fajada Butte in New Mexico \cite{FB}. The sampling has been rescaled in the range $[0, 0.45]$. The results and surfaces can be seen in Figure \ref{fig:FB}. As before, in figure (a) we have the jump classification for $\bar{L} = 8$ and in figure (b) we have mesh and surface of the fault driven strategy. In the caption we report the corresponding degrees of freedom and RMSE. In figure (c) we have a table schematizing the outcomes of the jump driven refinement, without and with anisotropic insertions, for different minimal jump classes. Also here, the reduction in the degrees of freedom is quite significant, despite producing approximation of the same quality of the fault guided approach, with respect to the RMSE. In figures (d) and (e) we show mesh and surfaces for $\bar{L} = 8$ in the isotropic and anisotropic, respectively, jump driven algorithms for a qualitative comparison with figure (b).
\begin{table}
\scalebox{0.85}{\renewcommand{\arraystretch}{2}\begin{tabular}{|c||c|c|c|c|c|c|c|c||c|c|c|c|c|c|c|c||}
\hline
 & \multicolumn{8}{c||}{Barringer Crater} & \multicolumn{8}{c||}{Fajada Butte} \\\hline
 &  \multicolumn{4}{c|}{$\bar{L} = 5$} & \multicolumn{4}{c||}{$\bar{L} = 7$} & \multicolumn{4}{c|}{$\bar{L} = 6$} & \multicolumn{4}{c||}{$\bar{L} = 8$} \\\hline
  & M+R & S & E & T & M+R & S & E & T & M+R & S & E & T & M+R & S & E & T\\\hline
JB & 22.9 & 10.7 & 10.0 & 14.02 & 6.82 & 4.86 & 4.73 & 5.64 & 16.9 & 7.98 & 8.57 & 11.54 & 4.15 & 2.60 & 2.73 & 3.31 \\\hline
AJB & 23.2 & 11.4 & 10.5 & 14.63 & 6.79 & 5.00 & 5.00 & 5.76 & 17.55 & 8.24 & 9.45 & 12.19 & 3.12 & 2.95 & 3.13 & 3.09\\\hline     
\end{tabular}}
\caption{Speedups of the isotropic jump based (JB) and anisotropic jump based (AJB) procedures compared to the purely fault driven (FD) method for the numerical tests on datasets sampled from the real world landscapes considered in Section \ref{sec:exm} and bidegree $(3, 3)$. More precisely, we measure the ratio of the times of execution of the fault driven over the jump driven approaches, that is, time(FD)/time(JB) and time(FD)/time(AJB) for the different steps of the adaptive cycle, i.e., Marking and Refining (M+R), Solving (S) and Evaluating/Estimating (E), as well as for the Total time (T). The times are taken for both the minimum and maximum $\bar{L}$ considered in the tests of the section. The code has been developed in Matlab R2023b and the tests have been run on a machine with a Intel i7-13700H CPU and 32GB RAM.}\label{tab:speedup}
\end{table}
\begin{table}
\scalebox{.65}{\renewcommand{\arraystretch}{2}\begin{tabular}{||c|c|c|c|c||c|c|c||c|c|c|c|c||c|c|c|c|c||}
\hline
\multicolumn{5}{||c||}{Function \eqref{eq:SS}} & \multicolumn{3}{c||}{Function \eqref{eq:BranchCuts}} & \multicolumn{5}{c||}{Barringer Crater} & \multicolumn{5}{c||}{Fajada Butte} \\\hline
& \multicolumn{2}{|c|}{isotropic} & \multicolumn{2}{c||}{anisotropic} & & \multicolumn{2}{|c||}{anisotropic}  & & \multicolumn{2}{|c|}{isotropic} & \multicolumn{2}{c||}{anisotropic}  & &\multicolumn{2}{|c|}{isotropic} & \multicolumn{2}{c||}{anisotropic}\\ \hline
$\bar{L}$        & ndofs & RMSE      & ndofs & RMSE & $\bar{L}$        & ndofs & RMSE  & $\bar{L}$        & ndofs & RMSE      & ndofs & RMSE & $\bar{L}$        & ndofs & RMSE      & ndofs & RMSE  \\ \hline
$-$ & $-$ & $-$ & $-$ & $-$  &  
1 			& 885    &  8.19e-02 & 
5                    & 6708 & 6.51e-04& 5519 & 6.54e-04 &
6                    & 12214 & 2.84e-03 & 10580 & 2.88e-03\\ \hline
$-$ & $-$ & $-$ & $-$ & $-$  & 
4			& 1127  &  8.15e-02 &
6			& 11547 & 4.02e-04 & 10119 & 4.31e-04 &
7 			& 20339 & 2.57e-03 & 17659 & 2.63e-03\\ \hline
9 	 		& 43906 & 3.17e-03 & 42089 & 3.19e-03  &
8			& 5258  & 6.65e-02   &
7 			& 22486 & 2.92e-04& 19810 & 3.07e-04 &
8			& 38171 & 2.43e-03 & 32995 & 2.47e-03\\ \hhline{||=|=|=|=|=||=|=|=||=|=|=|=|=||=|=|=|=|=||}
FD 			& 45240 & 3.17e-03 & $-$ & $-$ &
FD 			& 22068& 5.44e-02 & 
FD 			& 140689 & 2.41e-04 & $-$ & $-$ &
FD 			& 105336 & 2.34e-03 & $-$ & $-$ \\ \hline
\end{tabular}}
\caption{Numerical tests for bidegree $\m{p} = (2, 2)$. In the last row we report the outcomes for the pure fault driven (FD) method for comparison.}
\label{tab:bidegree22}
\end{table}

\section{Closure}\label{sec:closure}
In this work we have presented a procedure to estimate the jump of a function and its gradient across ordinary and gradient fault curves, respectively, when only values of such function at a scattered dataset are available. 
In the context of surface approximation, the larger the jump is, the more degrees of freedom are required to capture such behaviour. We proposed a new adaptive loop for surface reconstruction driven by the computed jump values rather than error estimates. Consequently, no intermediate approximations and evaluations are required. We have adopted LR spline spaces as adaptive spaces in order to gain the possibility for anisotropic refinements. This further feature allows to save more degrees of freedom when the fault curves are axis-aligned. 
As future developments, we consider of interest the problem of estimating the jump in case of fault intersections, both from a theoretical and a practical point of view. Furthermore, we shall extend the procedures to handle more general point clouds, such as those that cannot be parametrized as graph of a function, and data affected by noise and/or outliers, extending the preliminary study carried at the end of Section \ref{sec:funtests}.

\section*{Acknowledgments}
Supported by the National Recovery and Resilience Plan, Mission 4 Component 2 – Investment 1.4 – CN\_00000013 ``CENTRO NAZIONALE HPC, BIG DATA E QUANTUM COMPUTING'', spoke 6 (CUP B83C22002830001). The authors are members of INdAM-GNCS; CB and CG were partially supported by an INdAM-GNCS Project (CUP E53C22001930001). CB, CG, and AS also acknowledge the partial support of the Italian Ministry of University and Research (MUR) through the PRIN project COSMIC (No. 2022A79M75) and NOTES (No. P2022NC97R), funded by the European Union - Next Generation EU.

\bibliography{biblio1}

\end{document}